\documentclass{article}
\usepackage[english]{}
\usepackage[utf8]{inputenc}

\usepackage{amsmath}
\usepackage{amssymb}
\usepackage{graphicx}
\usepackage{color}
\usepackage{amsfonts}
\usepackage{setspace}
\usepackage[absolute,overlay]{textpos}
\usepackage{tikz}
\usepackage{multirow}
\usepackage{caption}
\usepackage{subcaption}
\usepackage{authblk}
\usepackage{placeins}

\usepackage[right]{lineno}
\usepackage{hyperref}
\usepackage{comment}
\title{Multicontinuum Homogenization for Poroelasticity Model}

\date{}

\author[a]{Dmitry Ammosov}
\author[a]{Mohammed Al Kobaisi\thanks{Corresponding author}}
\author[b]{Yalchin Efendiev}

\affil[a]{Chemical and Petroleum Engineering Department, Khalifa University of Science and Technology, Abu Dhabi, 127788, UAE}
\affil[b]{Department of Mathematics, Texas A\&M University, College Station, TX 77843, USA}
            
\begin{document}
\maketitle

\renewcommand{\thefootnote}{}%
\footnotetext{E-mail addresses: dmitrii.ammosov@ku.ac.ae (Dmitry Ammosov), mohammed.alkobaisi@ku.ac.ae (Mohammed Al Kobaisi), efendiev@math.tamu.edu (Yalchin Efendiev).}
\addtocounter{footnote}{-1}%

\begin{abstract}
In this paper, we derive multicontinuum poroelasticity models using the multicontinuum homogenization method. Poroelasticity models are widely used in many areas of science and engineering to describe coupled flow and mechanics processes in porous media. However, in many applications, the properties of poroelastic media possess high contrast, presenting serious computational challenges. It is well known that standard homogenization approaches often fail to give an accurate solution due to the lack of macroscopic parameters. Multicontinuum approaches allow us to consider such cases by defining several average states known as continua. In the field of poroelasticity, multiple-network models arising from the multiple porous media theory are representatives of these approaches. In this work, we extend previous findings by deriving the generalized multicontinuum poroelasticity model. We apply the recently developed multicontinuum homogenization method and provide a rigorous derivation of multicontinuum equations. For this purpose, we formulate coupled constraint cell problems in oversampled regions to consider different homogenized effects. Then, we obtain a multicontinuum expansion of the fine-scale fields and derive the multicontinuum model supposing the smoothness of macroscopic variables. We present the most general version of equations and the simplified ones based on our numerical experiments. Numerical results are presented for different heterogeneous media cases and demonstrate the high accuracy of our proposed multicontinuum models.

\end{abstract}

\section{Introduction}


Poroelasticity theory describes coupled flow and mechanics processes in porous media. It plays a significant role in many areas of science and engineering. For example, poroelasticity theory allows us to understand and predict land subsidence \cite{ferronato2006stochastic}. Modeling poroelastic effects is essential in the borehole stability analysis and preventing borehole failure \cite{ding2019borehole}. Poroelasticity theory can also explain the fluid-induced seismicity \cite{JIN2022229249}. Moreover, one can consider many biological tissues as poroelastic media, such as cardiac tissue \cite{chapelle2010poroelastic} and human brain \cite{mehrabian2011general}.

The original poroelasticity model was proposed by Biot in \cite{biot1941general} and further developed and extended in the subsequent works \cite{biot1955theory, biot1956theory, thigpen1985mechanics, zienkiewicz1984dynamic}. The basic structure of poroelasticity models consists of a coupled system of partial differential equations for mechanical displacements and pressure. The displacement equations are the elasticity equations, and the pressure equation is the parabolic equation of fluid diffusion in porous media. The coupling of the equations is governed by a mechanical body force proportional to the pressure gradient in the elasticity equations and the volumetric strain effects term in the fluid diffusion equation.


However, it should be noted that numerical simulation of many real-world applications based on poroelasticity models can present computational challenges. One of the main challenges is the multiscale nature of porous media. For example, it is known that reservoirs can possess high-contrast heterogeneous properties \cite{frias2004stochastic}. Moreover, many biological tissues also have multiscale features such as networks of fluid channels with different sizes and properties \cite{tully2011cerebral}. For accurate modeling in such cases, one must build very fine computational grids, which can significantly increase the computational cost. Given the coupled nature of the displacement and pressure equations, the cost of modeling can rise even further.


One of the popular approaches to handle multiscale problems is homogenization methods \cite{papanicolau1978asymptotic, jikov2012homogenization, lipton2006homogenization, wu2002analysis, durlofsky1991numerical}. The core idea of these methods is computing the effective properties to account for heterogeneities and building the macroscopic model. The effective properties are computed by solving cell problems at each macroscale grid block or point. Homogenization methods assume that a heterogeneous medium in a block can be replaced by a homogeneous one and often require scale separation. These methods have also been applied to elasticity and poroelasticity problems \cite{pinho2009asymptotic, vasilyeva2021machine}.


However, in complex cases with high contrast, homogenization methods may not provide an accurate solution. For such cases, macroscopic models require multiple effective properties per macroscale point. Multicontinuum approaches allow us to resolve this issue by defining multiple average states called continua and building multicontinuum models \cite{rubinvstein1948question, barenblatt1960basic, iecsan1997theory, bunoiu2019upscaling, arbogast1990derivation, chai2018efficient, bedford1972multi}. There are many works devoted to developing multicontinuum modeling. Many of these works are applied to fluid flow in porous media, with highly permeable fracture networks and low permeable porous matrix.


In terms of multicontinuum modeling of poroelastic media, multiple-network poroelastic theory (MPET) should be mentioned \cite{aifantis1980problem, wilson1982theory, bai1993multiporosity}. The MPET was developed as an extension of multiple porous media theory (double porosity, double porosity dual permeability, and multiple porosity multiple permeability) to the poroelastic case. One of the most popular applications of the multiple-network models is biological tissues such as the human brain, where one can define multiple fluid channels with different properties \cite{tully2011cerebral, guo2020multiple}. However, determining macroscopic coefficients for these models -- particularly the cross-porosity storage coefficients -- can be challenging due to the lack of experimental data to quantify them in a physiological sense \cite{guo2020multiple}.


Recently, the multicontinuum homogenization method was presented in \cite{efendiev2023multicontinuum, chung2024multicontinuum, leung2024some}. This method provides a rigorous and, at the same time, flexible methodology for deriving multicontinuum models. The main idea of the multicontinuum homogenization method is to formulate constraint cell problems in oversampled regions to account for different homogenized effects. Note that the oversampling strategy allows us to reduce boundary effects. One can obtain a homogenization-like expansion of the fine-scale solution over macroscopic variables (continua) by solving these cell problems. Then, using assumptions on the smoothness of the macroscopic variables, one can rigorously derive the multicontinuum model. The multicontinuum homogenization method has already been successfully applied for various problems with multiscale coefficients \cite{xie2025multicontinuum, ammosov2025multicontinuum, efendiev2024multicontinuum, ammosov2023multicontinuum}.


It should be noted that there are also multiscale methods for solving problems in high-contrast heterogeneous media \cite{jenny2003multi, efendiev2009multiscale, efendiev2013generalized, chung2018constraint, chaabi2024algorithmic}, including poroelasticity problems \cite{brown2016generalized, tyrylgin2020generalized, fu2020constraint}. In these methods, the fine-scale heterogeneities are captured by multiscale basis functions. These basis functions are computed by solving local problems. As a result, one can build a discrete macroscopic model on a coarse grid. Among the multiscale approaches, the multicontinuum homogenization is related to the CEM-GMsFEM \cite{chung2018constraint}, where the multiscale basis functions are also computed in oversampled regions. However, the resulting macroscopic models of the multiscale methods are not continuous and not in the form of macroscopic laws.


In this work, we extend previous findings and derive the multicontinuum poroelasticity models using the multicontinuum homogenization method. First, we formulate coupled cell problems in oversampled regions to account for different averages and gradient effects. Then, we obtain coupled multicontinuum expansions of pressure and displacement fields over macroscopic variables. We derive the corresponding multicontinuum displacement and pressure equations using these expansions. As a result, we obtain a general multicontinuum poroelasticity model for an arbitrary number of continua. Then, we present simplifications of these multicontinuum equations based on our numerical experiments. Note that our approach provides precise micro and macro relations, which overcomes the issue of measuring macroscopic parameters. We check our multicontinuum poroelasticity models by considering two-dimensional model problems with different microstructures. The numerical results show that our proposed multicontinuum approach provides high accuracy.


The paper has the following structure. Section \ref{sec:preliminaries} presents preliminaries about poroelasticity and multiple-network poroelasticity models. In Section \ref{sec:multicontinuum_homogenization}, we apply the multicontinuum homogenization method to derive multicontinuum poroelasticity equations. Section \ref{sec:numerical_results} presents numerical examples. In Section \ref{sec:conclusion}, we summarize the work.

\section{Preliminaries}\label{sec:preliminaries}

In this section, we present preliminaries on poroelasticity models. First, we outline the standard poroelasticity model and present its variational formulation. Then, we consider the multiple-network poroelasticity model, which is widely used in biological applications.


\subsection{Poroelasticity model and variational formulation}\label{sec:poroelasticity}

The poroelasticity model was initially proposed by Biot in \cite{biot1941general} and further developed in \cite{biot1955theory, biot1956theory, thigpen1985mechanics, zienkiewicz1984dynamic}. The model is described by a coupled system of partial differential equations for the displacement vector field $u$ and the pressure field $p$, which can be represented as follows
\begin{equation}
\begin{split}
&-\nabla\cdot \sigma(u)+\alpha\nabla p = f, \quad x \in \Omega,\\
&\alpha\frac{\partial \nabla \cdot u}{\partial t} + S \frac{\partial p}{\partial t} - \nabla \cdot \left(\kappa \nabla p\right) = g, \quad x \in \Omega,
\end{split}
\label{eq:system}
\end{equation}
where $\Omega \subset \mathbb{R}^2$ is a bounded computational domain, $\alpha$ is the Biot coefficient, $S = 1 / M$, $M$ is the Biot modulus, $\kappa$ is a permeability coefficient, $f$ is the body force, $g$ is the source/sink term, and $\sigma$ is the stress tensor defined as
\begin{equation*}
\sigma(u) = 2 \mu \varepsilon(u) + \lambda \nabla \cdot u \mathcal{I}, \quad 
\varepsilon(u) = \frac{1}{2} (\nabla u + (\nabla u)^T).
\end{equation*}
Here, $\lambda$ and $\mu$ are the Lam\'e parameters, $\mathcal{I}$ is the identity matrix, and $\varepsilon$ is the strain tensor. Note that, in our work, we suppose that we have multiscale porous media and the Lam\'e parameters $\lambda$ and $\mu$ and the permeability coefficient $\kappa$ possess high contrast, i.e., $\frac{\text{max}_{x \in \Omega} \lambda}{\text{min}_{x \in \Omega} \lambda} \gg 1$, $\frac{\text{max}_{x \in \Omega} \mu}{\text{min}_{x \in \Omega} \mu} \gg 1$, and $\frac{\text{max}_{x \in \Omega} \kappa}{\text{min}_{x \in \Omega} \kappa} \gg 1$.

The system \eqref{eq:system} is complemented with initial conditions $u|_{t = 0} = u_0$ and $p|_{t = 0} = p_0$ and some appropriate boundary conditions. For example, one can set the following mixed boundary conditions
\begin{equation}
\begin{split}
u = 0, \quad x \in \Gamma_D^u, \quad \sigma \cdot \nu = 0, \quad x \in \Gamma_N^u,\\
p = h, \quad x \in \Gamma_D^p, \quad -\kappa \frac{\partial p}{\partial \nu} = 0, \quad x \in \Gamma_N^p,\\
\end{split}
\end{equation}
where $\nu$ is a unit outward normal vector, and $\Gamma_D^u \cup \Gamma_N^u = \Gamma_D^p \cup \Gamma_N^p = \partial \Omega$.

Then, we can define the function spaces $W_p = \{q \in H^1(\Omega): q = h|_{\Gamma_D^p} \}$, $\hat{W}_p = \{q \in H^1(\Omega): q = 0|_{\Gamma_D^p} \}$, and $W_u = \hat{W}_u = \{v \in [H^1(\Omega)]^2: v = 0|_{\Gamma_D^u} \}$ and obtain the following variational formulation: Find $(u, p) \in W_u \times W_p$ such that
\begin{equation}
\begin{split}
&a_\Omega^u (u, v) + \bar{a}_\Omega^u (p, v) = L_\Omega^u (v),\\
&\bar{b}_\Omega^p (\frac{\partial u}{\partial t}, q) + b_\Omega^p (\frac{\partial p}{\partial t}, q) + a_\Omega^p (p, q) = L_\Omega^p (q),
\end{split}
\label{eq:weak}
\end{equation}
where $v \in \hat{W}_u$ and $q \in \hat{W}_p$ are arbitrary test functions. The bilinear and linear forms are defined as follows
\begin{equation*}
\begin{gathered}
a_\Omega^u (u, v) = \int_\Omega \sigma(u) : \varepsilon(v), \quad
\bar{a}_\Omega^u (p, v) = \int_\Omega \alpha \nabla p \cdot v, \quad
L_\Omega^u (v) = \int_\Omega f \cdot v,\\
b_\Omega^p (p, q) = \int_\Omega S p q, \quad
\bar{b}_\Omega^p (u, q) = \int_\Omega \alpha \nabla \cdot u q, \quad
a_\Omega^p (p, q) = \int_\Omega \kappa \nabla p \cdot \nabla q, \quad
L_\Omega^p (q) = \int_\Omega g q.
\end{gathered}
\end{equation*}

The subscripts of the bilinear forms indicate the respective domains of integration. Note that we can also represent the variational formulation in a coupled form
\begin{equation}
a_\Omega ((u, p), (v, q)) + b_\Omega (\frac{\partial}{\partial t}(u, p), q) = L_\Omega ((v, q)),
\label{eq:weak_2}
\end{equation}
where
\begin{equation}\label{eq:coupled_bilinear_forms}
\begin{gathered}
a_\Omega ((u, p), (v, q)) = a_\Omega^u (u, v) + \bar{a}_\Omega^u (p, v) + a_\Omega^p (p, q), \\
b_\Omega ((u, p), q) = b_\Omega^p (p, q) + \bar{b}_\Omega^p (u, q), \\
L_\Omega ((v, q)) = L_\Omega^u (v) + L_\Omega^p (q).
\end{gathered}
\end{equation}

Thus, we separate the spatial terms and the time-dependent terms into two bilinear forms.

\subsection{Multiple-network poroelasticity model}\label{sec:mpet}


It should be noted that numerical simulation using the standard poroelasticity model \eqref{eq:system} may present significant computational challenges in the case of multiscale porous media, which is considered in the present work. One of the popular examples of such multiscale media is biological tissues, where we have high-permeability arterial ($a$), arteriole/capillary ($c$), venous networks ($v$), and a low-permeability extracellular/CSF network ($e$) \cite{tully2011cerebral}. For such applications, multiple-network poroelasticity (MPET) models are widely used. One can write the MPET model in the following way \cite{guo2020multiple}
\begin{equation}
\begin{split}
&-\nabla\cdot \sigma(u) + \sum_i \alpha_i \nabla p_i = f, \quad x \in \Omega,\\
&\alpha_i \frac{\partial \nabla \cdot u}{\partial t} + S_i \frac{\partial p_i}{\partial t} + \sum_{j (i \neq j)} S_{ij} \frac{\partial p_j}{\partial t} - \nabla \cdot \left(\kappa_i \nabla p_i \right) - \sum_{j (i \neq j)} s_{ji} = g_i, \quad x \in \Omega,
\end{split}
\label{eq:mpet_system}
\end{equation}
where $i, j = a, c, e, v$ are the fluid compartment indices, $s_{ji}= \omega_{ji} (p_j - p_i)$ are the rates of fluid transfer between networks, $\omega_{ji}$ are the transfer coeﬃcients, and $S_{ij}$ are cross-porosity storage coefficients. However, the cross-porosity storage coefficients $S_{ij}$ are often neglected due to the lack of experimental data to quantify them in a physiological sense \cite{guo2020multiple}. 

In the next section, we will rigorously derive the multicontinuum poroelasticity models using the multicontinuum homogenization method. One will see that the MPET model can be considered as a particular case of the multicontinuum models. Our derivation starts at Darcy's scale, though we believe similar equations will be obtained if the microscopic model is formulated in perforated domains. We plan to study this in our future works. 

\section{Multicontinuum homogenization}\label{sec:multicontinuum_homogenization}


This section presents the derivation of the multicontinuum poroelasticity model using the multicontinuum homogenization method \cite{efendiev2023multicontinuum, chung2024multicontinuum,leung2024some}. First, we define multicontinuum expansions of the solution fields. Next, we formulate constraint cell problems in oversampled regions. Then, we derive the multicontinuum poroelasticity model. Finally, we present the general and simplified multicontinuum equations.

\subsection{Multicontinuum expansions}

\begin{figure}[hbt!]
\centering
\includegraphics[width=0.5\textwidth]{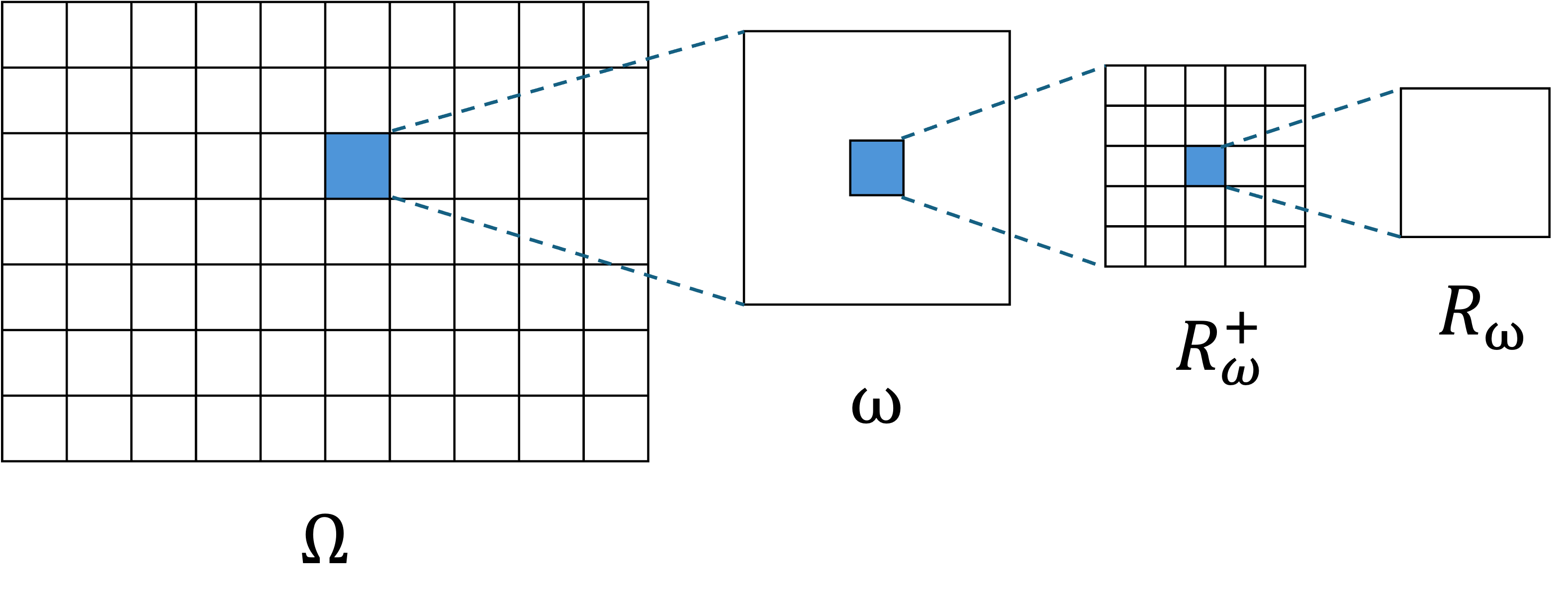}
\caption{Illustration of the domain $\Omega$, coarse block $\omega$, oversampled RVE $R_\omega^+$, and RVE $R_\omega$}
\label{fig:mh_scheme}
\end{figure}

We assume that $\Omega$ is partitioned into coarse blocks $\omega$. We suppose that there is a Representative Volume Element (RVE) $R_\omega$ inside each $\omega$. Moreover, we assume that we can define an oversampled RVE $R_\omega^+$ by extending $R_\omega$ with layers. Therefore, $R_\omega^+$ consists of RVEs $R_\omega^l$, where $l$ is a numbering, and $R_\omega^{l_0} = R_\omega$ (see Figure \ref{fig:mh_scheme}). We assume that there are $N$ components inside each RVE, which we call continua. For each continuum $j$ ($j = 1, ..., N$), we define a characteristic function $\psi_j$, which equals to $1$ in the continuum $j$ and $0$ otherwise.

Therefore, using the characteristic functions, we can define macroscopic variables, representing the average solutions in the corresponding subregions
\begin{equation}
U_{is}(x_\omega) = \frac{\int_{R_\omega} u_s \psi_i}{\int_{R_\omega} \psi_i}, \quad
P_{i}(x_\omega) = \frac{\int_{R_\omega} p \psi_i}{\int_{R_\omega} \psi_i},
\end{equation}
where $x_\omega$ is a point in $R_\omega$. Note that we assume that $U_{is}$ and $P_{i}$ are smooth functions.

The next important step of the multicontinuum homogenization is multicontinuum expansion. We expand the fine-scale solution fields over macroscopic variables in $R_\omega$
\begin{equation}
\begin{split}
&u|_{R_\omega} \approx \varphi_{is}^{uu} U_{is} + \varphi_{ism}^{uu} \nabla_m U_{is} + \varphi_{i}^{up} P_{i} + \varphi_{im}^{up} \nabla_m P_{i}, \\
&p|_{R_\omega} \approx \varphi_{is}^{pu} U_{is} + \varphi_{ism}^{pu} \nabla_m U_{is} + \varphi_{i}^{pp} P_{i} + \varphi_{im}^{pp} \nabla_m P_{i},
\end{split}
\label{eq:up_expansion}
\end{equation}
where $\varphi$ functions are the solutions of constraint cell problems. Hereafter, we assume summation over repeated indices. Note that we can also expand the test functions in a similar way
\begin{equation}
\begin{split}
&v|_{R_\omega} \approx \varphi_{jd}^{uu} V_{jd} + \varphi_{jdn}^{uu} \nabla_n V_{jd} + \varphi_{j}^{up} Q_{j} + \varphi_{jn}^{up} \nabla_n Q_{j}, \\
&q|_{R_\omega} \approx \varphi_{jd}^{pu} V_{jd} + \varphi_{jdn}^{pu} \nabla_n V_{jd} + \varphi_{j}^{pp} Q_{j} + \varphi_{jn}^{pp} \nabla_n Q_{j},
\end{split}
\label{eq:vq_expansion}
\end{equation}
where $V_{jd}$ and $Q_{j}$ are also smooth macroscopic variables.

\subsection{Cell problems}

Let us consider the cell problems for our poroelasticity problem. We suppose that $\varphi$ functions do not depend on time, and we can formulate the cell problems based only on the bilinear form $a$ defined in \eqref{eq:coupled_bilinear_forms}. Moreover, we construct the cell problems in an oversampled RVE $R_\omega^+$, i.e., we take $a_{R_{\omega^+}}$ instead of $a_{\Omega}$.

The first cell problems consider different averages in each continuum and represent the constants in the average behavior. Since we have two solution fields (displacement and pressure), we need to solve two cell problems.

\begin{itemize}
    \item Displacement cell problems for averages
    \begin{equation}
    \begin{split}
    &a_{R_\omega^+}((\varphi_{is}^{uu}, \varphi_{is}^{pu}), (v, q))
    -\sum_{j,l} \Gamma_{ijsl}^{uu} \int_{R_{\omega}^l}\psi_{j}^l v
    -\sum_{j,l} \Gamma_{ijsl}^{pu} \int_{R_{\omega}^l} \psi_{j}^l q = 0, \\
    &\int_{R_{\omega}^l} \varphi_{is}^{uu} \psi_j^l = \delta_{ij} e_s \int_{R_{\omega}^l} \psi_{j}^l,\\
    &\int_{R_{\omega}^l} \varphi_{is}^{pu} \psi_j^l = 0,
    \end{split}
    \end{equation}
    where $\Gamma$'s (in cell problems) denote the Lagrange multipliers, $e_s$ is the $s$th column of the identity matrix $I_2$.

    \item Pressure cell problems for averages
    \begin{equation}
    \begin{split}
    &a_{R_\omega^+}((\varphi_{i}^{up}, \varphi_{i}^{pp}), (v, q))
    -\sum_{j,l} \Gamma_{ijl}^{up} \int_{R_{\omega}^l}\psi_{j}^l v
    -\sum_{j,l} \Gamma_{ijl}^{pp} \int_{R_{\omega}^l} \psi_{j}^l q =0, \\
    &\int_{R_{\omega}^l} \varphi_{i}^{up} \psi_j^l = 0,\\
    &\int_{R_{\omega}^l} \varphi_{i}^{pp} \psi_j^l = \delta_{ij} \int_{R_{\omega}^l} \psi_{j}^l.
    \end{split}
    \end{equation}
\end{itemize}

The second cell problems consider gradient effects and impose constraints to represent the linear functions in the average behavior of each continuum

\begin{itemize}
    \item Displacement cell problems for gradient effects
    \begin{equation}
    \begin{split}
    &a_{R_\omega^+}((\varphi_{ism}^{uu}, \varphi_{ism}^{pu}), (v, q))
    -\sum_{j,l} \Gamma_{ijsml}^{uu} \int_{R_{\omega}^l}\psi_{j}^l v
    -\sum_{j,l} \Gamma_{ijsml}^{pu} \int_{R_{\omega}^l} \psi_{j}^l q =0, \\
    &\int_{R_{\omega}^l} \varphi_{ism}^{uu} \psi_j^l = \delta_{ij} e_s \int_{R_{\omega}^l} (x_m - c_{m})\psi_{j}^l,\\
    &\int_{R_{\omega}^l} \varphi_{ism}^{pu} \psi_j^l = 0,
    \end{split}
    \end{equation}
    where $c_m$ are some constants such that $\int_{R_\omega} (x_m - c_{m}) = 0$.

    \item Pressure cell problems for gradient effects
    \begin{equation}
    \begin{split}
    &a_{R_\omega^+}((\varphi_{im}^{up}, \varphi_{im}^{pp}), (v, q))
    -\sum_{j,l} \Gamma_{ijml}^{up} \int_{R_{\omega}^l}\psi_{j}^l v
    -\sum_{j,l} \Gamma_{ijml}^{pp} \int_{R_{\omega}^l} \psi_{j}^l q =0, \\
    &\int_{R_{\omega}^l} \varphi_{im}^{up} \psi_j^l = 0,\\
    &\int_{R_{\omega}^l} \varphi_{im}^{pp} \psi_j^l = \delta_{ij} \int_{R_{\omega}^l} (x_m - c_{m})\psi_{j}^l.
    \end{split}
    \end{equation}
\end{itemize}

Note that one can also consider higher-order cell problems, whose solutions behave as quadratic functions \cite{chung2024multicontinuum}.

\subsection{Derivation of multicontinuum equations}

Let us derive the multicontinuum poroelasticity model. First, we apply the property of the RVE to represent the whole coarse block in our variational formulation \eqref{eq:weak_2}
\begin{equation}\label{eq:appr_var_form_1}
\begin{split}
&a_{\Omega} ((u, p), (v, q)) + b_{\Omega} (\frac{\partial}{\partial t}(u, p), q) \approx \\
&\sum_{\omega} \frac{|\omega|}{|R_\omega|} \{a_{R_\omega} ((u, p), (v, q)) + b_{R_\omega} (\frac{\partial}{\partial t}(u, p), q)\} = \sum_{\omega} \frac{|\omega|}{|R_\omega|} L_{R_\omega} ((v, q)).
\end{split}
\end{equation}
Then, we need to substitute the expansions \eqref{eq:up_expansion} and \eqref{eq:vq_expansion} into this variational formulation. After that, we will obtain the following approximations of the bilinear forms 
\begin{equation}\label{eq:bilinear_forms_appr}
\begin{split}
&a_{R_\omega} ((u, p), (v, q)) \approx a_{R_\omega}^{UU} (U, V) + a_{R_\omega}^{PU} (P, V) + a_{R_\omega}^{UP} (U, Q) + a_{R_\omega}^{PP} (P, Q),\\
&b_{R_\omega} (\frac{\partial}{\partial t}(u, p), (v, q)) \approx b_{R_\omega}^{UU} (\frac{\partial U}{\partial t}, V) + b_{R_\omega}^{PU} (\frac{\partial P}{\partial t}, V) + b_{R_\omega}^{UP} (\frac{\partial U}{\partial t}, Q) + b_{R_\omega}^{PP} (\frac{\partial P}{\partial t}, Q),\\
\end{split}
\end{equation}
where $U = (U_1, ..., U_N)$ and $P = (P_1, ..., P_N)$. In these approximations, the first and second terms ($a_{R_\omega}^{UU}$, $a_{R_\omega}^{PU}$, $b_{R_\omega}^{UU}$, and $b_{R_\omega}^{PU}$) correspond to the multicontinuum displacement equations, and the third and fourth terms ($a_{R_\omega}^{UP}$, $a_{R_\omega}^{PP}$, $b_{R_\omega}^{UP}$, and $b_{R_\omega}^{PP}$) correspond to the multicontinuum pressure equations.

\subsubsection{Displacement equations}

One can see from \eqref{eq:vq_expansion} that the expansions of the test functions can be divided into $V$ and $Q$ parts. To derive the multicontinuum displacement equations, we need to consider only $V$ parts of the test functions' expansions. The expansions of the trial functions \eqref{eq:up_expansion} are taken in a complete form. First, let us consider the elliptic terms (based on the bilinear form $a_{R_\omega}$ \eqref{eq:coupled_bilinear_forms}) of the variational formulation. We obtain the following approximation terms
\begin{equation}\label{eq:displacement_appr_elliptic_terms}
\begin{split}
a_{R_\omega}^{UU} (U, V) =
&a_{R_\omega} ((\varphi_{ism}^{uu}, \varphi_{ism}^{pu}), (\varphi_{jdn}^{uu}, \varphi_{jdn}^{pu})) \nabla_m U_{is} \nabla_n V_{jd} +
a_{R_\omega} ((\varphi_{ism}^{uu}, \varphi_{ism}^{pu}), (\varphi_{jd}^{uu}, \varphi_{jd}^{pu})) \nabla_m U_{is} V_{jd} +\\
&a_{R_\omega} ((\varphi_{is}^{uu}, \varphi_{is}^{pu}), (\varphi_{jdn}^{uu}, \varphi_{jdn}^{pu})) U_{is} \nabla_n V_{jd} +
a_{R_\omega} ((\varphi_{is}^{uu}, \varphi_{is}^{pu}), (\varphi_{jd}^{uu}, \varphi_{jd}^{pu})) U_{is} V_{jd} =\\
&|R_\omega|\{A_{jdnism}^{uu} \nabla_m U_{is} \nabla_n V_{jd} + 
B_{jdism}^{uu} \nabla_m U_{is} V_{jd} +
\bar{B}_{jdnis}^{uu} U_{is} \nabla_n V_{jd} +
C_{jdis}^{uu} U_{is} V_{jd}\},\\
a_{R_\omega}^{PU} (P, V) =
&a_{R_\omega} ((\varphi_{im}^{up}, \varphi_{im}^{pp}), (\varphi_{jdn}^{uu}, \varphi_{jdn}^{pu})) \nabla_m P_{i} \nabla_n V_{jd} +
a_{R_\omega} ((\varphi_{im}^{up}, \varphi_{im}^{pp}), (\varphi_{jd}^{uu}, \varphi_{jd}^{pu})) \nabla_m P_{i} V_{jd} +\\
&a_{R_\omega} ((\varphi_{i}^{up}, \varphi_{i}^{pp}), (\varphi_{jdn}^{uu}, \varphi_{jdn}^{pu})) P_{i} \nabla_n V_{jd} +
a_{R_\omega} ((\varphi_{i}^{up}, \varphi_{i}^{pp}), (\varphi_{jd}^{uu}, \varphi_{jd}^{pu})) P_{i} V_{jd} =\\
&|R_\omega|\{A_{jdnim}^{pu} \nabla_m P_{i} \nabla_n V_{jd} + 
B_{jdim}^{pu} \nabla_m P_{i} V_{jd} +
\bar{B}_{jdni}^{pu} P_{i} \nabla_n V_{jd} +
C_{jdi}^{pu} P_{i} V_{jd}\}.
\end{split}
\end{equation}

Here, we have used the assumption about the smoothness of macroscopic variables $U_{is}$, $P_i$, and $V_{jd}$ and taken them out of the integrals over $R_\omega$. The effective properties can be represented as follows

\begin{equation*}
\begin{gathered}
A_{jdnism}^{uu} = \frac{1}{|R_\omega|} a_{R_\omega} ((\varphi_{ism}^{uu}, \varphi_{ism}^{pu}), (\varphi_{jdn}^{uu}, \varphi_{jdn}^{pu})), \quad
B_{jdism}^{uu} = \frac{1}{|R_\omega|} a_{R_\omega} ((\varphi_{ism}^{uu}, \varphi_{ism}^{pu}), (\varphi_{jd}^{uu}, \varphi_{jd}^{pu})), \\
\bar{B}_{jdnis}^{uu} = \frac{1}{|R_\omega|} a_{R_\omega} ((\varphi_{is}^{uu}, \varphi_{is}^{pu}), (\varphi_{jdn}^{uu}, \varphi_{jdn}^{pu})), \quad
C_{jdis}^{uu} = \frac{1}{|R_\omega|} a_{R_\omega} ((\varphi_{is}^{uu}, \varphi_{is}^{pu}), (\varphi_{jd}^{uu}, \varphi_{jd}^{pu})), \\
A_{jdnim}^{pu} = \frac{1}{|R_\omega|} a_{R_\omega} ((\varphi_{im}^{up}, \varphi_{im}^{pp}), (\varphi_{jdn}^{uu}, \varphi_{jdn}^{pu})), \quad
B_{jdim}^{pu} = \frac{1}{|R_\omega|} a_{R_\omega} ((\varphi_{im}^{up}, \varphi_{im}^{pp}), (\varphi_{jd}^{uu}, \varphi_{jd}^{pu})), \\
\bar{B}_{jdni}^{pu} = \frac{1}{|R_\omega|} a_{R_\omega} ((\varphi_{i}^{up}, \varphi_{i}^{pp}), (\varphi_{jdn}^{uu}, \varphi_{jdn}^{pu})), \quad
C_{jdi}^{pu} = \frac{1}{|R_\omega|} a_{R_\omega} ((\varphi_{i}^{up}, \varphi_{i}^{pp}), (\varphi_{jd}^{uu}, \varphi_{jd}^{pu})).
\end{gathered}
\end{equation*}

One can also derive the time derivative terms using the bilinear form $b_{R_\omega}$ (defined in \eqref{eq:coupled_bilinear_forms}) in a similar way
\begin{equation}\label{eq:displacement_appr_time_terms}
\begin{split}
b_{R_\omega}^{UU} (\frac{\partial U}{\partial t}, V) =
&|R_\omega|\{D_{jdnism}^{uu} \frac{\partial (\nabla_m U_{is})}{\partial t} \nabla_n V_{jd} + 
G_{jdism}^{uu} \frac{\partial (\nabla_m U_{is})}{\partial t} V_{jd} +
\bar{G}_{jdnis}^{uu} \frac{\partial U_{is}}{\partial t} \nabla_n V_{jd} +
H_{jdis}^{uu} \frac{\partial U_{is}}{\partial t} V_{jd}\},\\
b_{R_\omega}^{PU} (\frac{\partial P}{\partial t}, V) =
&|R_\omega|\{D_{jdnim}^{pu} \frac{\partial (\nabla_m P_{i})}{\partial t} \nabla_n V_{jd} + 
G_{jdim}^{pu} \frac{\partial (\nabla_m P_{i})}{\partial t} V_{jd} +
\bar{G}_{jdni}^{pu} \frac{\partial P_{i}}{\partial t} \nabla_n V_{jd} +
H_{jdi}^{pu} \frac{\partial P_{i}}{\partial t} V_{jd}\},
\end{split}
\end{equation}
where the effective properties are following
\begin{equation*}
\begin{gathered}
D_{jdnism}^{uu} = \frac{1}{|R_\omega|} b_{R_\omega} ((\varphi_{ism}^{uu}, \varphi_{ism}^{pu}), \varphi_{jdn}^{pu}), \quad
G_{jdism}^{uu} = \frac{1}{|R_\omega|} b_{R_\omega} ((\varphi_{ism}^{uu}, \varphi_{ism}^{pu}), \varphi_{jd}^{pu}), \\
\bar{G}_{jdnis}^{uu} = \frac{1}{|R_\omega|} b_{R_\omega} ((\varphi_{is}^{uu}, \varphi_{is}^{pu}), \varphi_{jdn}^{pu}), \quad
H_{jdis}^{uu} = \frac{1}{|R_\omega|} b_{R_\omega} ((\varphi_{is}^{uu}, \varphi_{is}^{pu}), \varphi_{jd}^{pu}), \\
D_{jdnim}^{pu} = \frac{1}{|R_\omega|} b_{R_\omega} ((\varphi_{im}^{up}, \varphi_{im}^{pp}), \varphi_{jdn}^{pu}), \quad
G_{jdim}^{pu} = \frac{1}{|R_\omega|} b_{R_\omega} ((\varphi_{im}^{up}, \varphi_{im}^{pp}), \varphi_{jd}^{pu}), \\
\bar{G}_{jdni}^{pu} = \frac{1}{|R_\omega|} b_{R_\omega} ((\varphi_{i}^{up}, \varphi_{i}^{pp}), \varphi_{jdn}^{pu}), \quad
H_{jdi}^{pu} = \frac{1}{|R_\omega|} b_{R_\omega} ((\varphi_{i}^{up}, \varphi_{i}^{pp}), \varphi_{jd}^{pu}).
\end{gathered}
\end{equation*}

Then, we sum the terms $a_{R_\omega}^{UU}$, $a_{R_\omega}^{PU}$ from \eqref{eq:displacement_appr_elliptic_terms} and $b_{R_\omega}^{UU}$, $b_{R_\omega}^{PU}$ from \eqref{eq:displacement_appr_time_terms} over all coarse blocks using the RVE property and obtain the weak formulation of the multicontinuum displacement equations
\begin{equation}\label{eq:mc_displacement_weak}
\begin{split}
\int_\Omega \{
&A_{jdnism}^{uu} \nabla_m U_{is} \nabla_n V_{jd} + 
B_{jdism}^{uu} \nabla_m U_{is} V_{jd} +
\bar{B}_{jdnis}^{uu} U_{is} \nabla_n V_{jd} +
C_{jdis}^{uu} U_{is} V_{jd} +\\
&A_{jdnim}^{pu} \nabla_m P_{i} \nabla_n V_{jd} + 
B_{jdim}^{pu} \nabla_m P_{i} V_{jd} +
\bar{B}_{jdni}^{pu} P_{i} \nabla_n V_{jd} +
C_{jdi}^{pu} P_{i} V_{jd} +\\
&D_{jdnism}^{uu} \frac{\partial (\nabla_m U_{is})}{\partial t} \nabla_n V_{jd} + 
G_{jdism}^{uu} \frac{\partial (\nabla_m U_{is})}{\partial t} V_{jd} +
\bar{G}_{jdnis}^{uu} \frac{\partial U_{is}}{\partial t} \nabla_n V_{jd} +
H_{jdis}^{uu} \frac{\partial U_{is}}{\partial t} V_{jd} +\\
&D_{jdnim}^{pu} \frac{\partial (\nabla_m P_{i})}{\partial t} \nabla_n V_{jd} + 
G_{jdim}^{pu} \frac{\partial (\nabla_m P_{i})}{\partial t} V_{jd} +
\bar{G}_{jdni}^{pu} \frac{\partial P_{i}}{\partial t} \nabla_n V_{jd} +
H_{jdi}^{pu} \frac{\partial P_{i}}{\partial t} V_{jd}\} = \int_\Omega F_{jd}^{u} V_{jd},
\end{split}
\end{equation}
where the multicontinuum body forces are defined as $F_{jd}^{u} = \frac{1}{|R_\omega|} L_{R_\omega} ((\varphi_{jd}^{uu}, \varphi_{jd}^{pu}))$.

Note that we can represent \eqref{eq:mc_displacement_weak} in the following strong form
\begin{equation}
\begin{split}
&-\nabla_n ( A_{jdnism}^{uu} \nabla_m U_{is} ) + 
B_{jdism}^{uu} \nabla_m U_{is} -
\nabla_n ( \bar{B}_{jdnis}^{uu} U_{is} ) +
C_{jdis}^{uu} U_{is} -\\
&\nabla_n ( A_{jdnim}^{pu} \nabla_m P_{i} ) + 
B_{jdim}^{pu} \nabla_m P_{i} -
\nabla_n ( \bar{B}_{jdni}^{pu} P_{i} ) +
C_{jdi}^{pu} P_{i} -\\
&\frac{\partial}{\partial t} ( \nabla_n ( D_{jdnism}^{uu} \nabla_m U_{is} ) ) + 
\frac{\partial}{\partial t} ( G_{jdism}^{uu} \nabla_m U_{is} ) -
\frac{\partial}{\partial t} ( \nabla_n ( \bar{G}_{jdnis}^{uu} U_{is} ) ) +
H_{jdis}^{uu} \frac{\partial U_{is}}{\partial t} -\\
&\frac{\partial}{\partial t} ( \nabla_n ( D_{jdnim}^{pu} \nabla_m P_{i} ) ) + 
\frac{\partial}{\partial t} ( G_{jdim}^{pu} \nabla_m P_{i} ) -
\frac{\partial}{\partial t} ( \nabla_n ( \bar{G}_{jdni}^{pu} P_{i} ) ) +
H_{jdi}^{pu} \frac{\partial P_{i}}{\partial t} = F_{jd}^{u}.
\end{split}
\end{equation}

Our numerical results show that we can neglect some terms and obtain simplified equations. For example, the time derivative terms can be neglected.

\subsubsection{Pressure equations} 

The derivation of the pressure equations are similar to the displacement equations. We consider $Q$ parts of the test functions' expansions \eqref{eq:vq_expansion} and take the trial functions' expansions in a complete form. Let us consider the elliptic term based on the bilinear form $a_{R_\omega}$. We obtain the following approximation terms

\begin{equation}\label{eq:pressure_appr_elliptic_terms}
\begin{split}
a_{R_\omega}^{UP} (U, Q)=
&a_{R_\omega} ((\varphi_{ism}^{uu}, \varphi_{ism}^{pu}), (\varphi_{jn}^{up}, \varphi_{jn}^{pp})) \nabla_m U_{is} \nabla_n Q_{j} +
a_{R_\omega} ((\varphi_{ism}^{uu}, \varphi_{ism}^{pu}), (\varphi_{j}^{up}, \varphi_{j}^{pp})) \nabla_m U_{is} Q_{j} +\\
&a_{R_\omega} ((\varphi_{is}^{uu}, \varphi_{is}^{pu}), (\varphi_{jn}^{up}, \varphi_{jn}^{pp})) U_{is} \nabla_n Q_{j} +
a_{R_\omega} ((\varphi_{is}^{uu}, \varphi_{is}^{pu}), (\varphi_{j}^{up}, \varphi_{j}^{pp})) U_{is} Q_{j} =\\
&|R_\omega|\{A_{jnism}^{up} \nabla_m U_{is} \nabla_n Q_{j} + 
B_{jism}^{up} \nabla_m U_{is} Q_{j} +
\bar{B}_{jnis}^{up} U_{is} \nabla_n Q_{j} +
C_{jis}^{up} U_{is} Q_{j}\},\\
a_{R_\omega}^{PP} (P, Q)=
&a_{R_\omega} ((\varphi_{im}^{up}, \varphi_{im}^{pp}), (\varphi_{jn}^{up}, \varphi_{jn}^{pp})) \nabla_m P_{i} \nabla_n Q_{j} +
a_{R_\omega} ((\varphi_{im}^{up}, \varphi_{im}^{pp}), (\varphi_{j}^{up}, \varphi_{j}^{pp})) \nabla_m P_{i} Q_{j} +\\
&a_{R_\omega} ((\varphi_{i}^{up}, \varphi_{i}^{pp}), (\varphi_{jn}^{up}, \varphi_{jn}^{pp})) P_{i} \nabla_n Q_{j} +
a_{R_\omega} ((\varphi_{i}^{up}, \varphi_{i}^{pp}), (\varphi_{j}^{up}, \varphi_{j}^{pp})) P_{i} Q_{j} =\\
&|R_\omega|\{A_{jnim}^{pp} \nabla_m P_{i} \nabla_n Q_{j} + 
B_{jim}^{pp} \nabla_m P_{i} Q_{j} +
\bar{B}_{jni}^{pp} P_{i} \nabla_n Q_{j} +
C_{ji}^{pp} P_{i} Q_{j}\},
\end{split}
\end{equation}
where the effective properties are following
\begin{equation*}
\begin{gathered}
A_{jnism}^{up} = \frac{1}{|R_\omega|} a_{R_\omega} ((\varphi_{ism}^{uu}, \varphi_{ism}^{pu}), (\varphi_{jn}^{up}, \varphi_{jn}^{pp})), \quad
B_{jism}^{up} = \frac{1}{|R_\omega|} a_{R_\omega} ((\varphi_{ism}^{uu}, \varphi_{ism}^{pu}), (\varphi_{j}^{up}, \varphi_{j}^{pp})), \\
\bar{B}_{jnis}^{up} = \frac{1}{|R_\omega|} a_{R_\omega} ((\varphi_{is}^{uu}, \varphi_{is}^{pu}), (\varphi_{jn}^{up}, \varphi_{jn}^{pp})), \quad
C_{jis}^{up} = \frac{1}{|R_\omega|} a_{R_\omega} ((\varphi_{is}^{uu}, \varphi_{is}^{pu}), (\varphi_{j}^{up}, \varphi_{j}^{pp})), \\
A_{jnim}^{pp} = \frac{1}{|R_\omega|} a_{R_\omega} ((\varphi_{im}^{up}, \varphi_{im}^{pp}), (\varphi_{jn}^{up}, \varphi_{jn}^{pp})), \quad
B_{jim}^{pp} = \frac{1}{|R_\omega|} a_{R_\omega} ((\varphi_{im}^{up}, \varphi_{im}^{pp}), (\varphi_{j}^{up}, \varphi_{j}^{pp})), \\
\bar{B}_{jni}^{pp} = \frac{1}{|R_\omega|} a_{R_\omega} ((\varphi_{i}^{up}, \varphi_{i}^{pp}), (\varphi_{jn}^{up}, \varphi_{jn}^{pp})), \quad
C_{ji}^{pp} = \frac{1}{|R_\omega|} a_{R_\omega} ((\varphi_{i}^{up}, \varphi_{i}^{pp}), (\varphi_{j}^{up}, \varphi_{j}^{pp})).
\end{gathered}
\end{equation*}

We can also derive the time derivative terms based on the bilinear form $b_{R_\omega}$ as follows
\begin{equation}\label{eq:pressure_appr_time_terms}
\begin{split}
b_{R_\omega}^{UP} (\frac{\partial U}{\partial t}, Q) =
&|R_\omega|\{D_{jnism}^{up} \frac{\partial (\nabla_m U_{is})}{\partial t} \nabla_n Q_{j} + 
G_{jism}^{up} \frac{\partial (\nabla_m U_{is})}{\partial t} Q_{j} +
\bar{G}_{jnis}^{up} \frac{\partial U_{is}}{\partial t} \nabla_n Q_{j} +
H_{jis}^{up} \frac{\partial U_{is}}{\partial t} Q_{j}\},\\
b_{R_\omega}^{PP} (\frac{\partial P}{\partial t}, Q) =
&|R_\omega|\{D_{jnim}^{pp} \frac{\partial (\nabla_m P_{i})}{\partial t} \nabla_n Q_{j} + 
G_{jim}^{pp} \frac{\partial (\nabla_m P_{i})}{\partial t} Q_{j} +
\bar{G}_{jni}^{pp} \frac{\partial P_{i}}{\partial t} \nabla_n Q_{j} +
H_{ji}^{pp} \frac{\partial P_{i}}{\partial t} Q_{j}\},
\end{split}
\end{equation}
where
\begin{equation*}
\begin{gathered}
D_{jnism}^{up} = \frac{1}{|R_\omega|} b_{R_\omega} \left((\varphi_{ism}^{uu}, \varphi_{ism}^{pu}), \varphi_{jn}^{pp}\right), \quad
G_{jism}^{up} = \frac{1}{|R_\omega|} b_{R_\omega} \left((\varphi_{ism}^{uu}, \varphi_{ism}^{pu}), \varphi_{j}^{pp}\right), \\
\bar{G}_{jnis}^{up} = \frac{1}{|R_\omega|} b_{R_\omega} \left((\varphi_{is}^{uu}, \varphi_{is}^{pu}), \varphi_{jn}^{pp}\right), \quad
H_{jis}^{up} = \frac{1}{|R_\omega|} b_{R_\omega} \left((\varphi_{is}^{uu}, \varphi_{is}^{pu}), \varphi_{j}^{pp}\right), \\
D_{jnim}^{pp} = \frac{1}{|R_\omega|} b_{R_\omega} \left((\varphi_{im}^{up}, \varphi_{im}^{pp}), \varphi_{jn}^{pp}\right), \quad
G_{jim}^{pp} = \frac{1}{|R_\omega|} b_{R_\omega} \left((\varphi_{im}^{up}, \varphi_{im}^{pp}), \varphi_{j}^{pp}\right), \\
\bar{G}_{jni}^{pp} = \frac{1}{|R_\omega|} b_{R_\omega} \left((\varphi_{i}^{up}, \varphi_{i}^{pp}), \varphi_{jn}^{pp}\right), \quad
H_{ji}^{pp} = \frac{1}{|R_\omega|} b_{R_\omega} \left((\varphi_{i}^{up}, \varphi_{i}^{pp}), \varphi_{j}^{pp}\right).
\end{gathered}
\end{equation*}

Again, we sum the terms $a_{R_\omega}^{UP}$, $a_{R_\omega}^{PP}$ from \eqref{eq:pressure_appr_elliptic_terms} and $b_{R_\omega}^{UP}$, $b_{R_\omega}^{PP}$ from \eqref{eq:pressure_appr_time_terms} over all coarse blocks and obtain the following weak formulation of the multicontinuum pressure equations
\begin{equation}\label{eq:mc_pressure_weak}
\begin{split}
\int_\Omega \{
&A_{jnism}^{up} \nabla_m U_{is} \nabla_n Q_{j} + 
B_{jism}^{up} \nabla_m U_{is} Q_{j} +
\bar{B}_{jnis}^{up} U_{is} \nabla_n Q_{j} +
C_{jis}^{up} U_{is} Q_{j} +\\
&A_{jnim}^{pp} \nabla_m P_{i} \nabla_n Q_{j} + 
B_{jim}^{pp} \nabla_m P_{i} Q_{j} +
\bar{B}_{jni}^{pp} P_{i} \nabla_n Q_{j} +
C_{ji}^{pp} P_{i} Q_{j} +\\
&D_{jnism}^{up} \frac{\partial (\nabla_m U_{is})}{\partial t} \nabla_n Q_{j} + 
G_{jism}^{up} \frac{\partial (\nabla_m U_{is})}{\partial t} Q_{j} +
\bar{G}_{jnis}^{up} \frac{\partial U_{is}}{\partial t} \nabla_n Q_{j} +
H_{jis}^{up} \frac{\partial U_{is}}{\partial t} Q_{j} +\\
&D_{jnim}^{pp} \frac{\partial (\nabla_m P_{i})}{\partial t} \nabla_n Q_{j} + 
G_{jim}^{pp} \frac{\partial (\nabla_m P_{i})}{\partial t} Q_{j} +
\bar{G}_{jni}^{pp} \frac{\partial P_{i}}{\partial t} \nabla_n Q_{j} +
H_{ji}^{pp} \frac{\partial P_{i}}{\partial t} Q_{j} \} = \int_\Omega F_{j}^{p} Q_{j},
\end{split}
\end{equation}
where the source terms are defined as $F_{j}^{p} = \frac{1}{|R_\omega|} L_{R_\omega} ((\varphi_{j}^{up}, \varphi_{j}^{pp}))$.

Therefore, we have the following general multicontinuum equations for the pressure field
\begin{equation}
\begin{split}
&-\nabla_n ( A_{jnism}^{up} \nabla_m U_{is} ) + 
B_{jism}^{up} \nabla_m U_{is} -
\nabla_n ( \bar{B}_{jnis}^{up} U_{is} ) +
C_{jis}^{up} U_{is} -\\
&\nabla_n ( A_{jnim}^{pp} \nabla_m P_{i} ) + 
B_{jim}^{pp} \nabla_m P_{i} -
\nabla_n ( \bar{B}_{jni}^{pp} P_{i} ) +
C_{ji}^{pp} P_{i} -\\
&\frac{\partial}{\partial t} ( \nabla_n ( D_{jnism}^{up} \nabla_m U_{is} ) ) + 
\frac{\partial}{\partial t} ( G_{jism}^{up} \nabla_m U_{is} ) -
\frac{\partial}{\partial t} ( \nabla_n ( \bar{G}_{jnis}^{up} U_{is} ) ) +
H_{jis}^{up} \frac{\partial U_{is}}{\partial t} -\\
&\frac{\partial}{\partial t} ( \nabla_n ( D_{jnim}^{pp} \nabla_m P_{i} ) ) + 
\frac{\partial}{\partial t} ( G_{jim}^{pp} \nabla_m P_{i} ) -
\frac{\partial}{\partial t} ( \nabla_n ( \bar{G}_{jni}^{pp} P_{i} ) ) +
H_{ji}^{pp} \frac{\partial P_{i}}{\partial t} = F_{j}^{p}.
\end{split}
\end{equation}

One can neglect some terms and obtain simplified equations. In the next subsection, we will present the full multicontinuum poroelasticity model and its simplified versions.

\subsection{Multicontinuum poroelasticity models}

Finally, let us summarize the obtained multicontinuum models. We have derived general multicontinuum displacement and pressure equations for an arbitrary number of continua. Based on our numerical experiments, we have presented simplified versions of these equations.
\begin{itemize}
    \item Full model
    \begin{itemize}
        \item Multicontinuum displacement equations
        \begin{equation}\label{eq:full_mc_model_u_eq}
        \begin{split}
        &-\nabla_n ( A_{jdnism}^{uu} \nabla_m U_{is} ) + 
        B_{jdism}^{uu} \nabla_m U_{is} -
        \nabla_n ( \bar{B}_{jdnis}^{uu} U_{is} ) +
        C_{jdis}^{uu} U_{is} -\\
        &\nabla_n ( A_{jdnim}^{pu} \nabla_m P_{i} ) + 
        B_{jdim}^{pu} \nabla_m P_{i} -
        \nabla_n ( \bar{B}_{jdni}^{pu} P_{i} ) +
        C_{jdi}^{pu} P_{i} -\\
        &\frac{\partial}{\partial t} ( \nabla_n ( D_{jdnism}^{uu} \nabla_m U_{is} ) ) + 
        \frac{\partial}{\partial t} ( G_{jdism}^{uu} \nabla_m U_{is} ) -
        \frac{\partial}{\partial t} ( \nabla_n ( \bar{G}_{jdnis}^{uu} U_{is} ) ) +
        H_{jdis}^{uu} \frac{\partial U_{is}}{\partial t} -\\
        &\frac{\partial}{\partial t} ( \nabla_n ( D_{jdnim}^{pu} \nabla_m P_{i} ) ) + 
        \frac{\partial}{\partial t} ( G_{jdim}^{pu} \nabla_m P_{i} ) -
        \frac{\partial}{\partial t} ( \nabla_n ( \bar{G}_{jdni}^{pu} P_{i} ) ) +
        H_{jdi}^{pu} \frac{\partial P_{i}}{\partial t} = F_{jd}^{u}.
        \end{split}
        \end{equation}

        \item Multicontinuum pressure equations
        \begin{equation}\label{eq:full_mc_model_p_eq}
        \begin{split}
        &-\nabla_n ( A_{jnism}^{up} \nabla_m U_{is} ) + 
        B_{jism}^{up} \nabla_m U_{is} -
        \nabla_n ( \bar{B}_{jnis}^{up} U_{is} ) +
        C_{jis}^{up} U_{is} -\\
        &\nabla_n ( A_{jnim}^{pp} \nabla_m P_{i} ) + 
        B_{jim}^{pp} \nabla_m P_{i} -
        \nabla_n ( \bar{B}_{jni}^{pp} P_{i} ) +
        C_{ji}^{pp} P_{i} -\\
        &\frac{\partial}{\partial t} ( \nabla_n ( D_{jnism}^{up} \nabla_m U_{is} ) ) + 
        \frac{\partial}{\partial t} ( G_{jism}^{up} \nabla_m U_{is} ) -
        \frac{\partial}{\partial t} ( \nabla_n ( \bar{G}_{jnis}^{up} U_{is} ) ) +
        H_{jis}^{up} \frac{\partial U_{is}}{\partial t} -\\
        &\frac{\partial}{\partial t} ( \nabla_n ( D_{jnim}^{pp} \nabla_m P_{i} ) ) + 
        \frac{\partial}{\partial t} ( G_{jim}^{pp} \nabla_m P_{i} ) -
        \frac{\partial}{\partial t} ( \nabla_n ( \bar{G}_{jni}^{pp} P_{i} ) ) +
        H_{ji}^{pp} \frac{\partial P_{i}}{\partial t} = F_{j}^{p}.
        \end{split}
        \end{equation}
    \end{itemize}

    \item Simplified model 1
    \begin{itemize}
        \item Multicontinuum displacement equations
        \begin{equation}\label{eq:simplified_mc_model_1_u_eq}
        \begin{split}
        &-\nabla_n ( A_{jdnism}^{uu} \nabla_m U_{is} ) + 
        B_{jdism}^{uu} \nabla_m U_{is} -
        \nabla_n ( \bar{B}_{jdnis}^{uu} U_{is} ) +
        C_{jdis}^{uu} U_{is} +\\
        &B_{jdim}^{pu} \nabla_m P_{i} -
        \nabla_n ( \bar{B}_{jdni}^{pu} P_{i} ) = F_{jd}^{u}.
        \end{split}
        \end{equation}

        \item Multicontinuum pressure equations
        \begin{equation}\label{eq:simplified_mc_model_1_p_eq}
        \begin{split}
        &B_{jism}^{up} \nabla_m U_{is} -
        \nabla_n ( \bar{B}_{jnis}^{up} U_{is} ) -
        \nabla_n ( A_{jnim}^{pp} \nabla_m P_{i} ) +
        C_{ji}^{pp} P_{i} +\\
        &\frac{\partial}{\partial t} ( G_{jism}^{up} \nabla_m U_{is} ) -
        \frac{\partial}{\partial t} ( \nabla_n ( \bar{G}_{jnis}^{up} U_{is} ) ) -
        \frac{\partial}{\partial t} ( \nabla_n ( D_{jnim}^{pp} \nabla_m P_{i} ) ) +
        H_{ji}^{pp} \frac{\partial P_{i}}{\partial t} = F_{j}^{p}.
        \end{split}
        \end{equation}
    \end{itemize}

    \item Simplified model 2
    \begin{itemize}
        \item Multicontinuum displacement equations
        \begin{equation}\label{eq:simplified_mc_model_2_u_eq}
        \begin{split}
        &-\nabla_n ( A_{jdnism}^{uu} \nabla_m U_{is} ) +
        C_{jdis}^{uu} U_{is} +
        B_{jdim}^{pu} \nabla_m P_{i} = F_{jd}^{u}.
        \end{split}
        \end{equation}

        \item Multicontinuum pressure equations
        \begin{equation}\label{eq:simplified_mc_model_2_p_eq}
        \begin{split}
        -\nabla_n ( A_{jnim}^{pp} \nabla_m P_{i} ) +
        C_{ji}^{pp} P_{i} +
        \frac{\partial}{\partial t} ( G_{jism}^{up} \nabla_m U_{is} ) +
        H_{ji}^{pp} \frac{\partial P_{i}}{\partial t} = F_{j}^{p}.
        \end{split}
        \end{equation}
    \end{itemize}
\end{itemize}

One can see that the MPET model \eqref{eq:mpet_system} is a particular case of our second simplified model with some neglected terms (such as additional diffusion terms) and single continuum for displacement. Note that the first simplified model contains new coupling terms between pressures and displacements such as the first two terms in the pressure equations. In general, the influence of these terms depends on the microstructure and material parameters.

\section{Numerical examples}\label{sec:numerical_results}


This section presents representative numerical examples of our proposed multicontinuum poroelasticity models. We consider three model problems with different heterogeneous microstructures and boundary conditions. We set $\Omega = \Omega_1 \cup \Omega_2 = [0, 1]\times [0, 1]$ as the computational domain, where $\Omega_1$ is the region of the first continuum, and $\Omega_2$ is the region of the second continuum. We assume that the Lam\'e parameters $\lambda$ and $\mu$ and the permeability coefficient $\kappa$ possess high contrast.

For the first two model problems, we use a uniform fine grid with $160 801$ vertices and $320 000$ triangular cells. We construct an unstructured fine grid with $85 012$ vertices and $169 019$ triangular cells in the third model problem. In all the numerical examples, we consider two uniform coarse grids: (a) $10 \times 10$ (with $101$ vertices and $100$ rectangular cells) and (b) $20 \times 20$ (with $401$ vertices and $400$ rectangular cells). 

For simplicity, we take each coarse block as an RVE itself. We construct the oversampled RVE as an extension of the RVE by $l$ layers. For this purpose, we set $l = 5$ layers for the coarse grid $10 \times 10$ and $l = 6$ layers for the coarse grid $20 \times 20$.
 
We solve the model problems using our proposed multicontinuum poroelasticity models (full, simplified 1, and simplified 2) on the coarse grids. As reference solutions, we consider the finite-element solutions of the fine-scale poroelasticity model using linear basis functions on the fine grid. To compare the homogenized and reference solutions, we use the following relative $L_2$ errors for pressure and displacements
\begin{equation*}
e^{(i)}_p = \sqrt{\frac{\sum_K | \frac{1}{|K|} \int_K P_i \,dx - \frac{1}{|K \cap \Omega_i|} \int_{K \cap \Omega_i} p \,dx|^2}{\sum_K | \frac{1}{|K \cap \Omega_i|} \int_{K \cap \Omega_i} p \,dx|^2}}, \quad
e^{(i)}_u = \sqrt{\frac{\sum_K | \frac{1}{|K|} \int_K U_i \,dx - \frac{1}{|K \cap \Omega_i|} \int_{K \cap \Omega_i} u \,dx|^2}{\sum_K | \frac{1}{|K \cap \Omega_i|} \int_{K \cap \Omega_i} u \,dx|^2}},
\end{equation*}
where $i = 1, 2$ denotes the continuum number, and $K$ is the RVE taken as the coarse block.

Numerical implementation is based on the FEniCS computational package \cite{logg2012automated}. We use the mesh generator Gmsh to generate the unstructured fine grid \cite{geuzaine2009gmsh}. Visualization is performed using the ParaView software \cite{ayachit2015paraview}. 

\subsection{Example 1}

\begin{figure}[hbt!]
\centering
\includegraphics[width=0.3\textwidth]{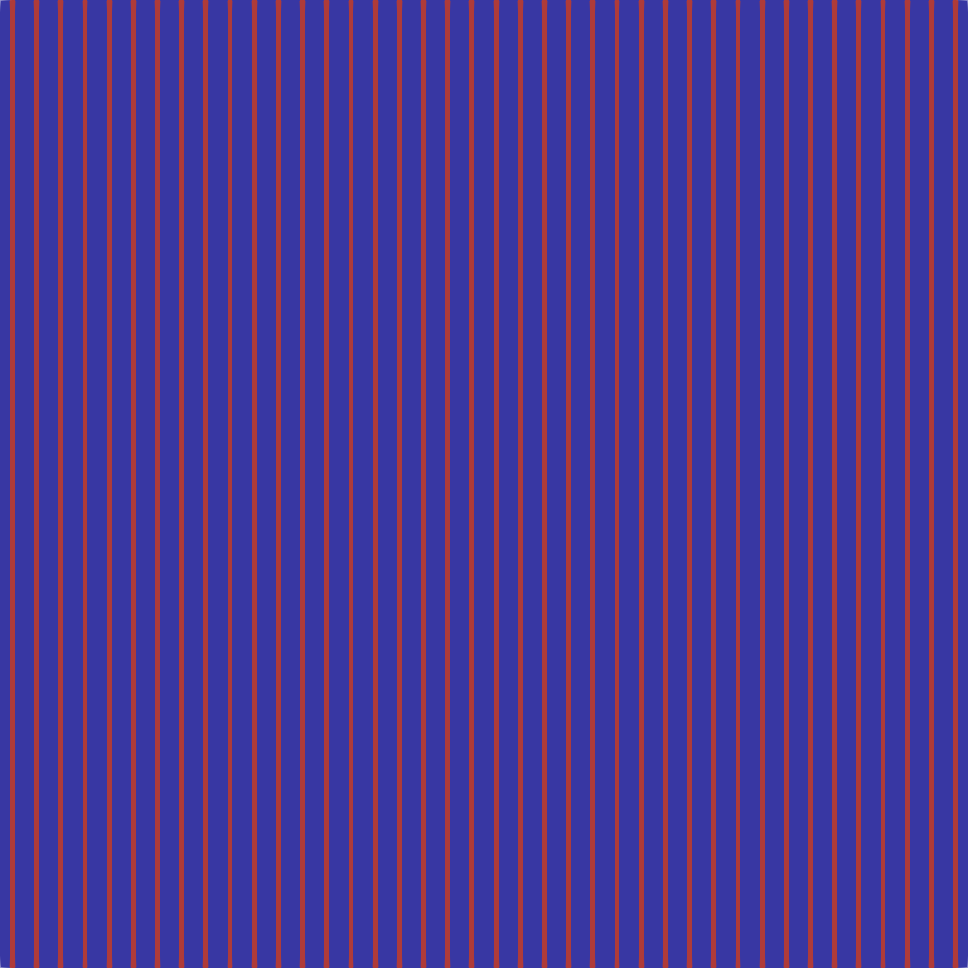}
\caption{Microstructure ($\Omega_1$ is blue, $\Omega_2$ is red) for Example 1.}
\label{fig:microstructure_ex_1}
\end{figure}

The first numerical example considers an anisotropic high-contrast heterogeneous medium. We present the microstructure in Figure \ref{fig:microstructure_ex_1}, where the region $\Omega_1$ is blue, and the region $\Omega_2$ is red. We set the following heterogeneous Lam\'e parameters $\lambda$ and $\mu$ and the permeability coefficient $\kappa$

\begin{equation}\label{eq:material_parameters_ex_1}
\lambda = \begin{cases}
    10^{9}, \quad x \in \Omega_1,\\
    10^{5}, \quad x \in \Omega_2
\end{cases}, \quad
\mu = \begin{cases}
    10^{9}, \quad x \in \Omega_1,\\
    10^{5}, \quad x \in \Omega_2
\end{cases}, \quad
\kappa = \begin{cases}
    10^{-10}, \quad x \in \Omega_1,\\
    10^{-6}, \quad x \in \Omega_2
\end{cases}.
\end{equation}

For other parameters, we set the Biot coefficient $\alpha = 0.8$ and the Biot modulus $M = 10^6$. We set the duration time $t_{max} = 5$ with 50 time steps using an implicit time scheme. Note that all these parameters are model ones. For right-hand sides of displacements and pressure equations, we set the following heterogeneous functions
\begin{equation}\label{eq:sources_ex_1}
\begin{gathered}
f_{1} = - \sin\left(2\pi\, x_1 \right)\sin\left(\pi\, x_2\right) \, 10^4, \quad
f_{2} = \sin\left(\pi\, x_1\right)\sin\left(\pi\, x_2\right) \, 10^7,\\
g = 0.15 \, \exp(-40[(x_1-0.5)^2+(x_2-0.5)^2]),
\end{gathered}
\end{equation}
where $f_{1}$ and $f_{2}$ are the body forces in $x_1$ and $x_2$ directions, respectively, and $g$ is the source term for pressure.

We set the initial value $u_0 = (0, 0)$ for displacements and $p_0 = 10^6$ for pressure. For boundary conditions, we fix displacements in $x_1$ on the right boundary and $x_2$ on the bottom boundary and set free traction on the other boundaries. For pressure, we set no-flow on all the boundaries. Therefore, we have the following boundary conditions
\begin{equation}\label{eq:bcs_ex_1}
\begin{split}
&u_1 = 0, \quad (\sigma \cdot \nu)_2 = 0, \quad x \in \Gamma_R,\\
&u_2 = 0, \quad (\sigma \cdot \nu)_1 = 0, \quad x \in \Gamma_B,\\
&\sigma \cdot \nu = 0, \quad x \in \Gamma_L \cup \Gamma_T,\\
&-\kappa \frac{\partial p}{\partial \nu} = 0, \quad x \in \partial \Omega,
\end{split}
\end{equation}
where $\nu$ is a unit outward normal vector, and $\Gamma_L$, $\Gamma_R$, $\Gamma_B$, and $\Gamma_T$ denote the left, right, bottom, and top boundaries, respectively, such that $\partial \Omega = \Gamma_L \cup \Gamma_R \cup \Gamma_B \cup \Gamma_T$.

\begin{figure}[hbt!]
\centering
\includegraphics[height=0.24\textheight]{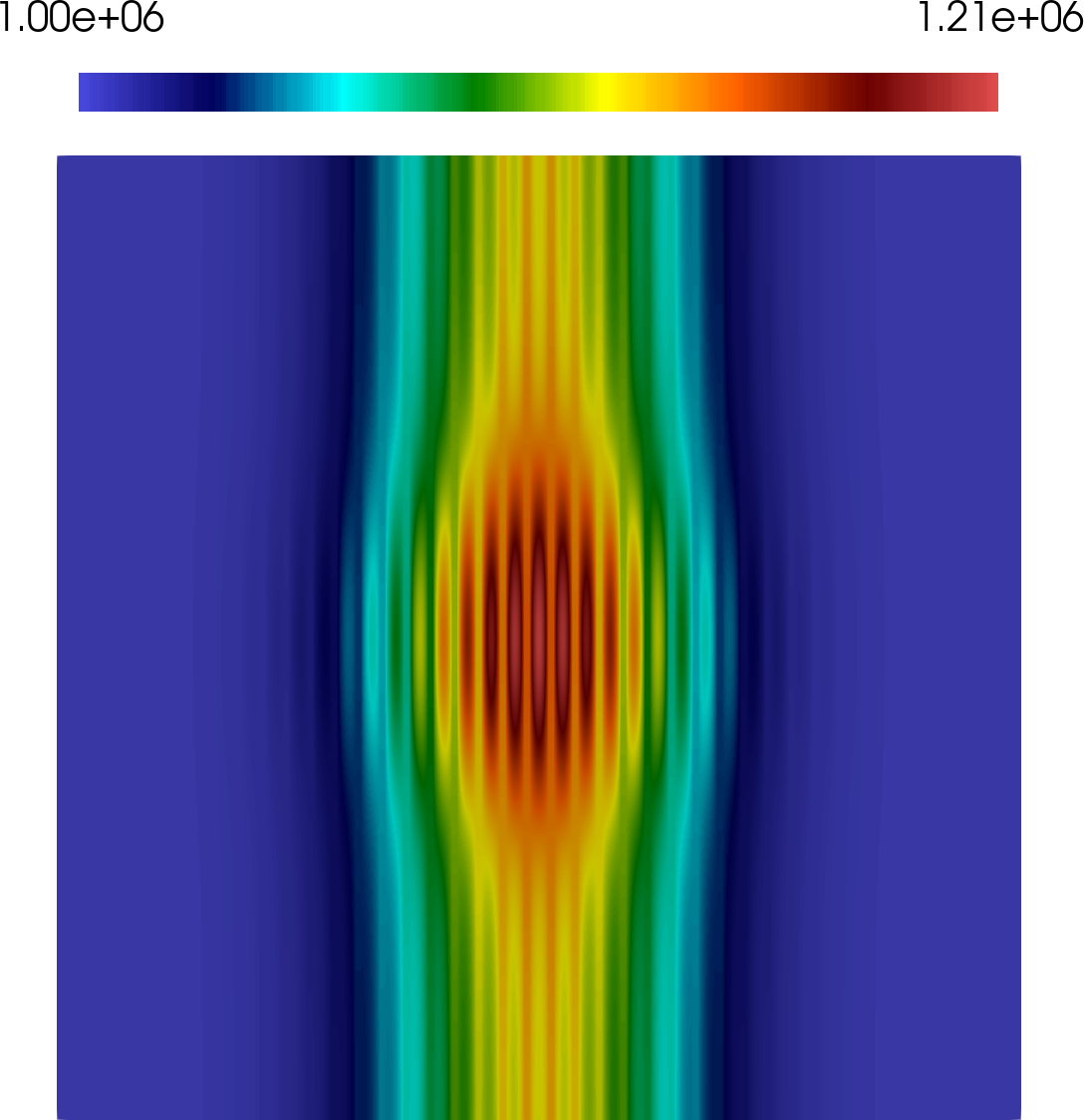}
\includegraphics[height=0.24\textheight]{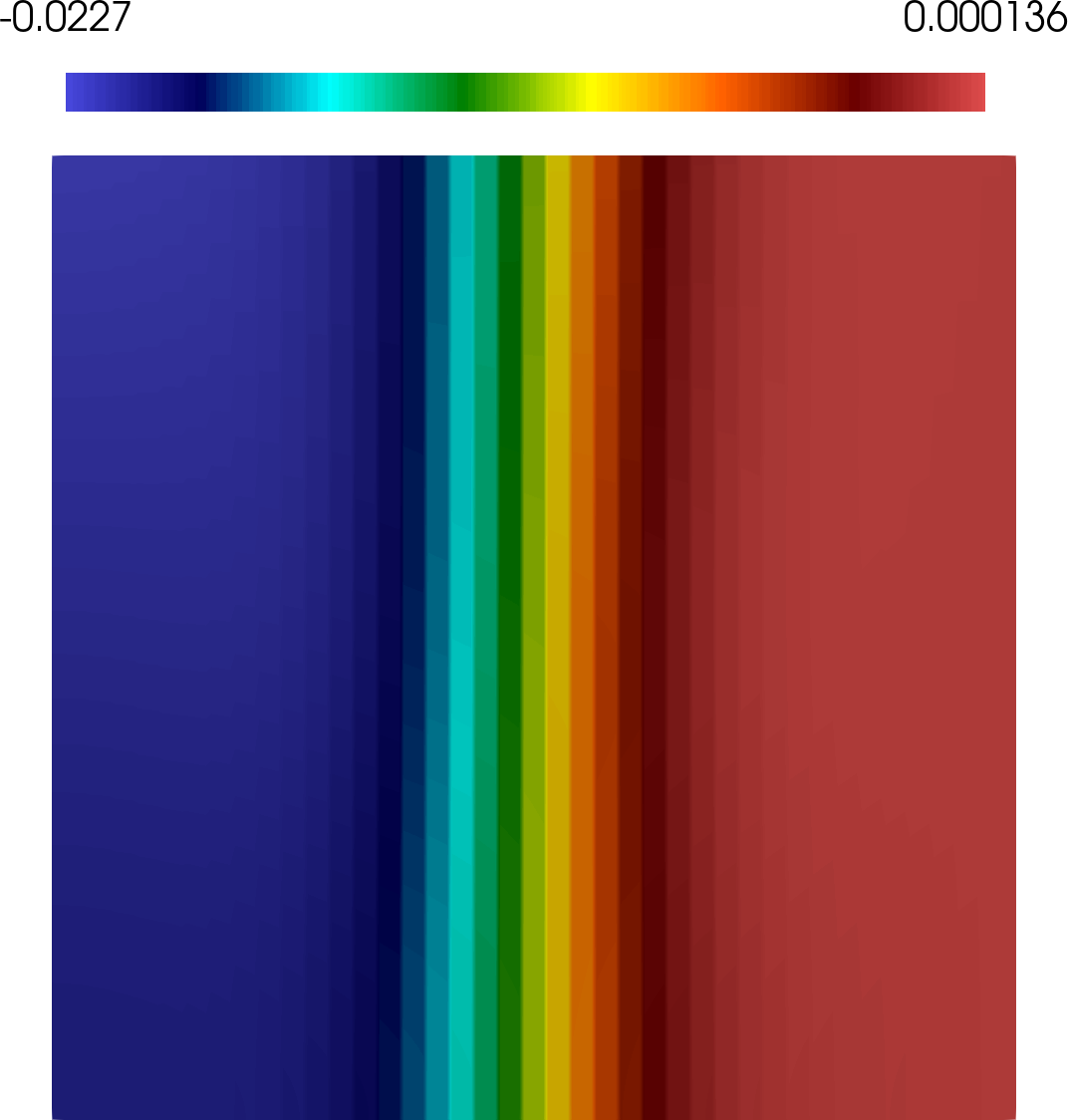}
\includegraphics[height=0.24\textheight]{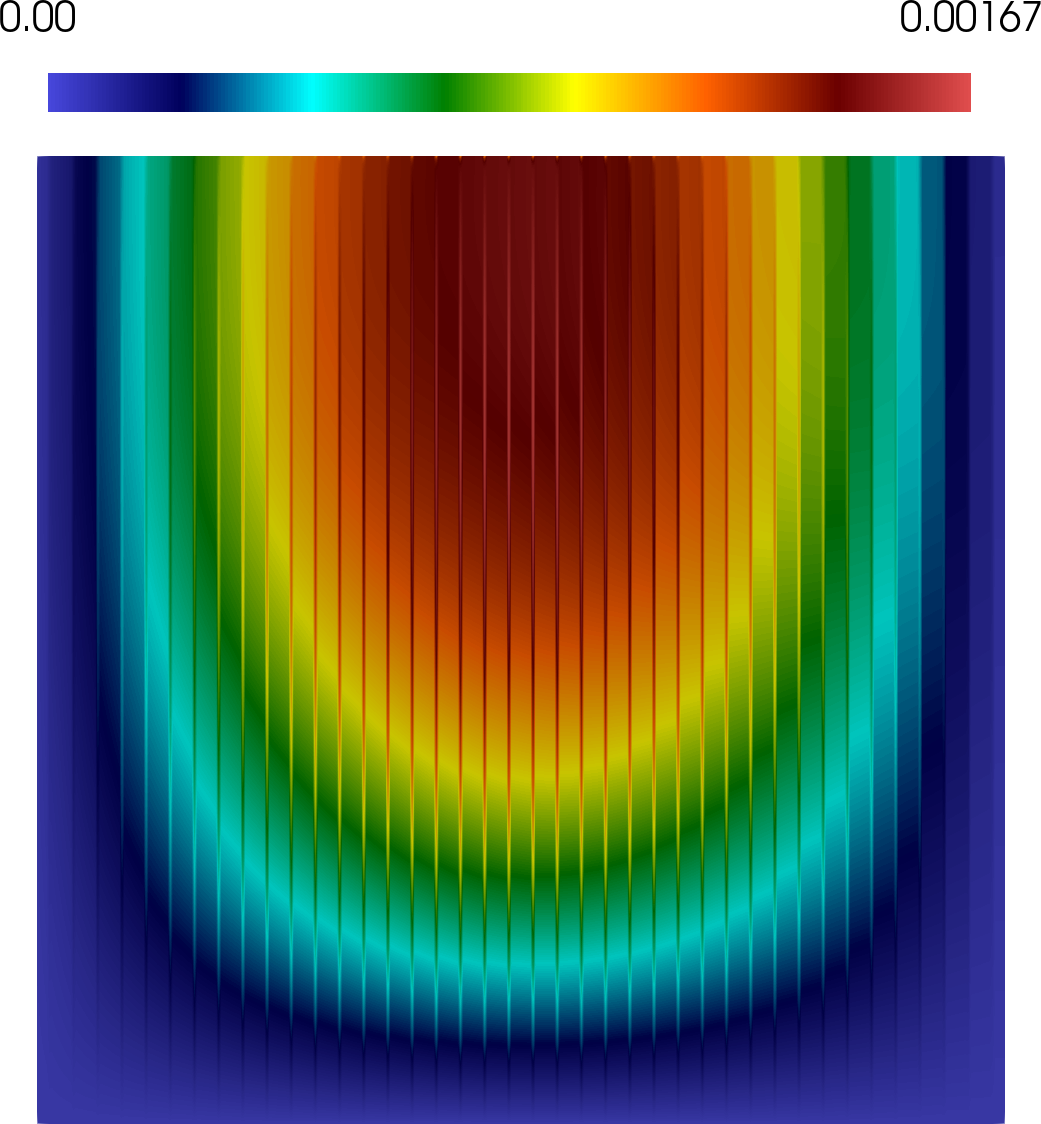}
\caption{Distributions of pressure and displacements in $x_1$ and $x_2$ directions (from left to right) at the final time on the fine grid for Example 1.}
\label{fig:fine_results_ex_1}
\end{figure}

Figure \ref{fig:fine_results_ex_1} depicts distributions of pressure and displacements in $x_1$ and $x_2$ directions (from left to right) at the final time on the fine grid. One can see the influence of the microstructure on the obtained solutions. We observe the fluid diffusion through vertical channels from the source in the middle of the domain. However, it also slowly diffuses to the sides. Displacement in the $x_1$ direction corresponds to the horizontal stretching due to the body force and pressure gradient. The pattern of the stretching is also affected by stiff and soft channels and free traction on the left side. Displacement in the $x_2$ direction depicts vertical upward stretching. In general, the obtained fine-scale solutions demonstrate an anisotropic nature and correspond to the modeled process.

\begin{figure}[hbt!]
\centering
\begin{subfigure}{\textwidth}
\centering
\includegraphics[height=0.24\textheight]{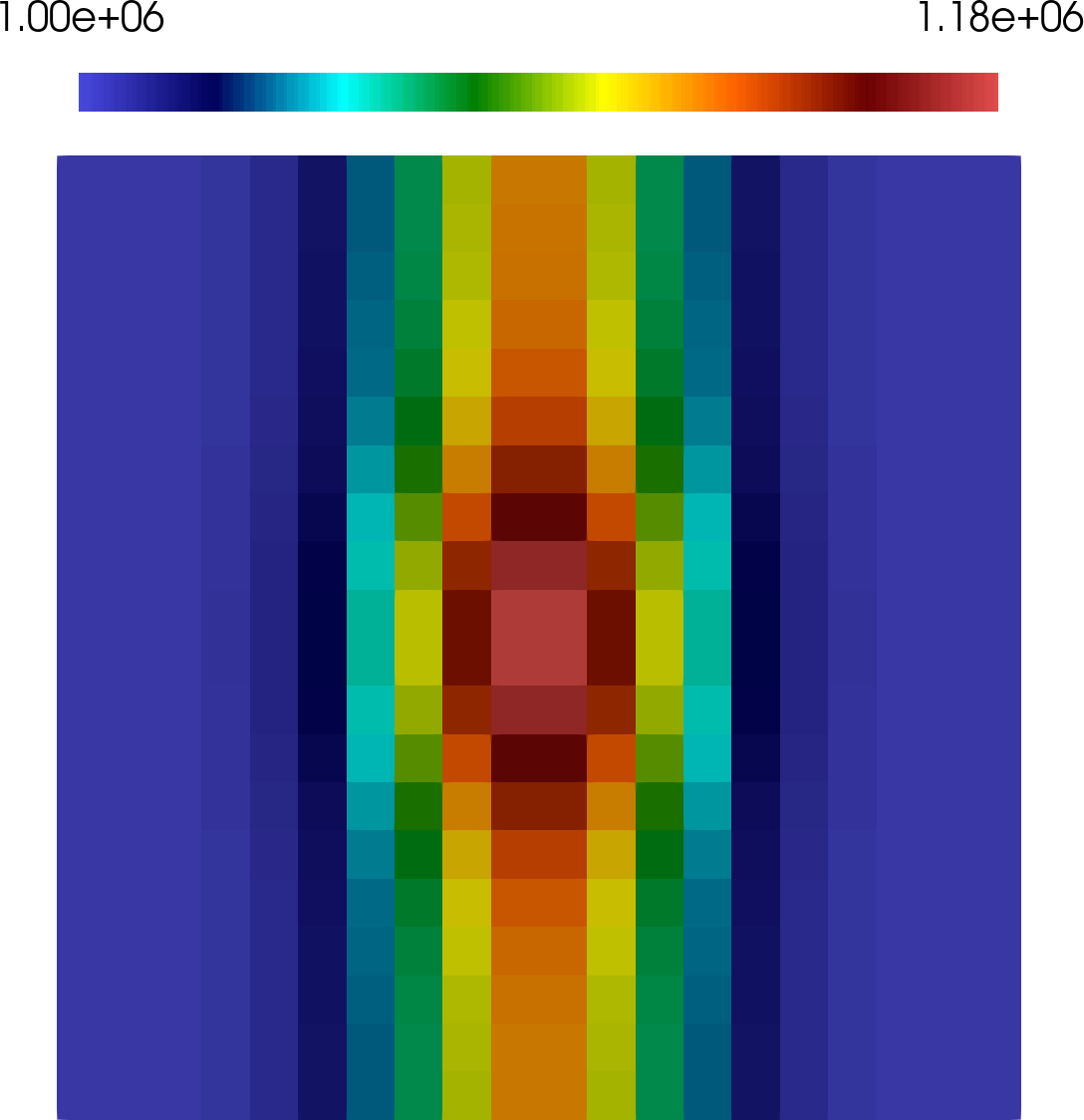}
\includegraphics[height=0.24\textheight]{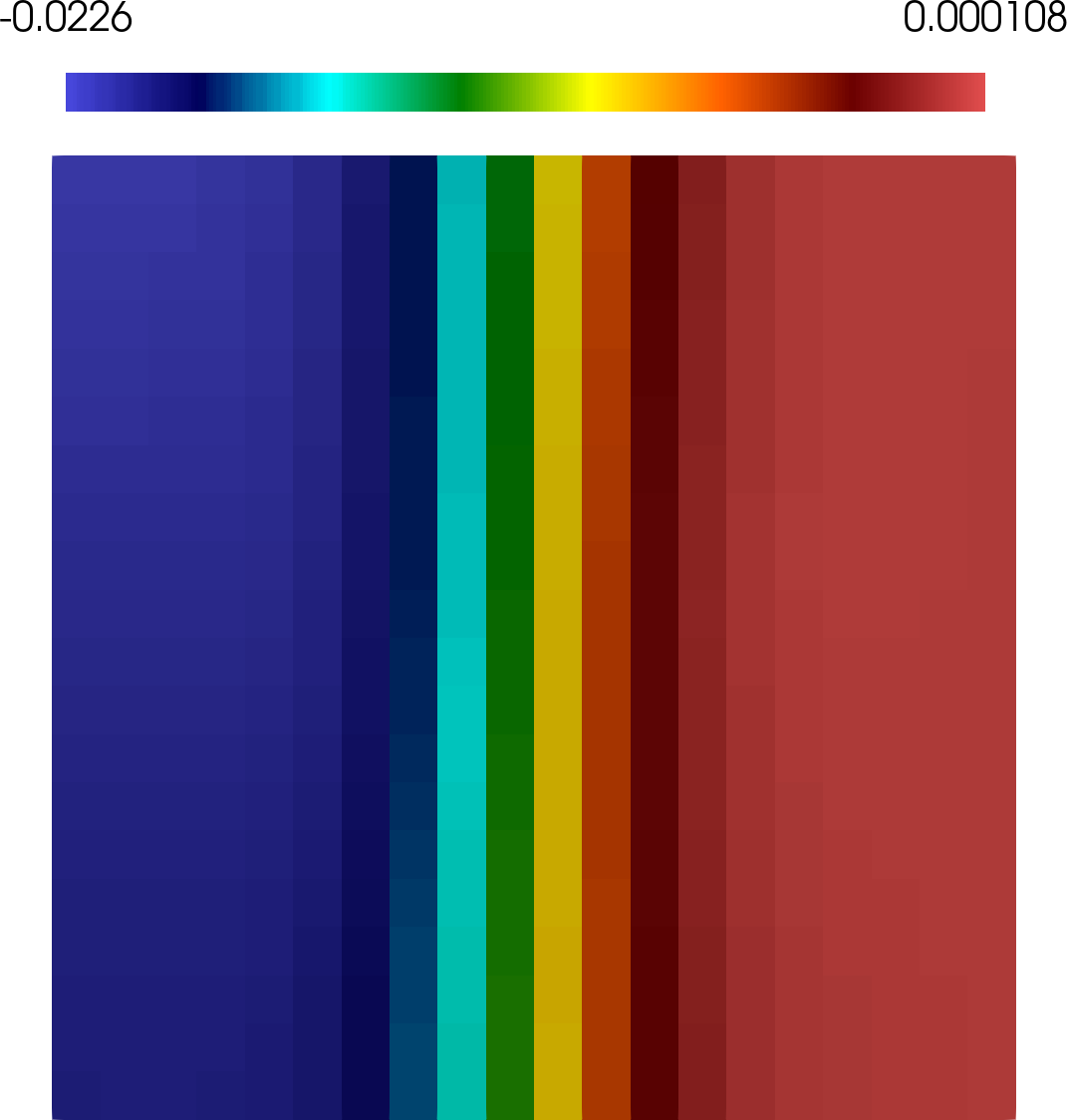}
\includegraphics[height=0.24\textheight]{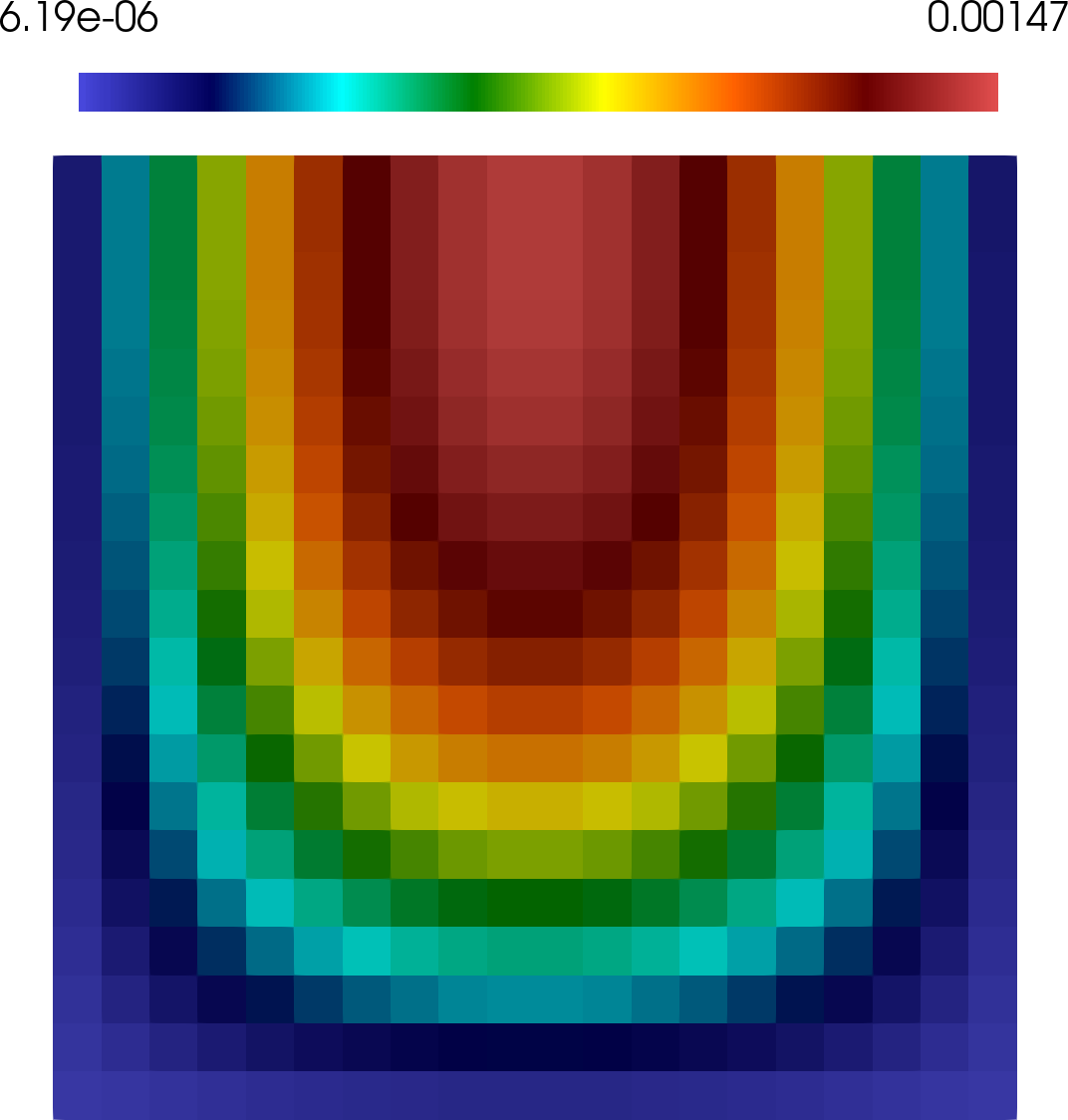}
\caption{Reference average solution}
\end{subfigure}
\begin{subfigure}{\textwidth}
\centering
\includegraphics[height=0.24\textheight]{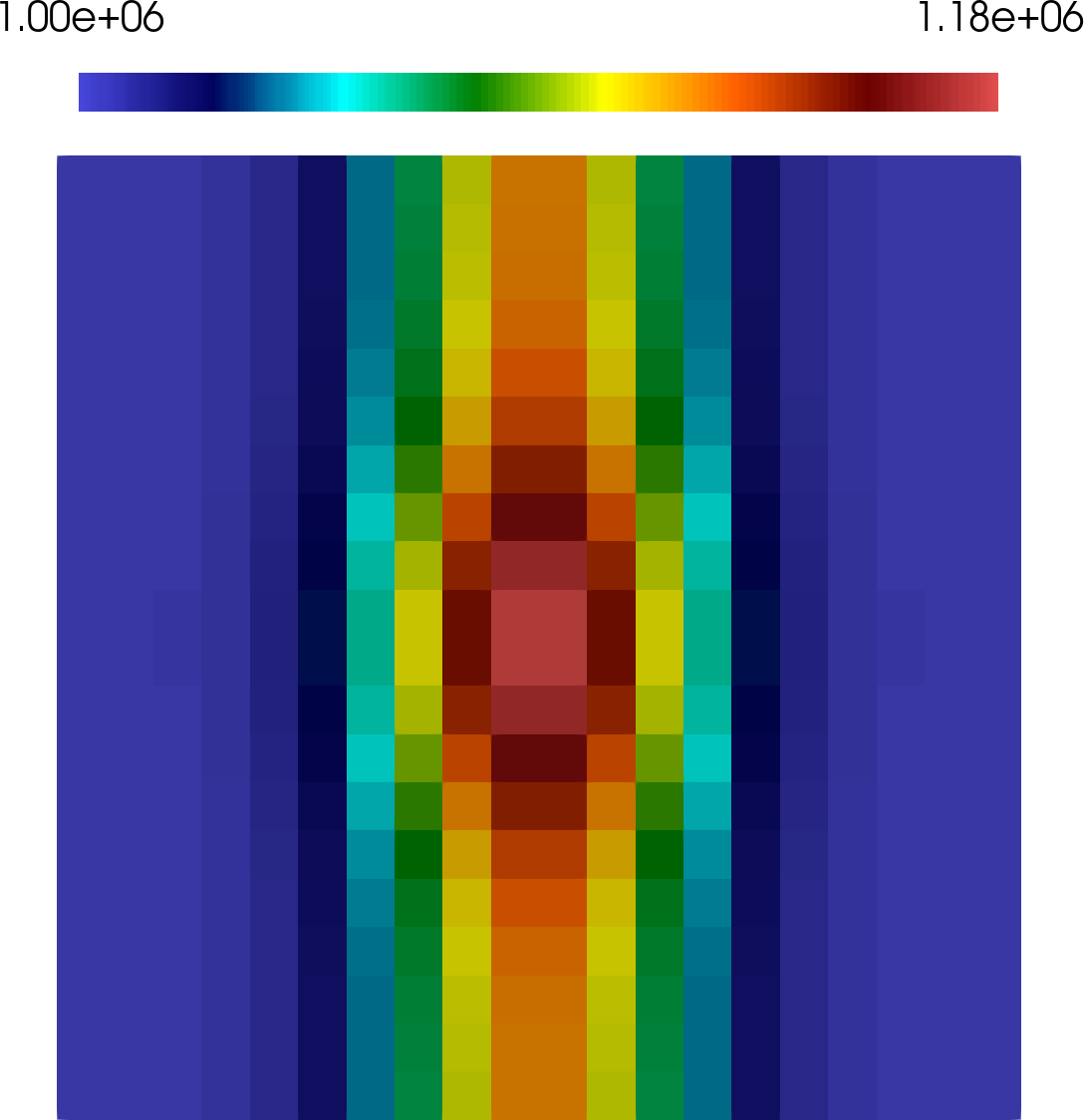}
\includegraphics[height=0.24\textheight]{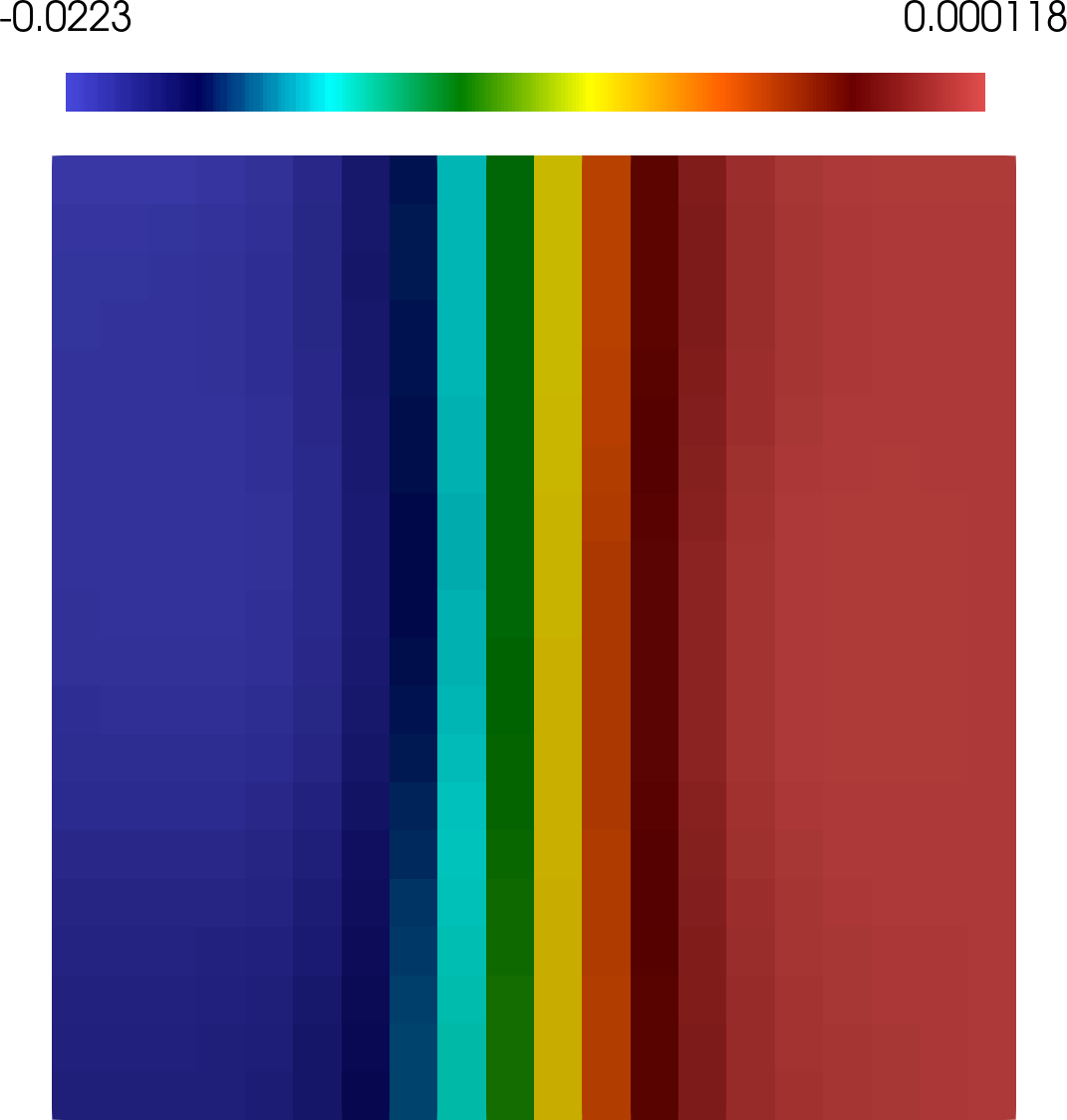}
\includegraphics[height=0.24\textheight]{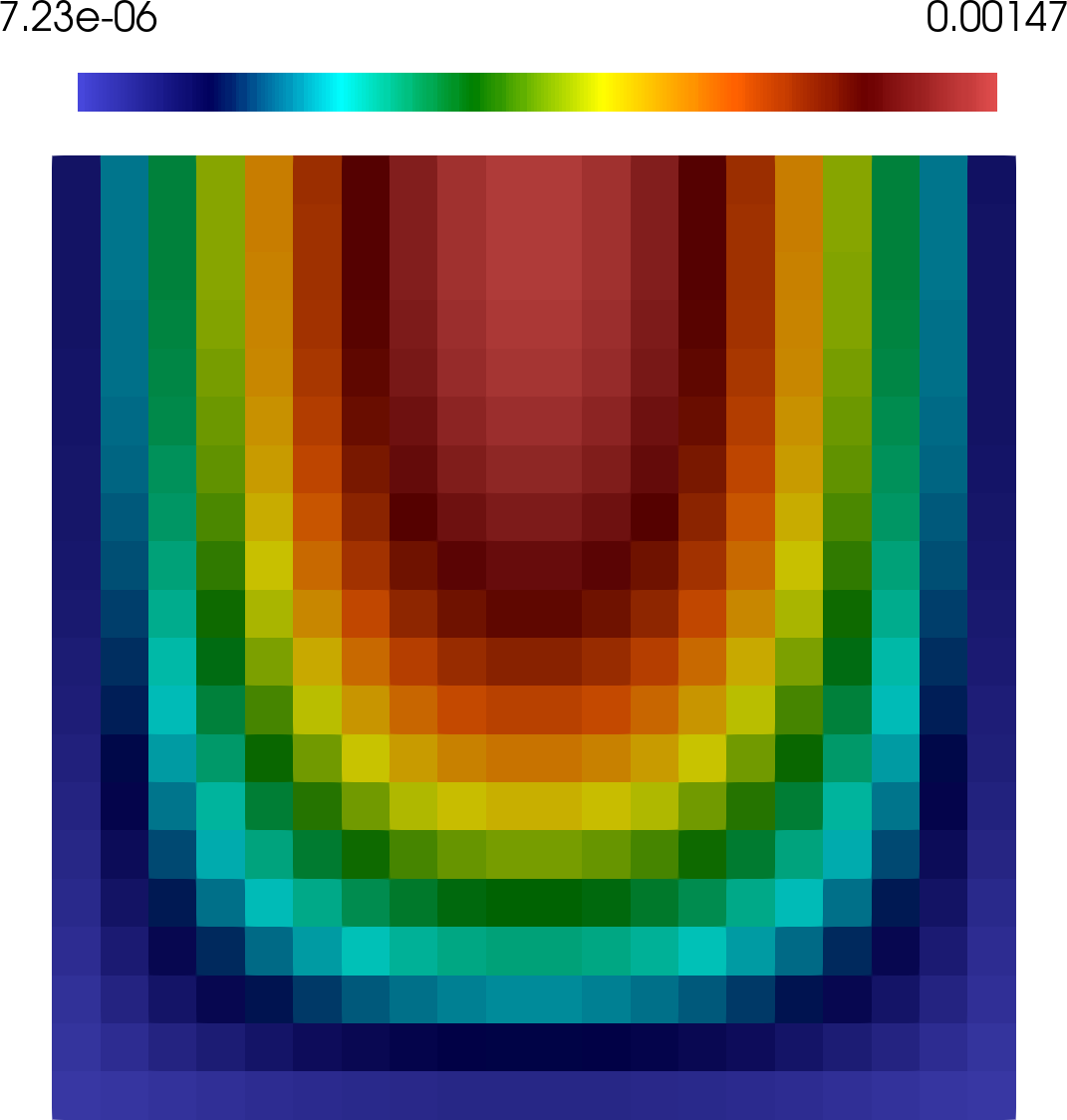}
\caption{Homogenized average solution using the full multicontinuum model}
\end{subfigure}
\begin{subfigure}{\textwidth}
\centering
\includegraphics[height=0.24\textheight]{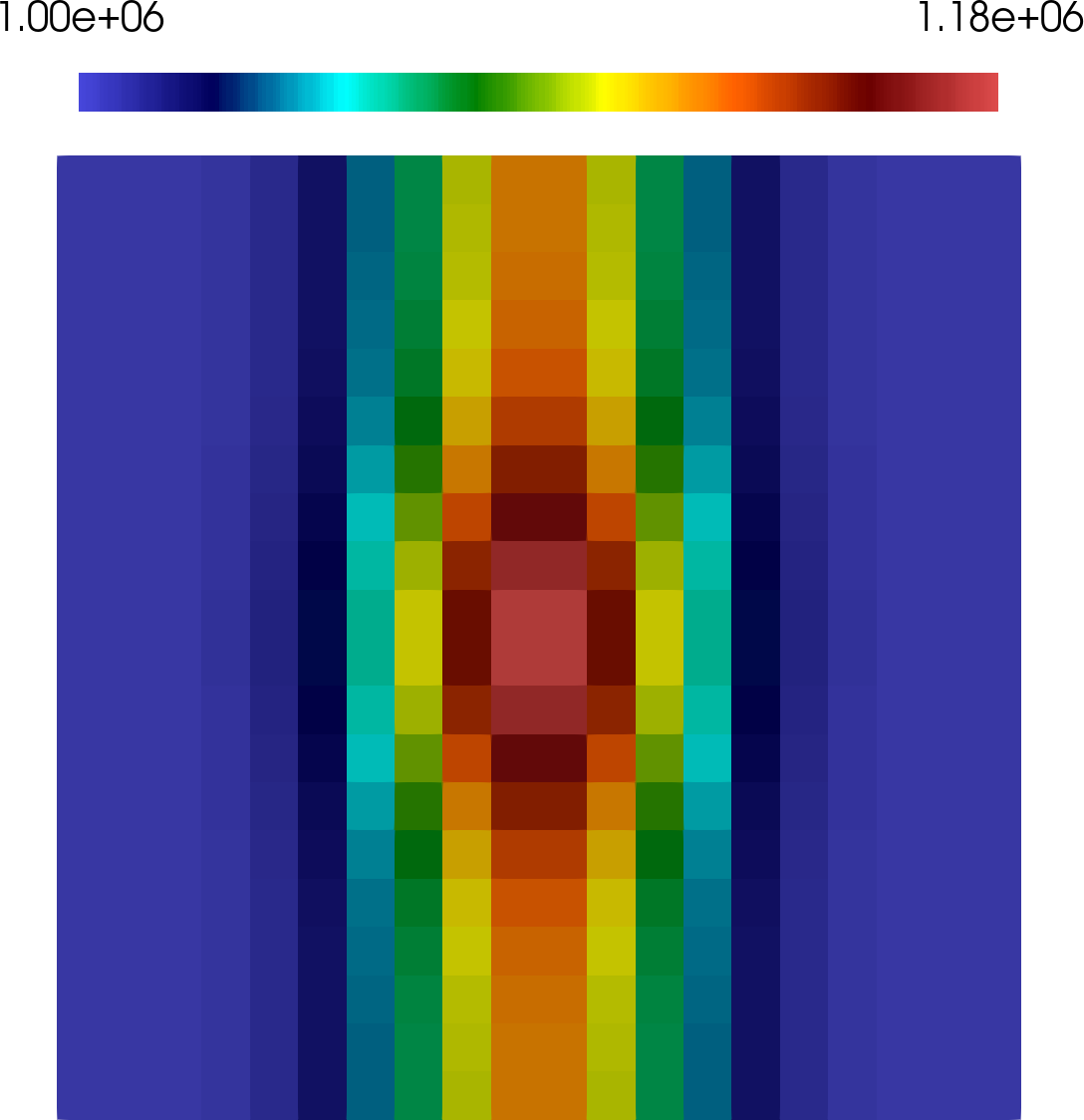}
\includegraphics[height=0.24\textheight]{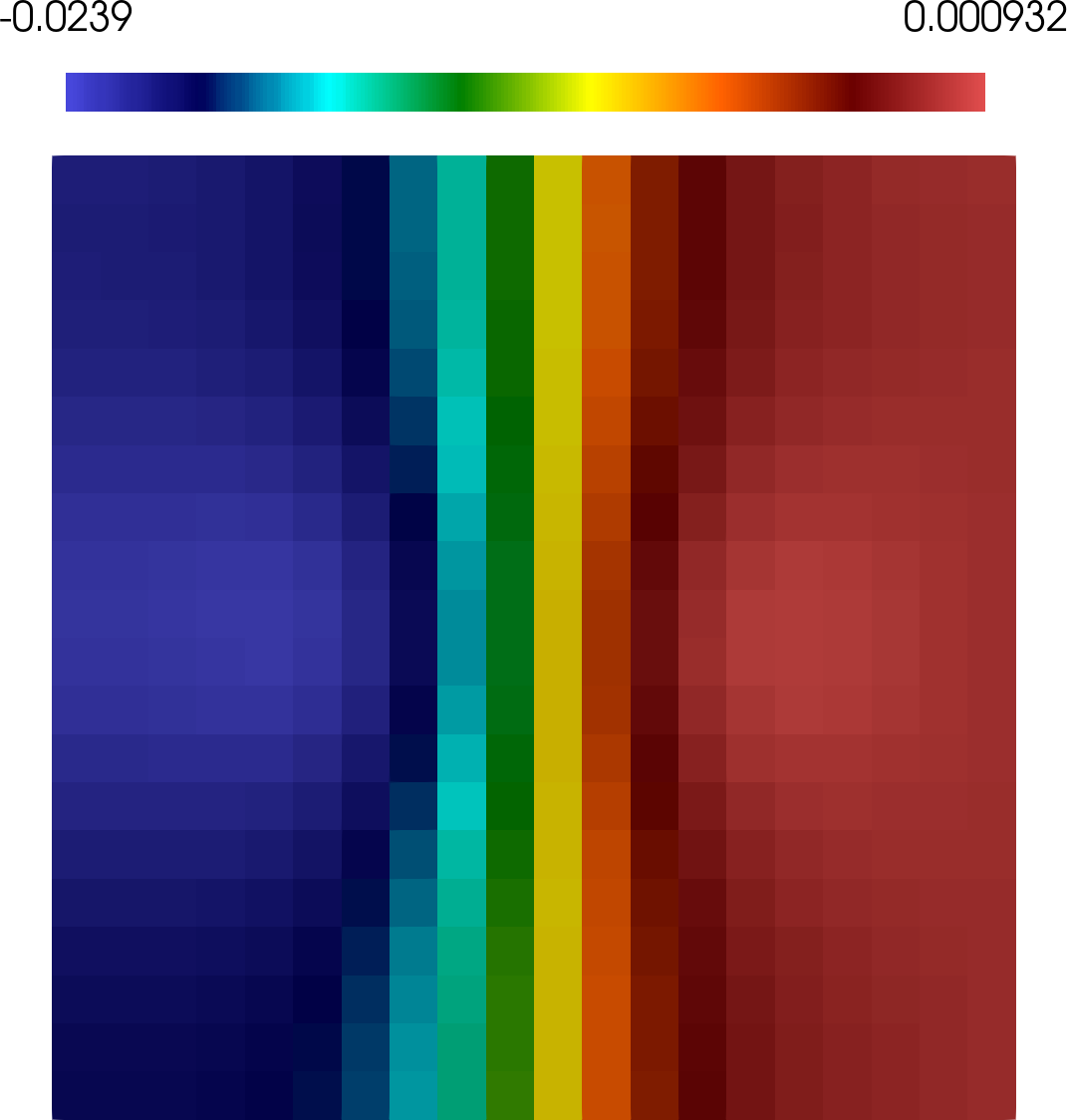}
\includegraphics[height=0.24\textheight]{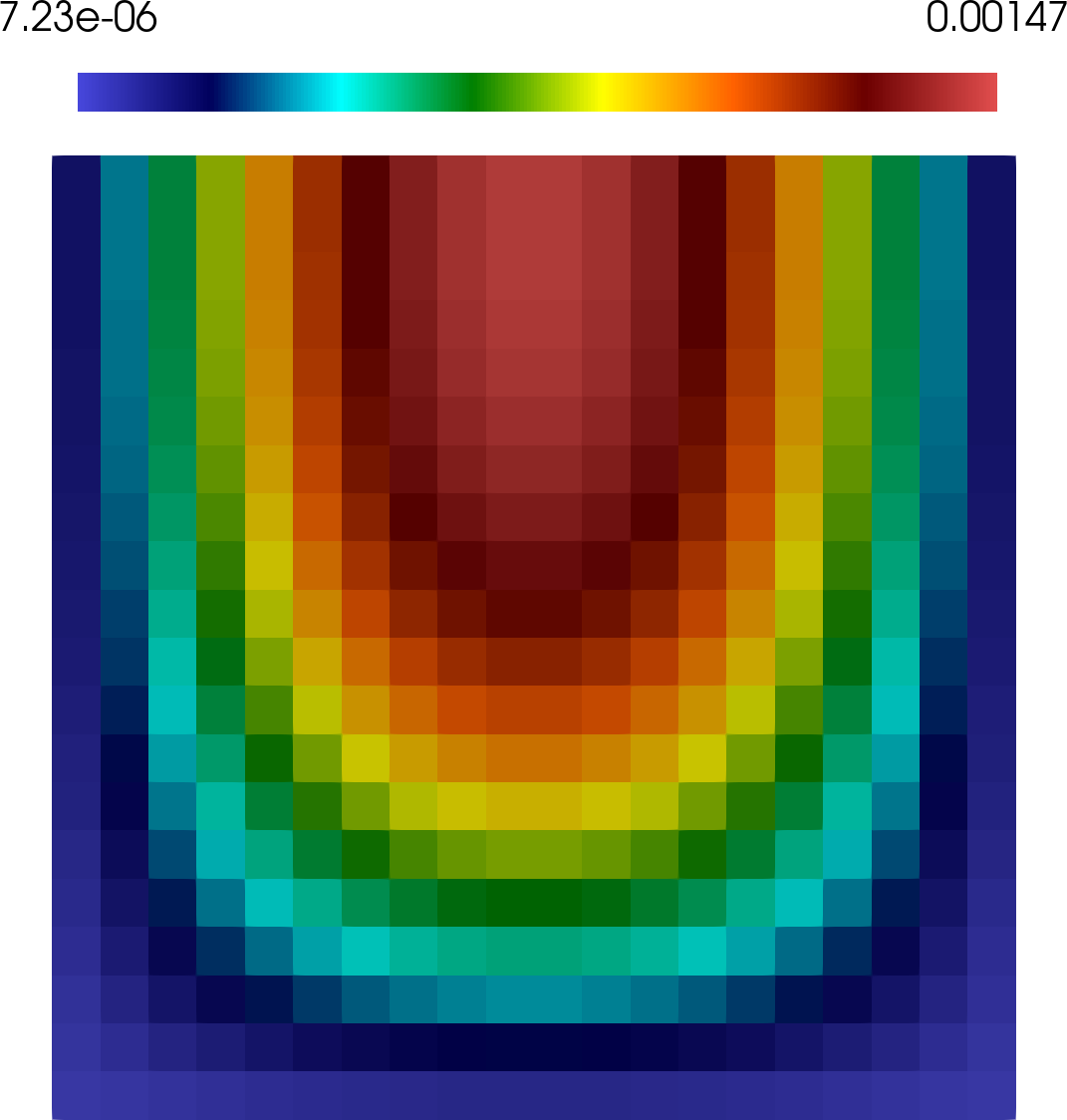}
\caption{Homogenized average solution using the simplified multicontinuum model 2}
\end{subfigure}
\caption{Distributions of average pressure and displacements in $x_1$ and $x_2$ directions (from left to right) in $\Omega_1$ at the final time on the coarse grid $20 \times 20$ for Example 1.}
\label{fig:coarse_results_1_ex_1}
\end{figure}

\begin{figure}[hbt!]
\centering
\begin{subfigure}{\textwidth}
\centering
\includegraphics[height=0.24\textheight]{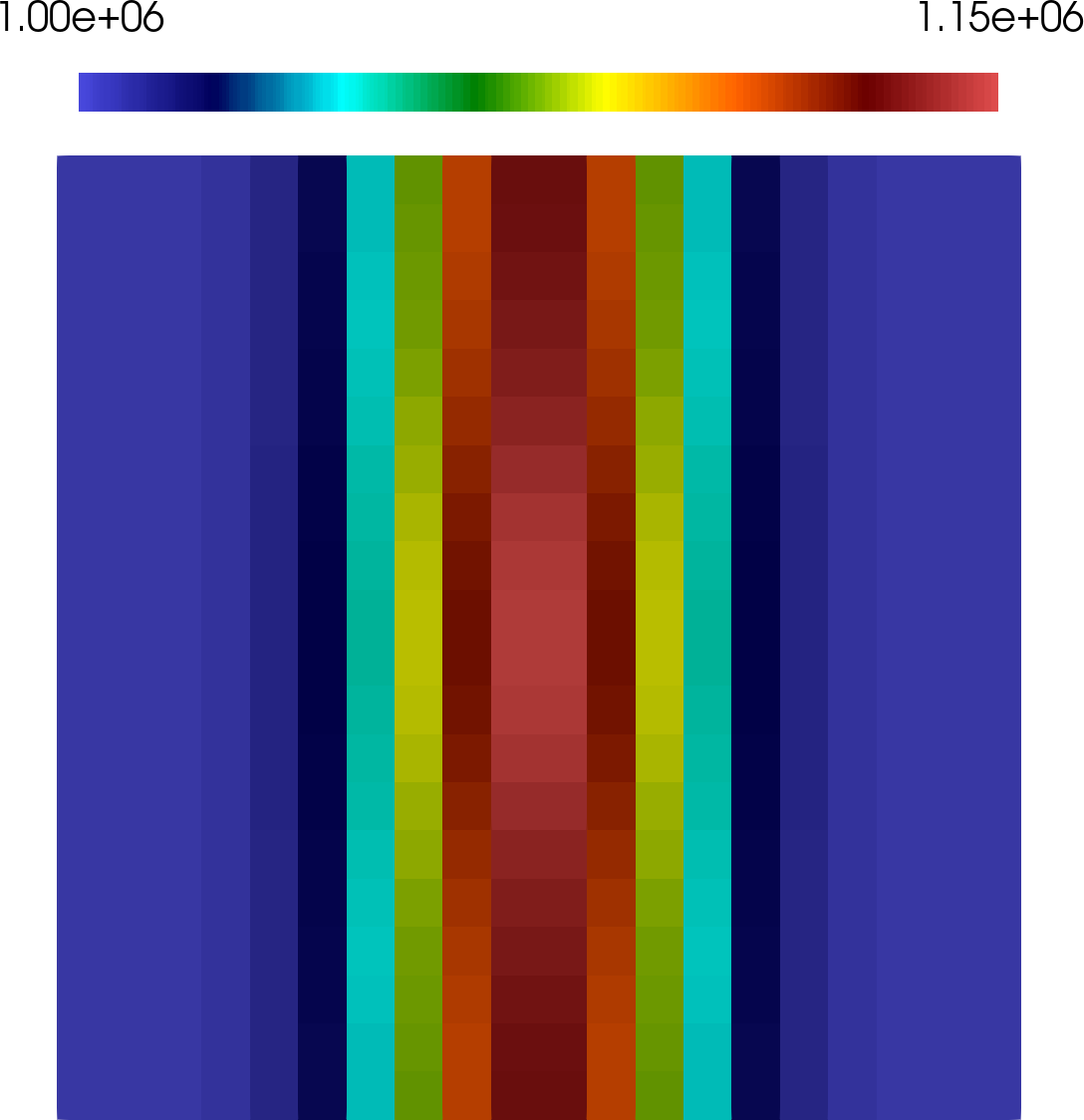}
\includegraphics[height=0.24\textheight]{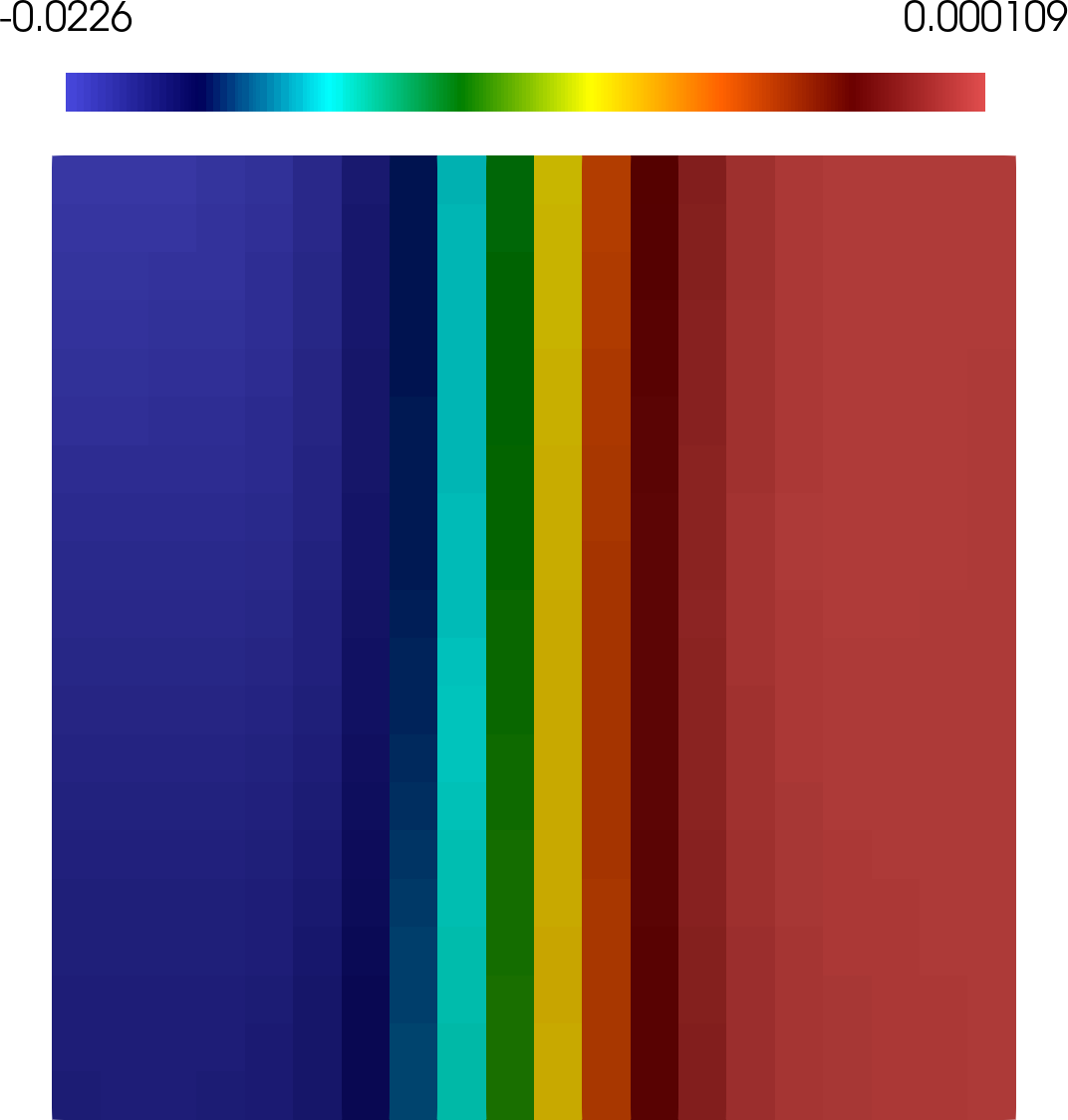}
\includegraphics[height=0.24\textheight]{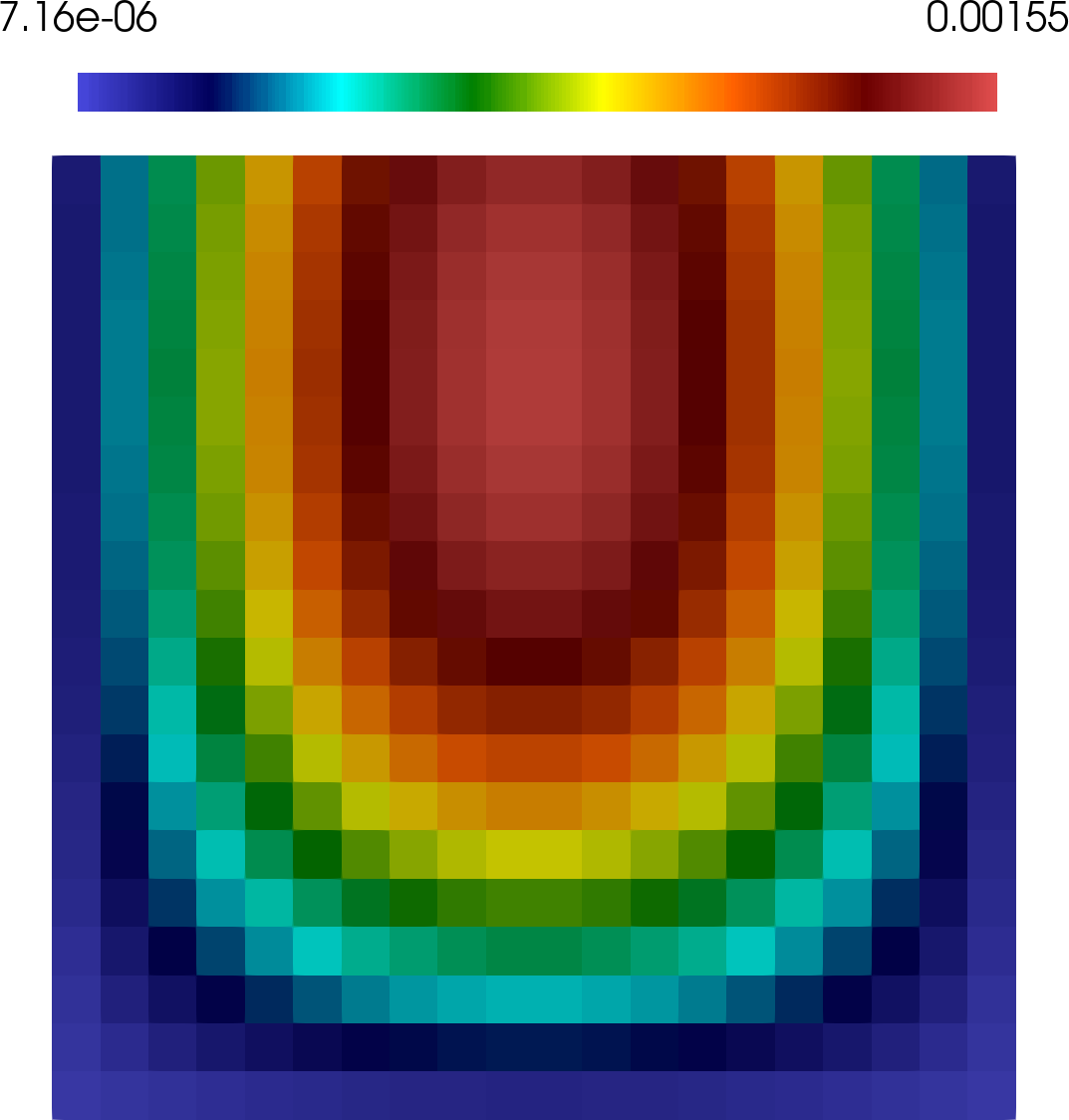}
\caption{Reference average solution}
\end{subfigure}
\begin{subfigure}{\textwidth}
\centering
\includegraphics[height=0.24\textheight]{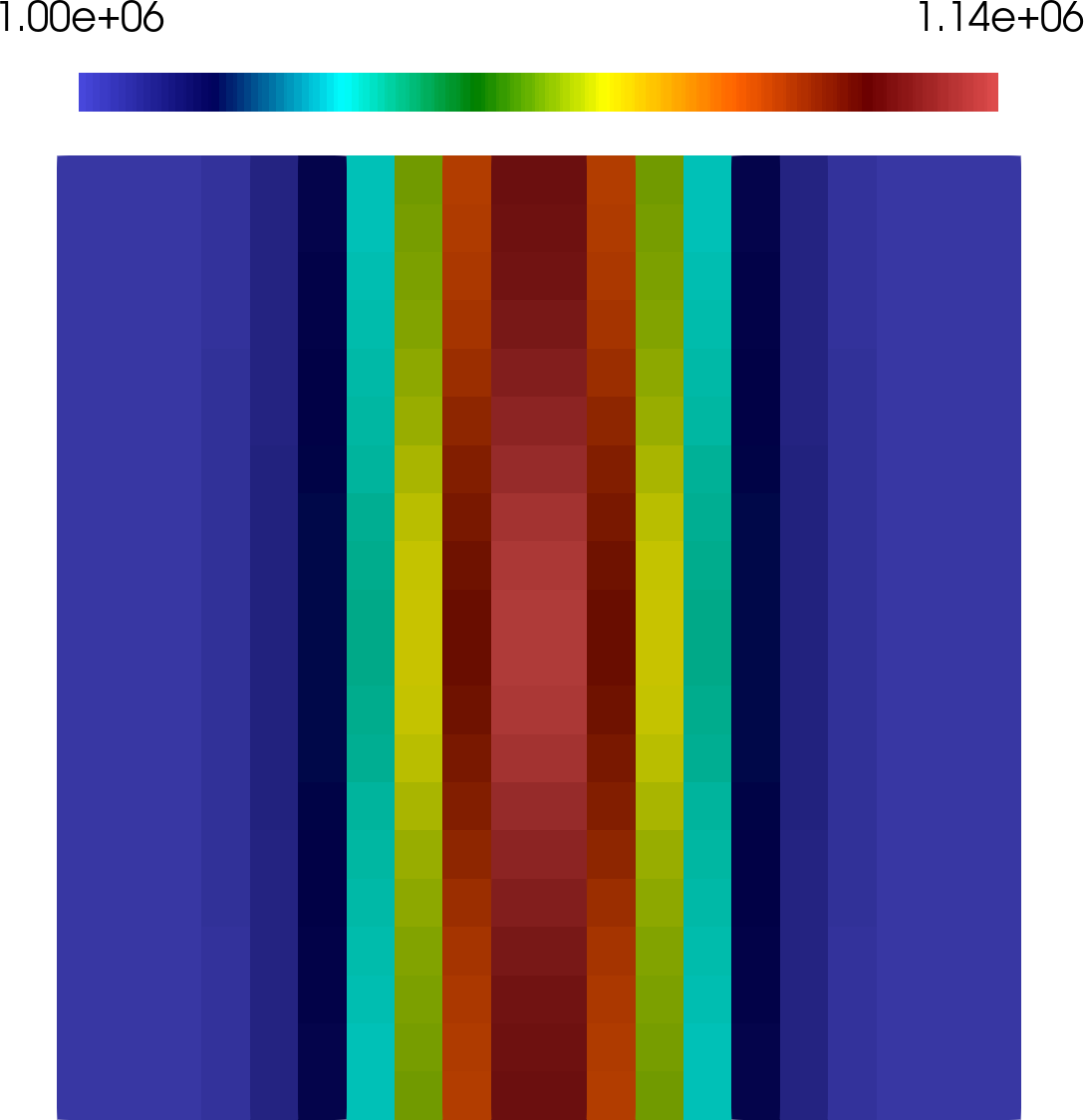}
\includegraphics[height=0.24\textheight]{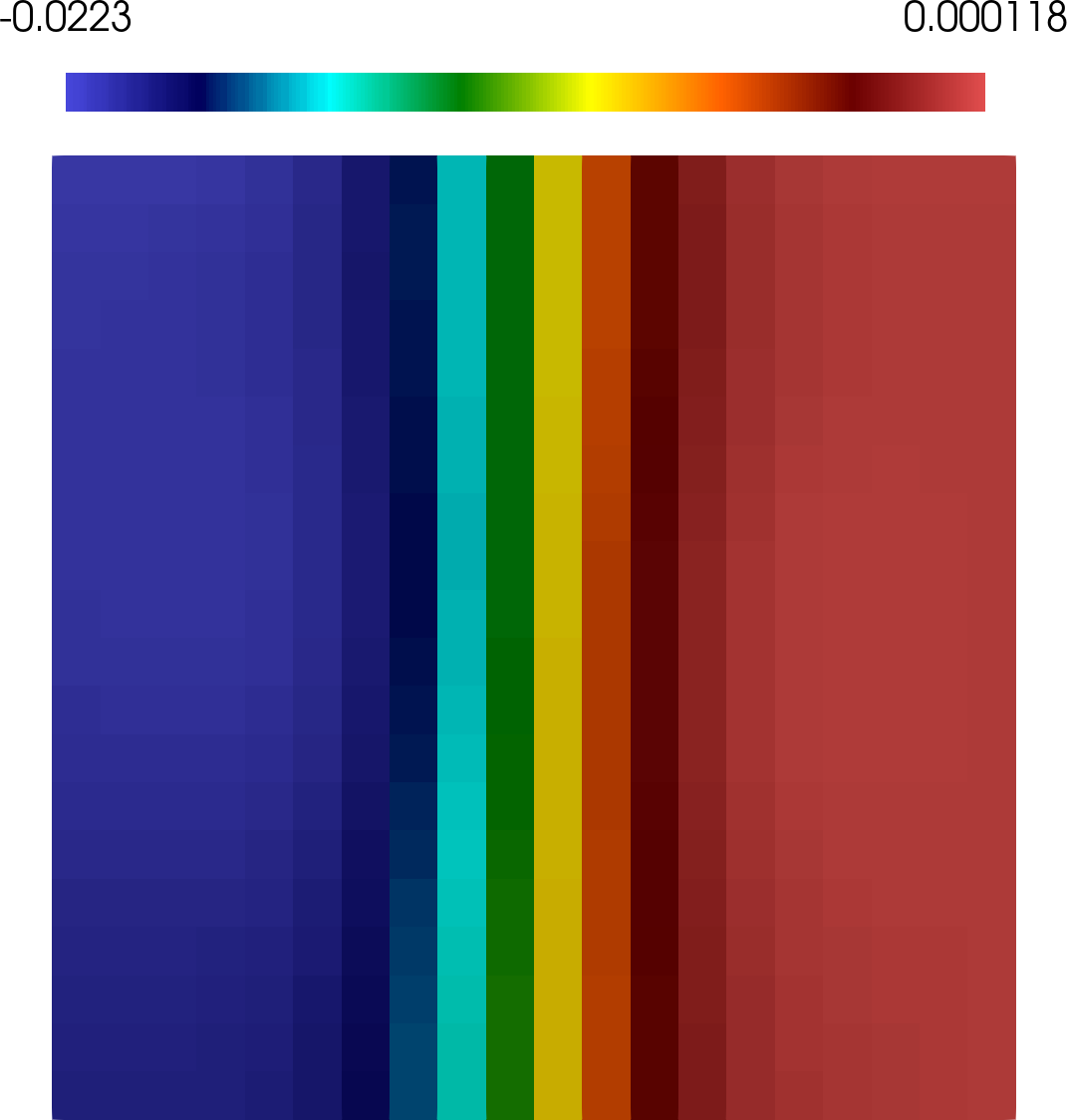}
\includegraphics[height=0.24\textheight]{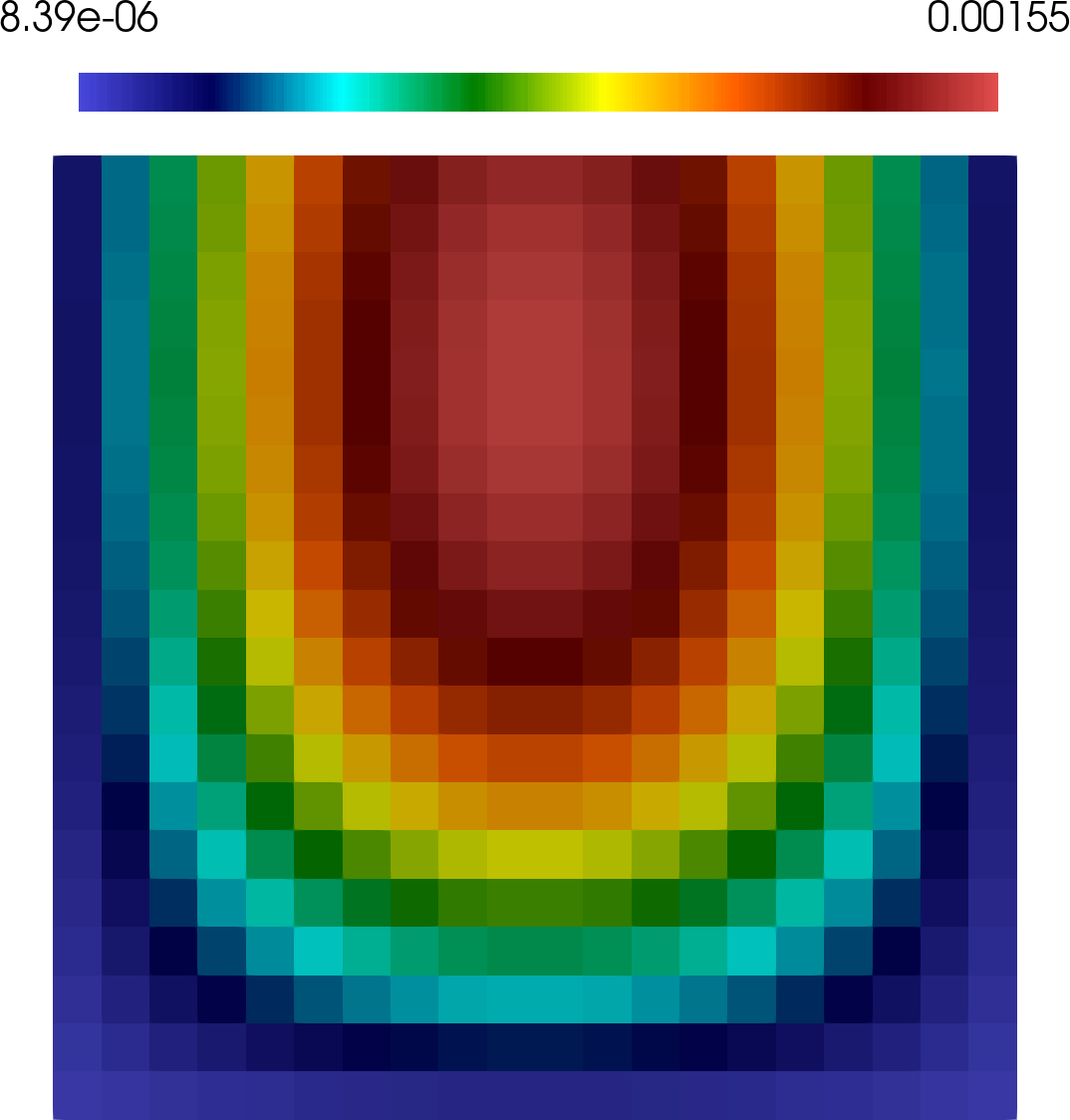}
\caption{Homogenized average solution using the full multicontinuum model}
\end{subfigure}
\begin{subfigure}{\textwidth}
\centering
\includegraphics[height=0.24\textheight]{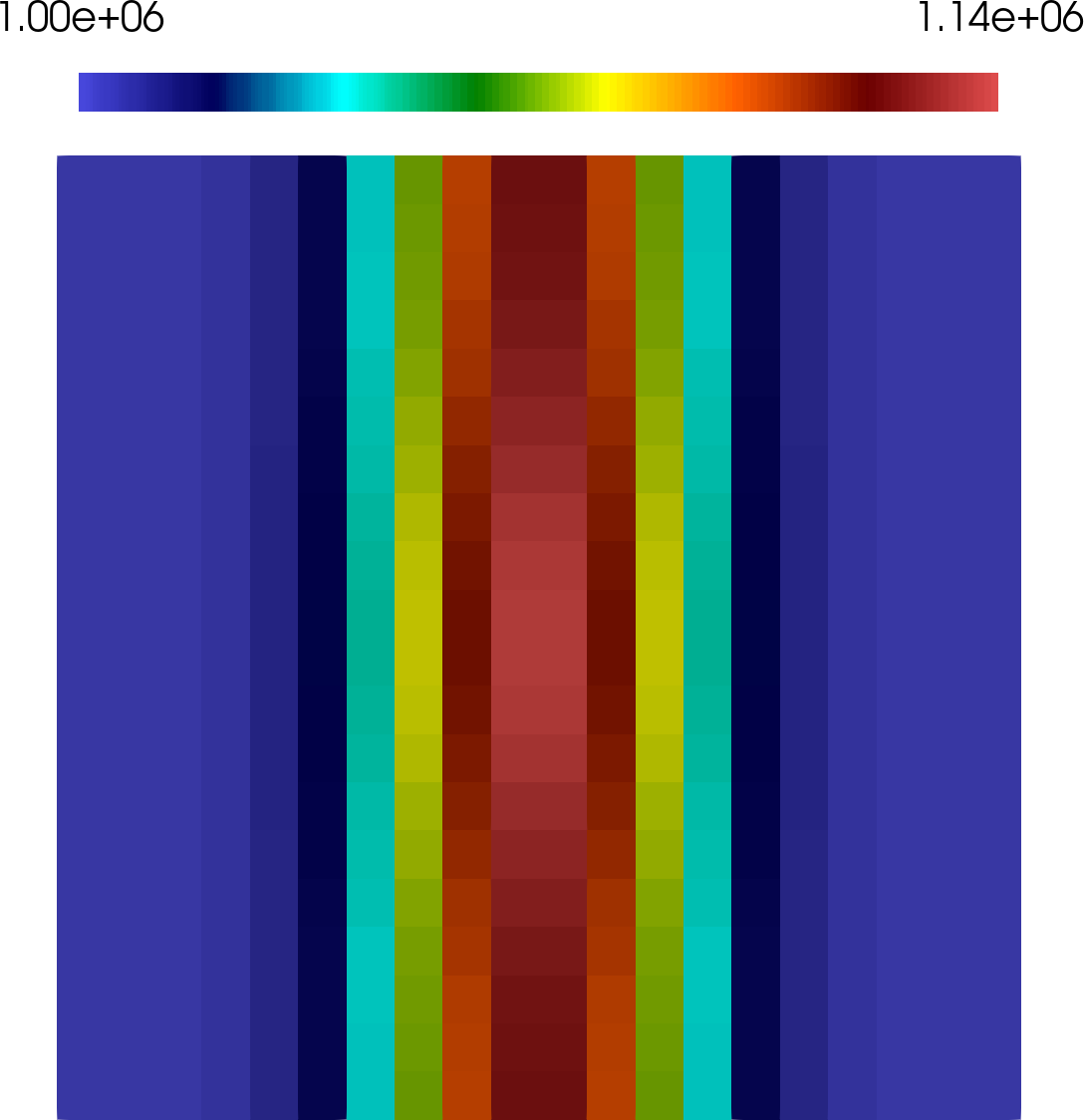}
\includegraphics[height=0.24\textheight]{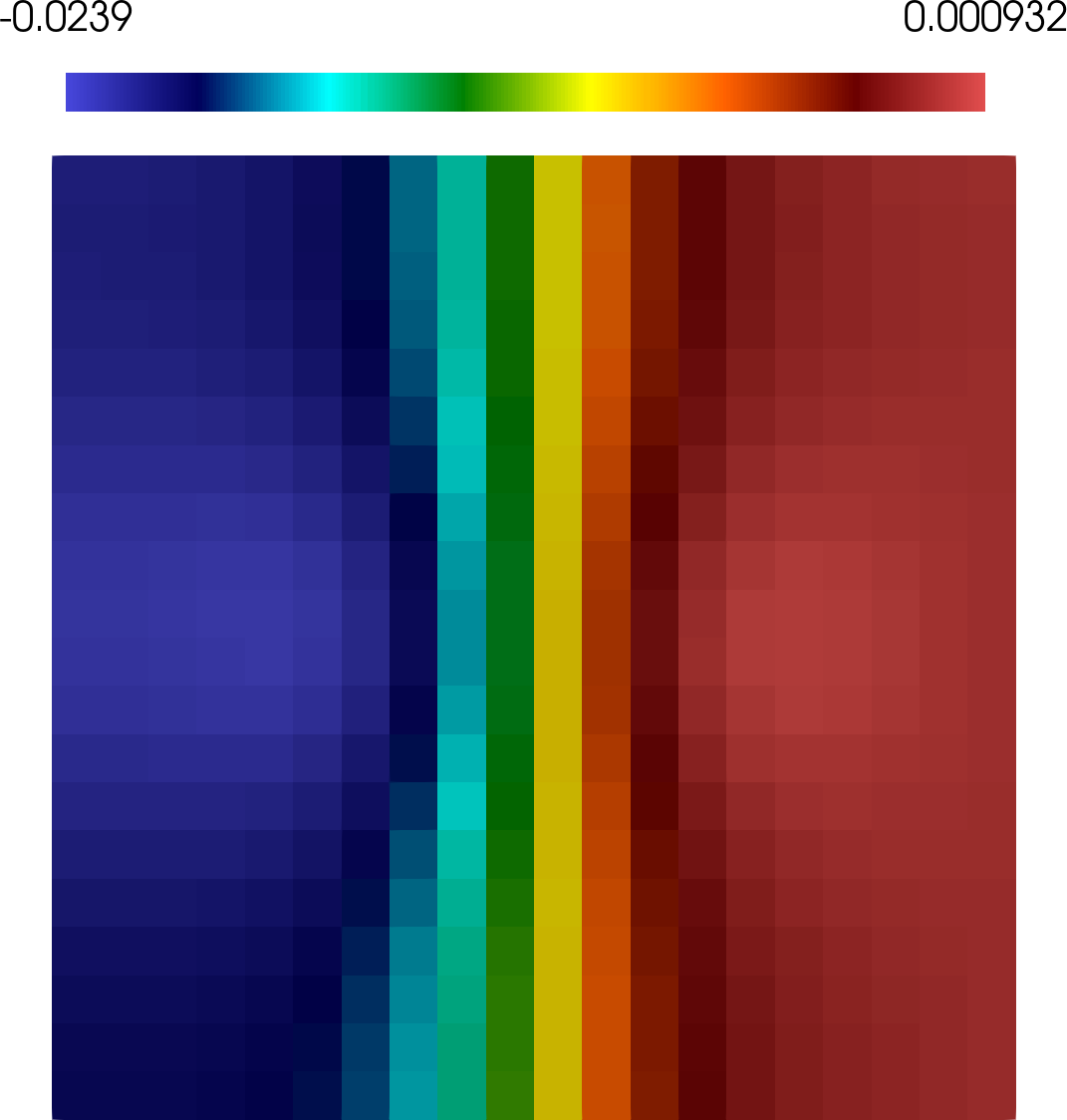}
\includegraphics[height=0.24\textheight]{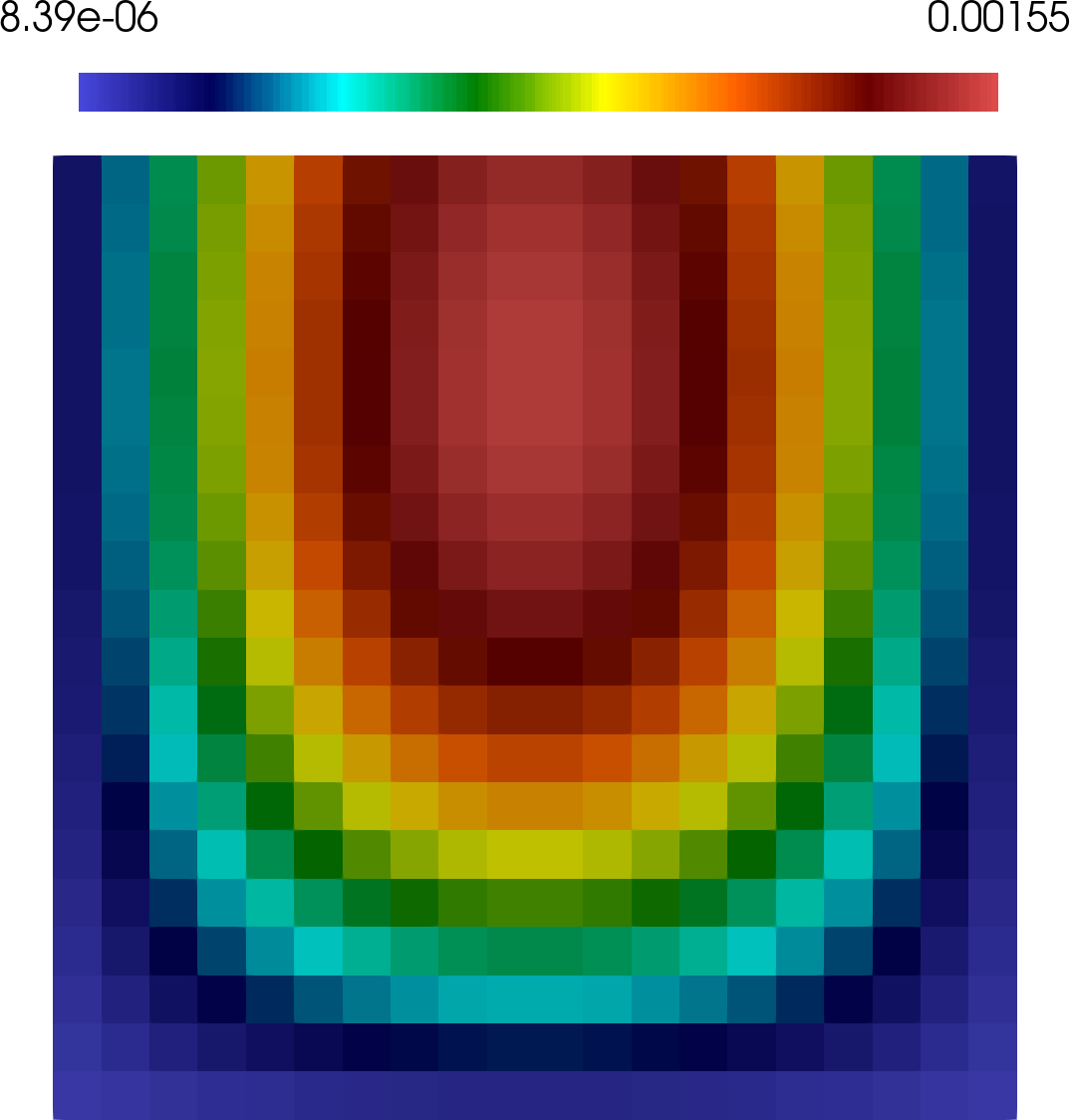}
\caption{Homogenized average solution using the simplified multicontinuum model 2}
\end{subfigure}
\caption{Distributions of average pressure and displacements in $x_1$ and $x_2$ directions (from left to right) in $\Omega_2$ at the final time on the coarse grid $20 \times 20$ for Example 1.}
\label{fig:coarse_results_2_ex_1}
\end{figure}

Figures \ref{fig:coarse_results_1_ex_1}-\ref{fig:coarse_results_2_ex_1} present the reference and homogenized average solutions at the final time in $\Omega_1$ and $\Omega_2$, respectively, using the coarse grid $20 \times 20$. From left to right, we depict distributions of average pressure and displacements in $x_1$ and $x_2$ directions. From top to bottom, we present the reference average solution and the homogenized average solutions using the full model and the second simplified model. The results of the first simplified model are similar to the full model. One can see that all the solutions obtained are similar, indicating high accuracy of our proposed multicontinuum models. However, one can notice the differences in distributions of average displacements in the $x_1$ direction obtained using the second simplified model. We can also see that the second simplified model overestimates average displacements in the $x_1$ direction.

Regarding the modeled process, we observe the fluid diffusing in the vertical directions. At the same time, the average pressure in $\Omega_2$ diffuses faster than the average pressure in $\Omega_1$, as expected. The displacements in the $x_1$ direction describe horizontal stretching. The displacements in the $x_2$ direction describe an upward stretching with a slightly different distribution and larger values in $\Omega_2$ due to softer elastic properties. In general, the results obtained are in agreement with those on the fine grid.

Let us consider the errors of our proposed multicontinuum models. Table \ref{tab:errors_ex_1} presents the relative $L_2$ errors for pressures and displacements at the final time using different multicontinuum models and coarse grids. One can see that all the errors are small. However, the second simplified model has more significant errors for displacements, which corresponds to Figures \ref{fig:coarse_results_1_ex_1}-\ref{fig:coarse_results_2_ex_1}. In contrast, the first simplified model provides the same accuracy as the full model.

\begin{table}[hbt!]
\centering
\begin{tabular}{c|c|c|c|c}
Multicontinuum model & $e^{(1)}_p$ & $e^{(2)}_p$ & $e^{(1)}_u$ & $e^{(2)}_u$\\ \hline
\multicolumn{5}{c}{Coarse grid $10 \times 10$} \\ \hline
Full         & 3.89e-03 & 3.46e-03 & 2.32e-02 & 2.32e-02 \\
Simplified 1 & 3.89e-03 & 3.46e-03 & 2.32e-02 & 2.32e-02 \\
Simplified 2 & 2.20e-03 & 1.86e-03 & 5.53e-02 & 5.53e-02 \\ \hline
\multicolumn{5}{c}{Coarse grid $20 \times 20$} \\ \hline
Full         & 9.54e-04 & 8.33e-04 & 1.04e-02 & 1.04e-02 \\
Simplified 1 & 9.54e-04 & 8.33e-04 & 1.04e-02 & 1.04e-02 \\
Simplified 2 & 9.71e-04 & 9.06e-04 & 5.84e-02 & 5.84e-02 \\
\end{tabular}
\caption{Relative $L_2$ errors for different multicontinuum models and coarse grids. Example 1.}
\label{tab:errors_ex_1}
\end{table}

Therefore, our proposed multicontinuum models can provide high accuracy for the anisotropic medium case. While the full and first simplified models have smaller errors compared to the second simplified model.

\subsection{Example 2}

\begin{figure}[hbt!]
\centering
\includegraphics[width=0.3\textwidth]{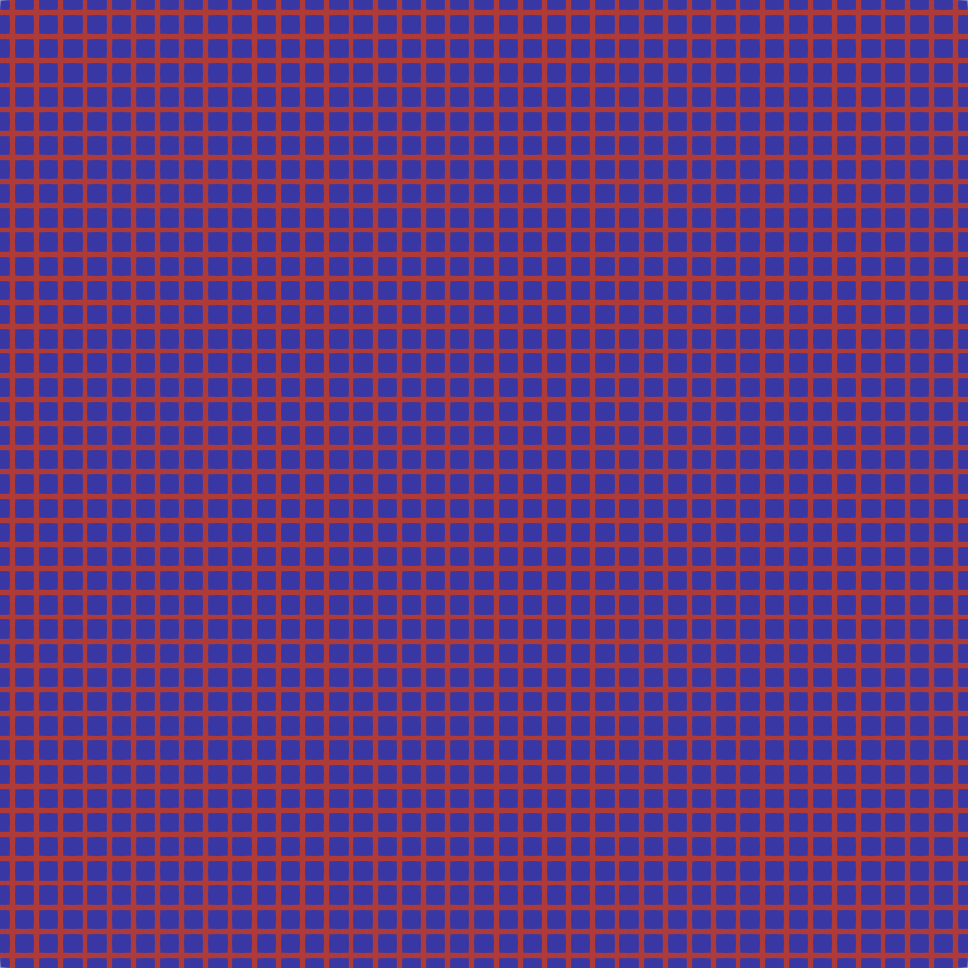}
\caption{Microstructure ($\Omega_1$ is blue, $\Omega_2$ is red) for Example 2.}
\label{fig:microstructure_ex_2}
\end{figure}

In this example, we consider an isotropic high-contrast heterogeneous medium. We depict the microstructure in Figure \ref{fig:microstructure_ex_2}, where $\Omega_1$ is blue, and $\Omega_2$ is red. For our numerical test, we set the following heterogeneous Lam\'e parameters $\lambda$ and $\mu$ and the permeability coefficient $\kappa$
\begin{equation}\label{eq:material_parameters_ex_2}
\lambda = \begin{cases}
    10^{5}, \quad x \in \Omega_1,\\
    10^{9}, \quad x \in \Omega_2
\end{cases}, \quad
\mu = \begin{cases}
    10^{5}, \quad x \in \Omega_1,\\
    10^{9}, \quad x \in \Omega_2
\end{cases}, \quad
\kappa = \begin{cases}
    10^{-10}, \quad x \in \Omega_1,\\
    10^{-6}, \quad x \in \Omega_2
\end{cases}.
\end{equation}
Therefore, we switch the values of $\lambda$ and $\mu$ in $\Omega_1$ and $\Omega_2$ compared to the previous example.

We use the same Biot coefficient $\alpha = 0.8$ and Biot modulus $M = 10^6$. For the time duration, we set $t_{max} = 1$ with 50 time steps. We set zero right-hand sides $f_1 = 0$, $f_2 = 0$, and $g = 0$. We use the same $u_0 = (0, 0)$ and $p_0 = 10^6$ for initial conditions. For boundary conditions, we fix displacements in $x_1$ on the left and right boundaries and $x_2$ on the bottom boundary and set free traction on the other boundaries. We set input and output pressures on the top and bottom boundaries, respectively, and no-flow on the other boundaries. We can write the boundary conditions in the following way
\begin{equation}\label{eq:bcs_ex_2}
\begin{split}
&u_1 = 0, \quad (\sigma \cdot \nu)_2 = 0, \quad x \in (\Gamma_L \cup \Gamma_R),\\
&u_2 = 0, \quad (\sigma \cdot \nu)_1 = 0, \quad x \in \Gamma_B,\\
&\sigma \cdot \nu = 0, \quad x \in \Gamma_T,\\
&p = h, \quad x \in \Gamma_T, \quad p = p_0, \quad x \in \Gamma_B,\\
&-\kappa \frac{\partial p}{\partial \nu} = 0, \quad x \in \Gamma_L \cup \Gamma_R,
\end{split}
\end{equation}
where $h = 1.6 \,p_0 + 0.4 \,p_0 \cos(2 \pi x_1)$. Therefore, we simulate fluid injection from the top boundary, and displacements are forced by pressure gradient.

\begin{figure}[hbt!]
\centering
\includegraphics[height=0.24\textheight]{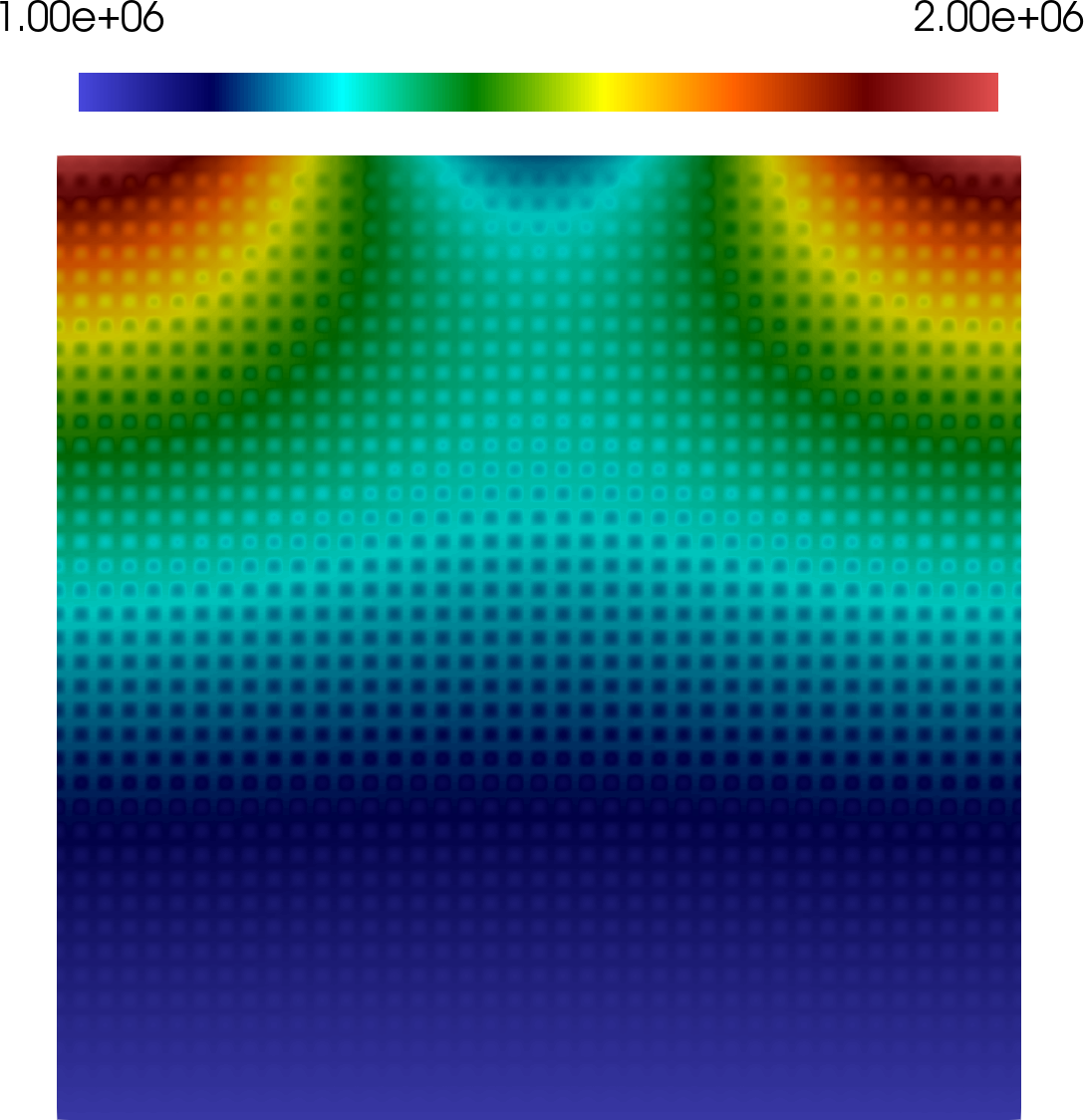}
\includegraphics[height=0.24\textheight]{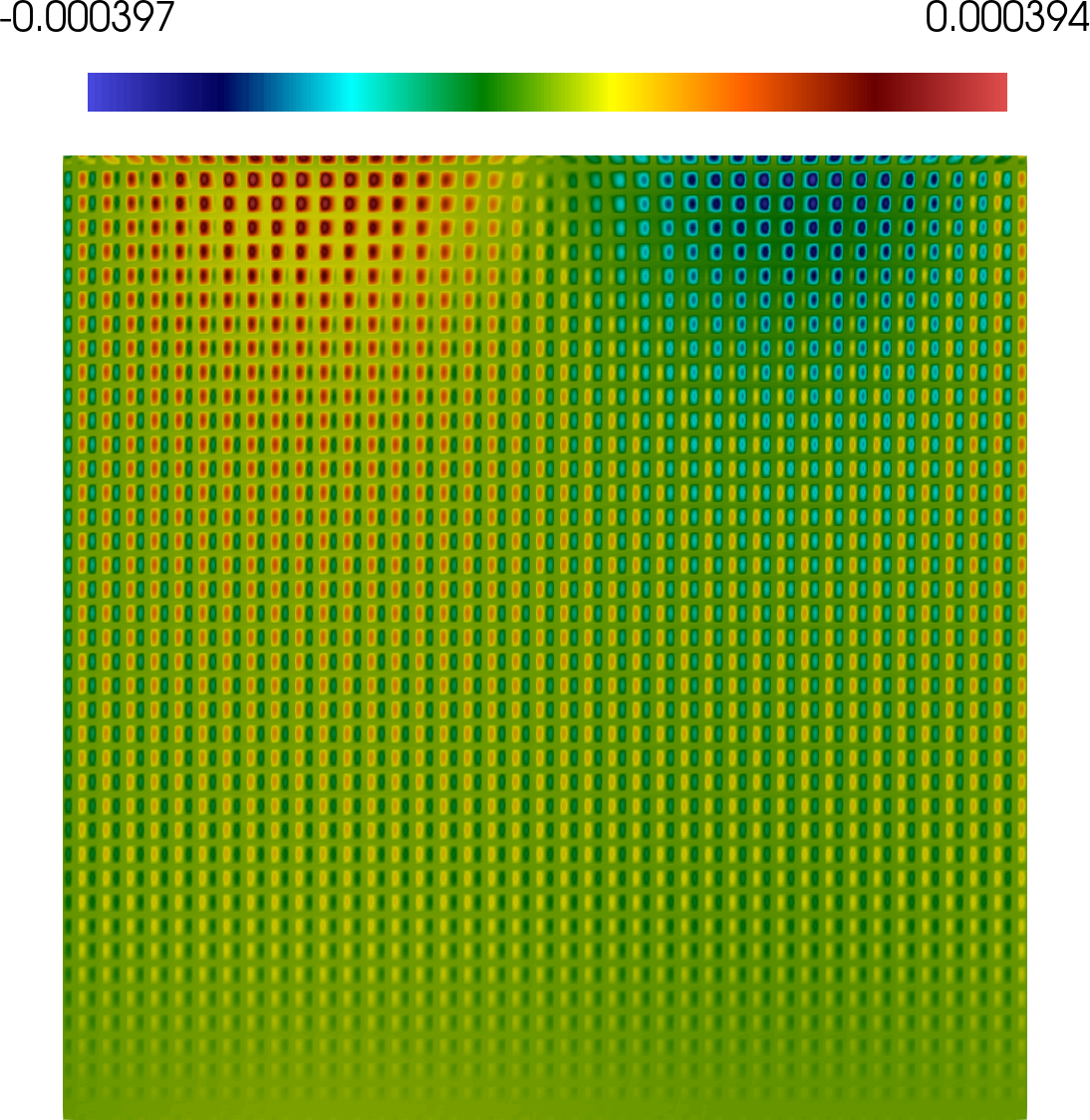}
\includegraphics[height=0.24\textheight]{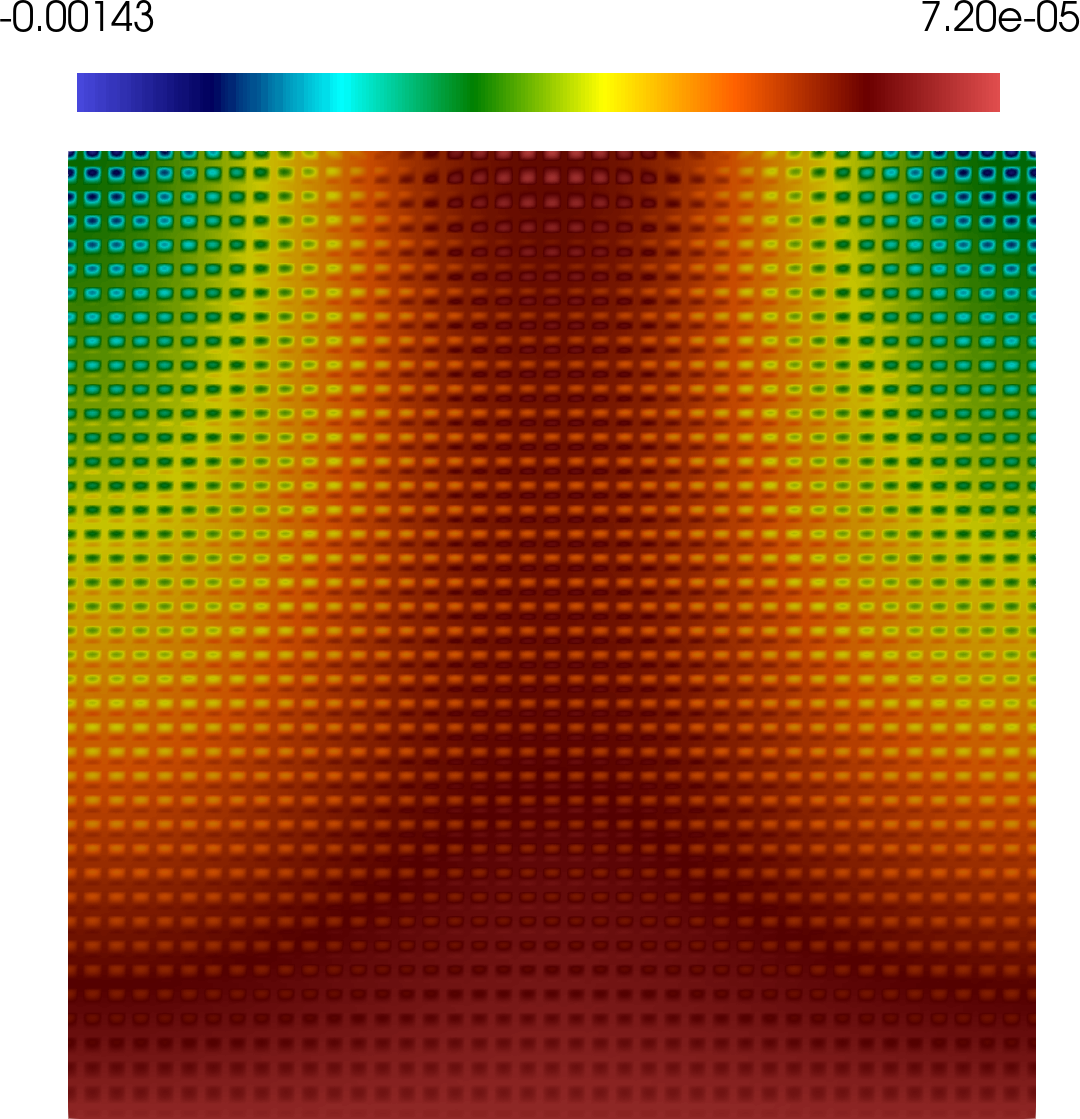}
\caption{Distributions of pressure and displacements in $x_1$ and $x_2$ directions (from left to right) at the final time on the fine grid for Example 2.}
\label{fig:fine_results_ex_2}
\end{figure}

Figure \ref{fig:fine_results_ex_2} presents distributions of pressure, displacements in $x_1$ and $x_2$ directions (from left to right) at the final time on the fine grid. We observe the fluid injected from the top boundary with different intensities; the highest pressure is in the corners. The fluid diffuses through the domain in an isotropic manner. One can notice the influence of the microstructure on the diffusion process. We also observe deformations due to the fluid injection, mainly on the upper part of the domain. The displacement distributions are also heterogeneous due to the influence of the microstructure. The obtained solution has an isotropic nature and corresponds to the modeled process.

\begin{figure}[hbt!]
\centering
\begin{subfigure}{\textwidth}
\centering
\includegraphics[height=0.24\textheight]{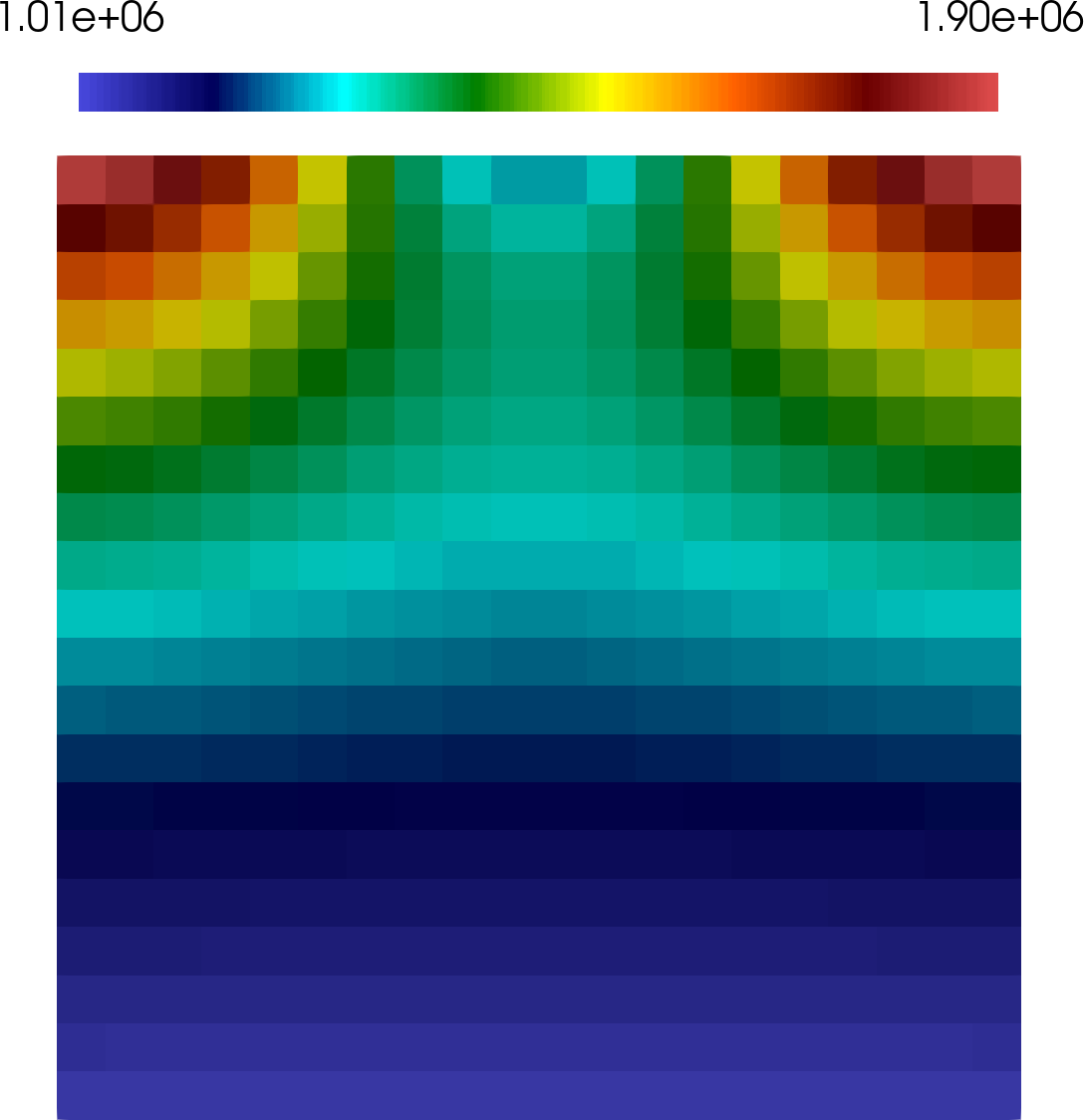}
\includegraphics[height=0.24\textheight]{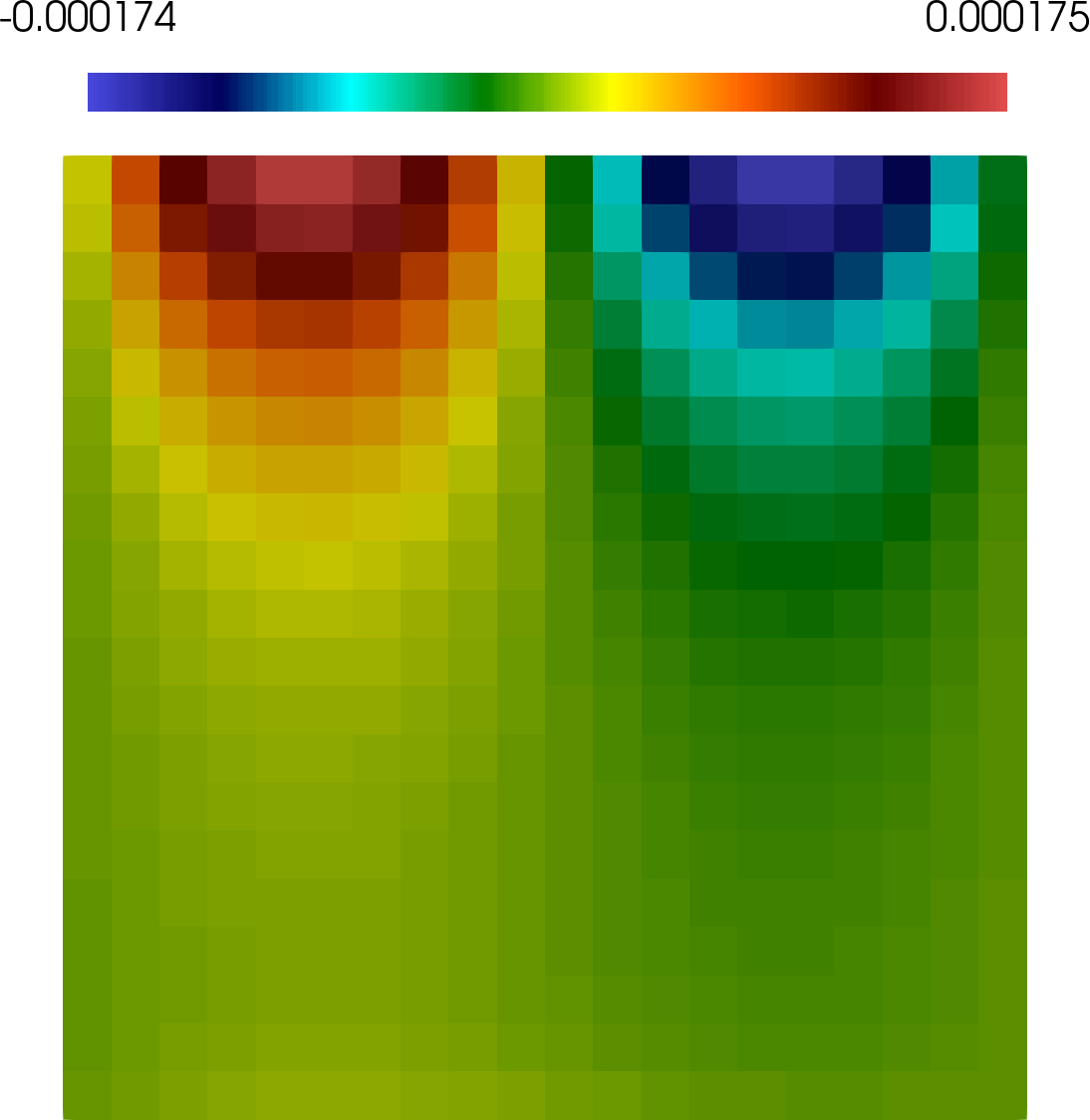}
\includegraphics[height=0.24\textheight]{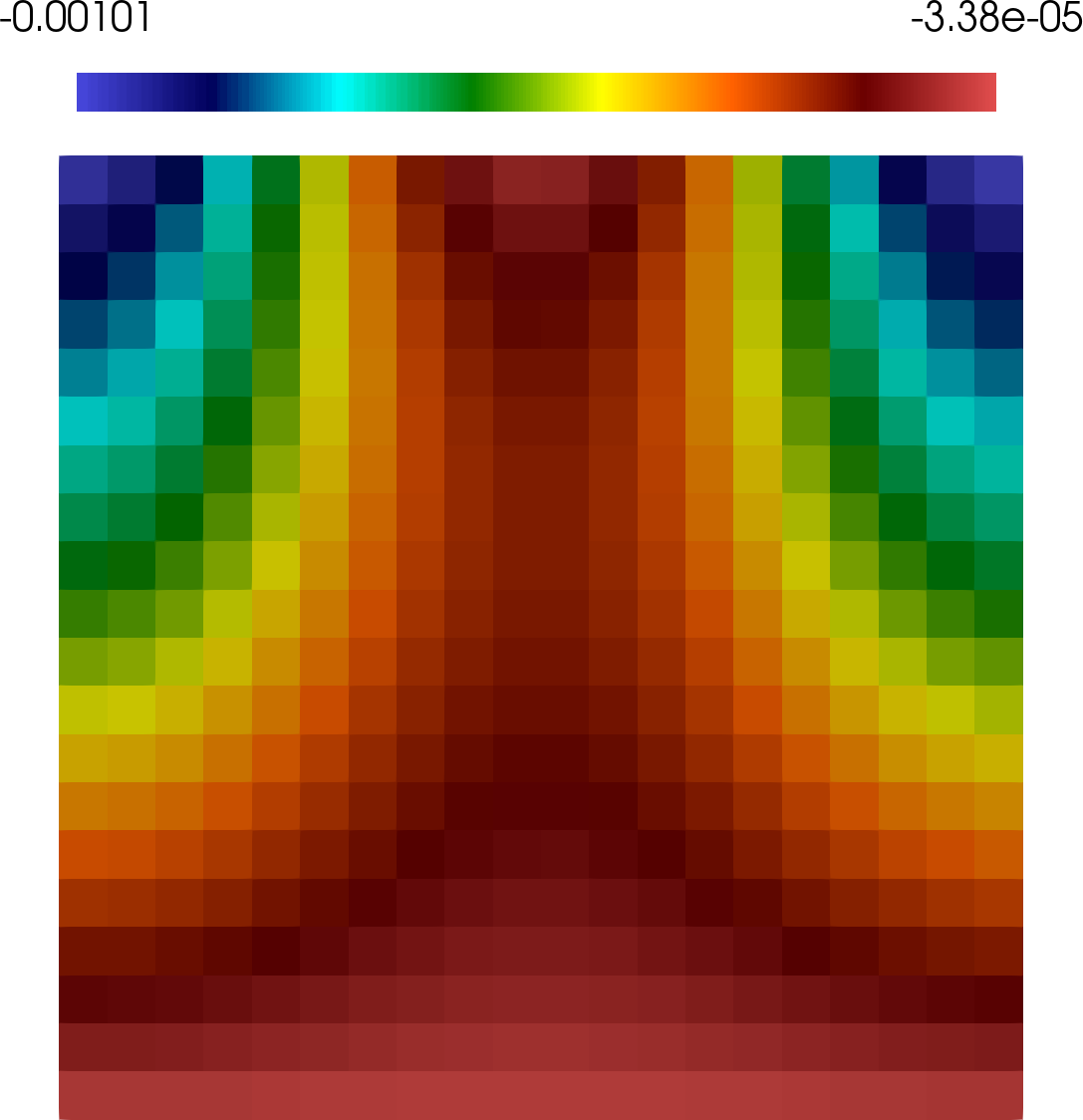}
\caption{Reference average solution}
\end{subfigure}
\begin{subfigure}{\textwidth}
\centering
\includegraphics[height=0.24\textheight]{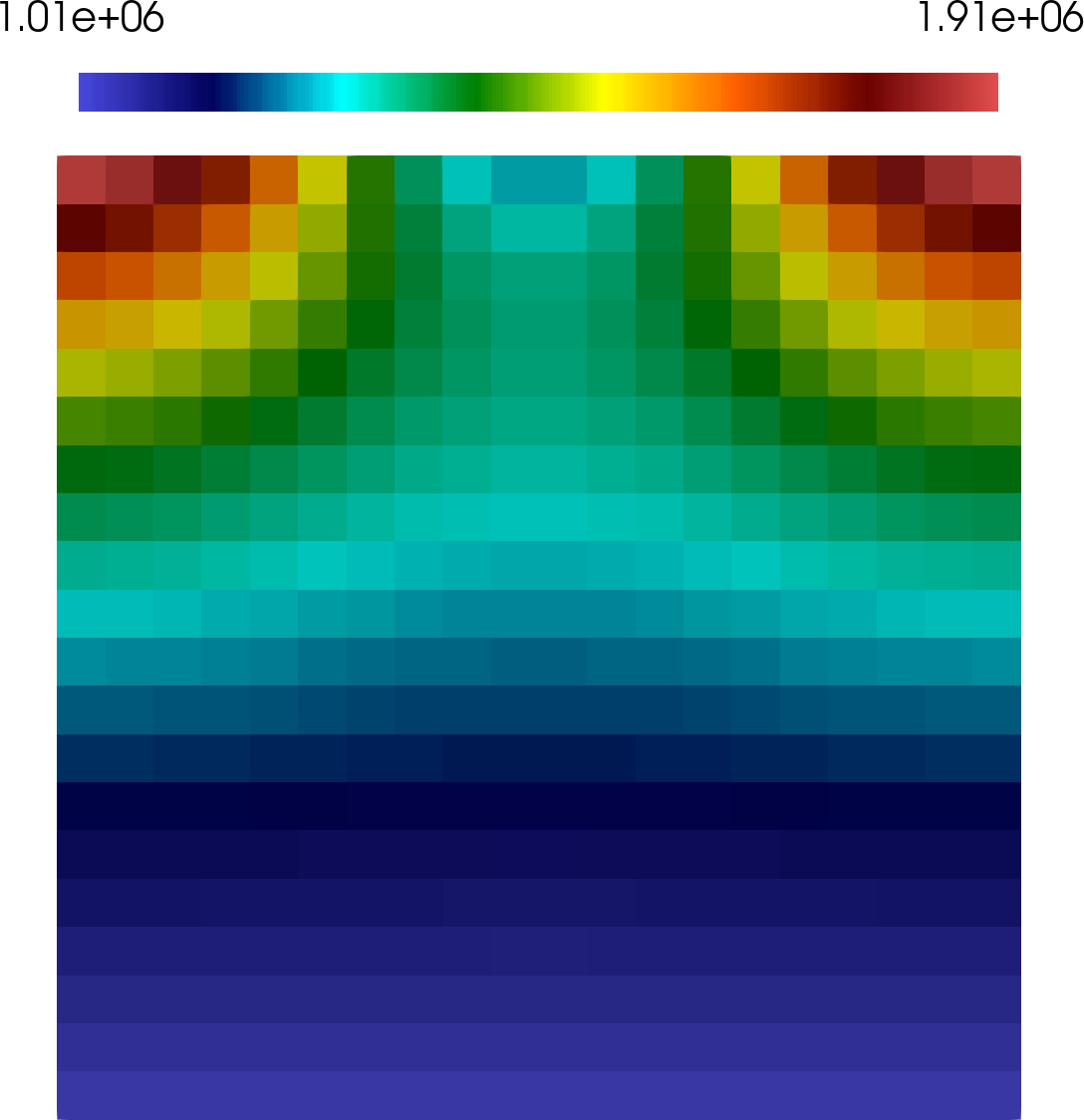}
\includegraphics[height=0.24\textheight]{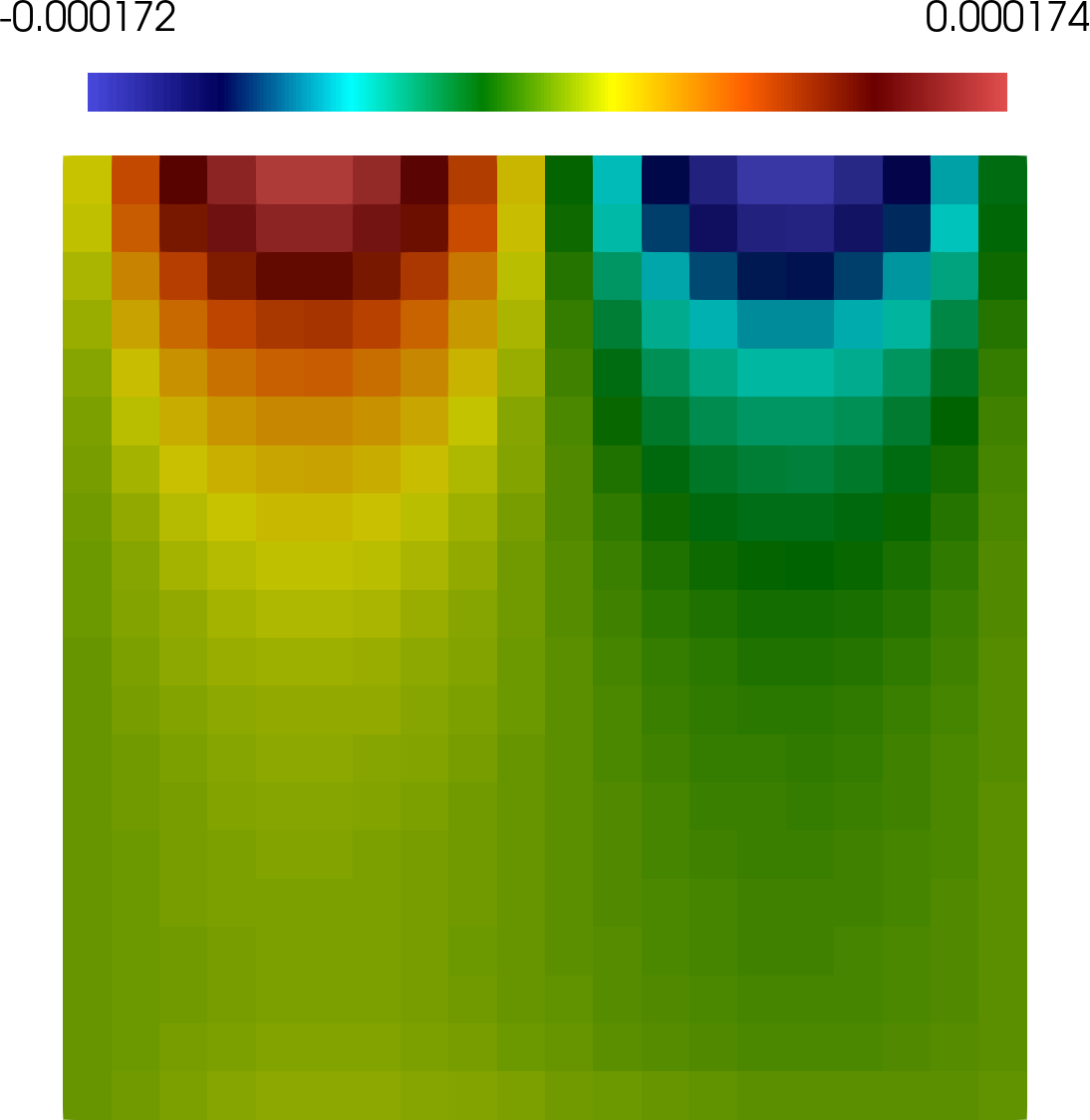}
\includegraphics[height=0.24\textheight]{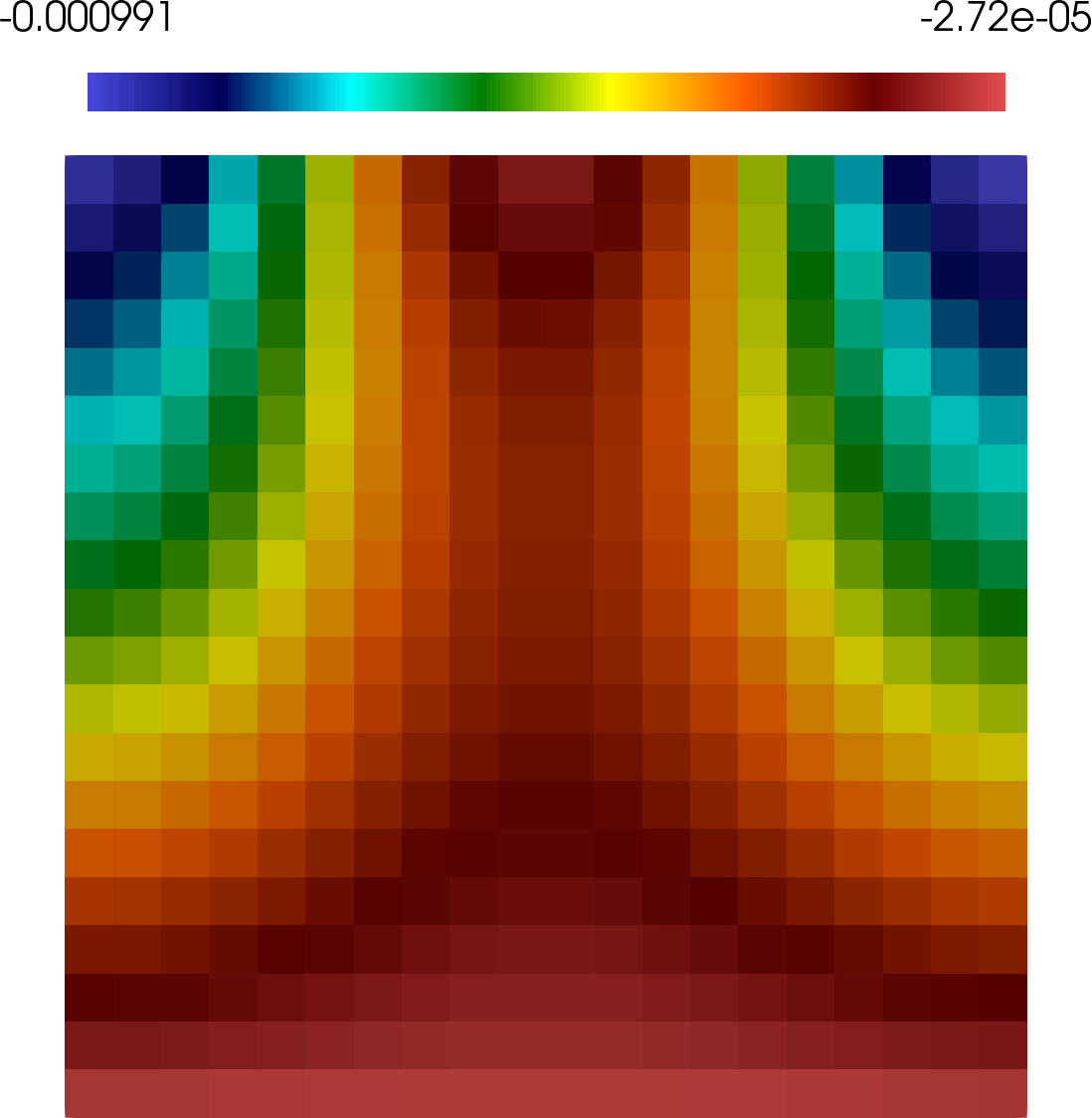}
\caption{Homogenized average solution using the full multicontinuum model}
\end{subfigure}
\begin{subfigure}{\textwidth}
\centering
\includegraphics[height=0.24\textheight]{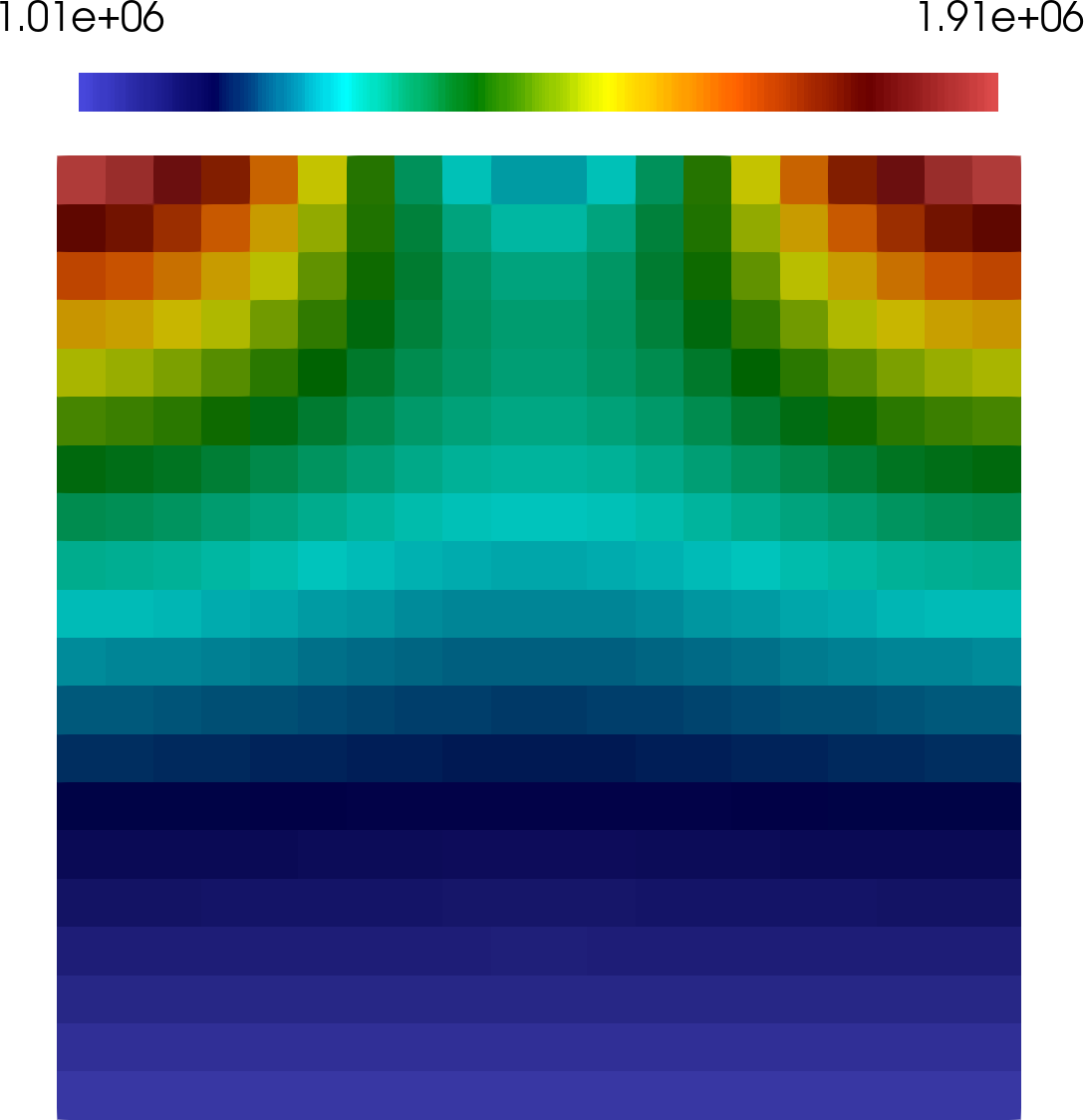}
\includegraphics[height=0.24\textheight]{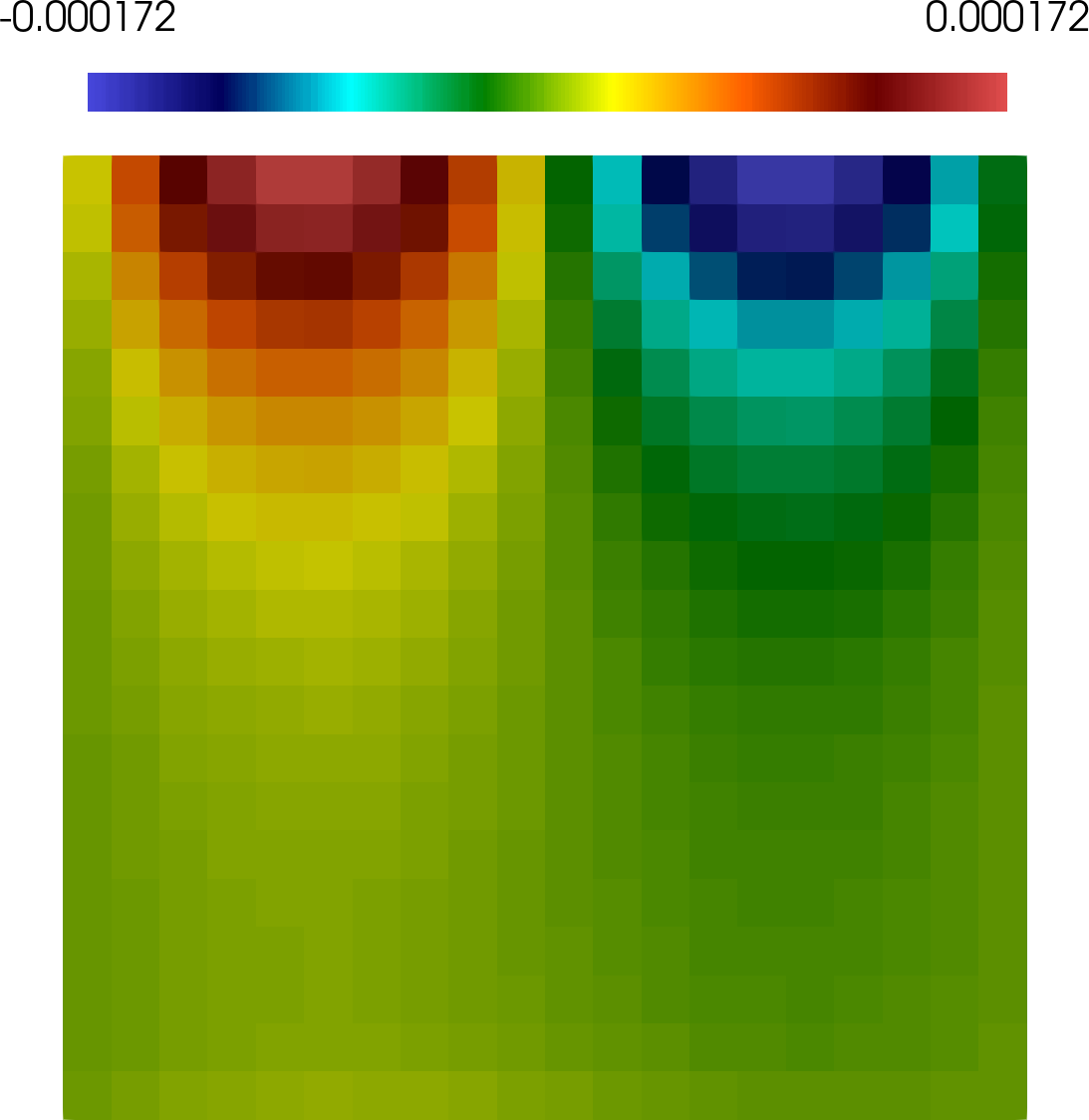}
\includegraphics[height=0.24\textheight]{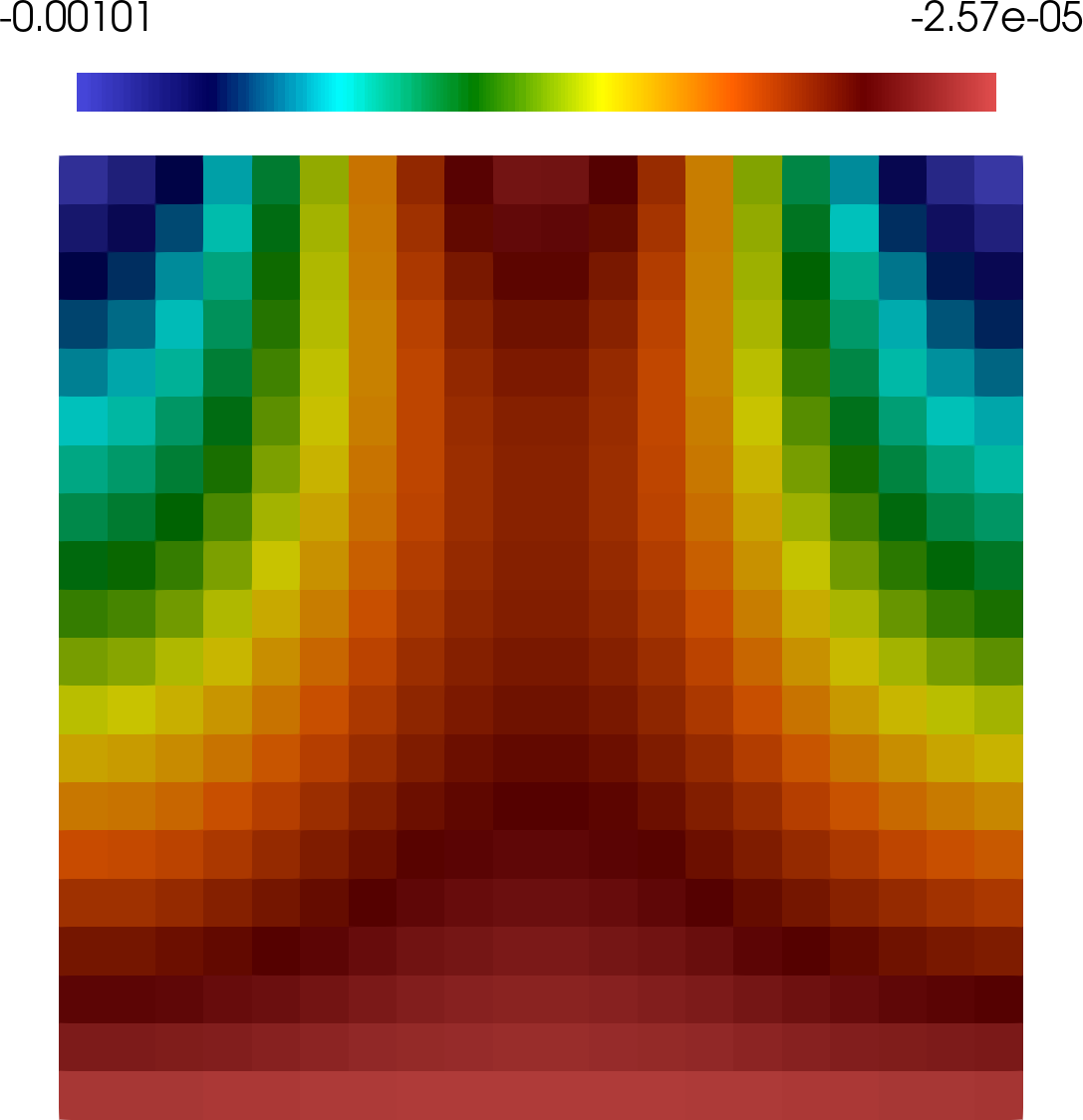}
\caption{Homogenized average solution using the simplified multicontinuum model 2}
\end{subfigure}
\caption{Distributions of average pressure and displacements in $x_1$ and $x_2$ directions (from left to right) in $\Omega_1$ at the final time on the coarse grid $20 \times 20$ for Example 2.}
\label{fig:coarse_results_1_ex_2}
\end{figure}

\begin{figure}[hbt!]
\centering
\begin{subfigure}{\textwidth}
\centering
\includegraphics[height=0.24\textheight]{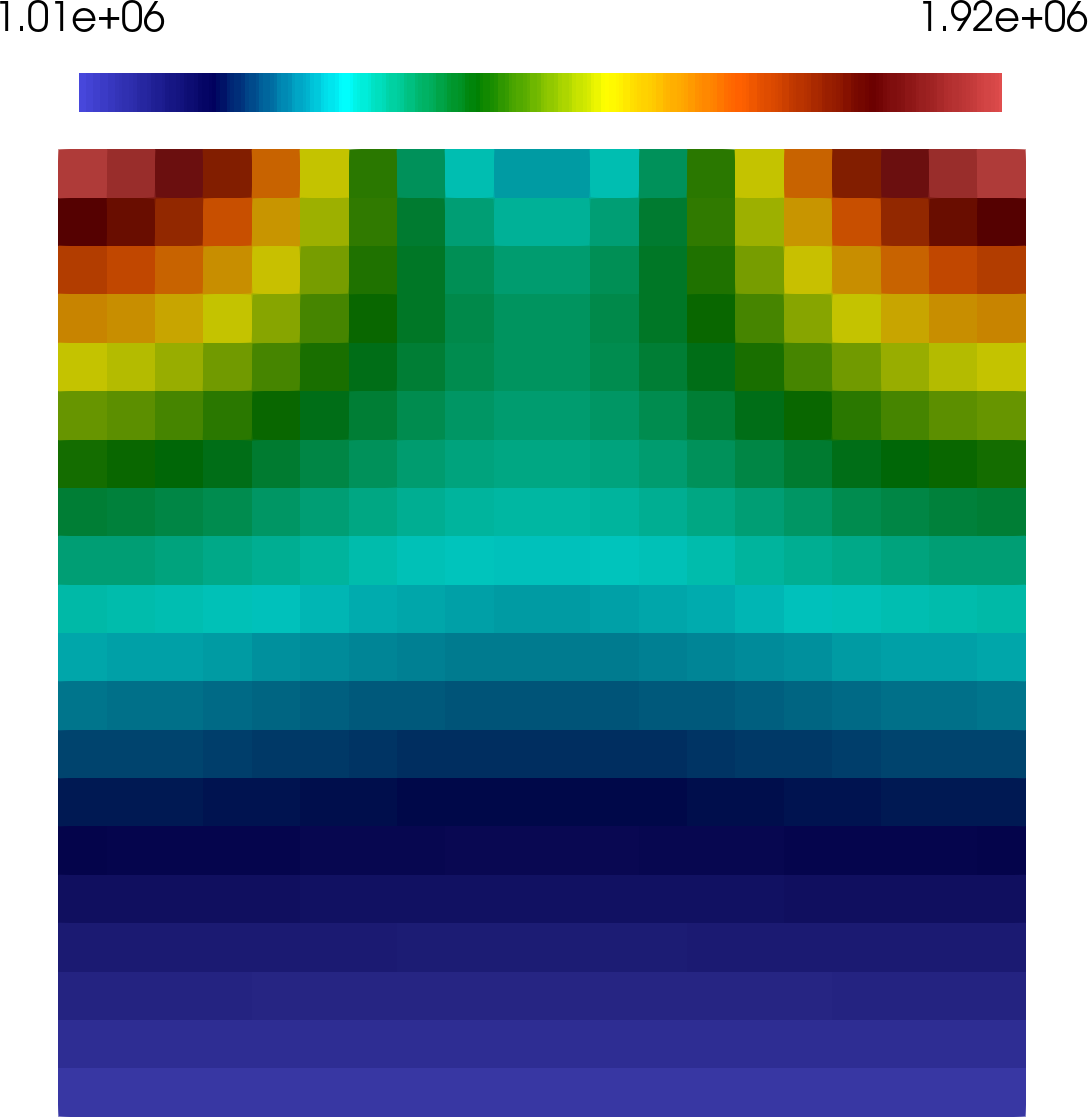}
\includegraphics[height=0.24\textheight]{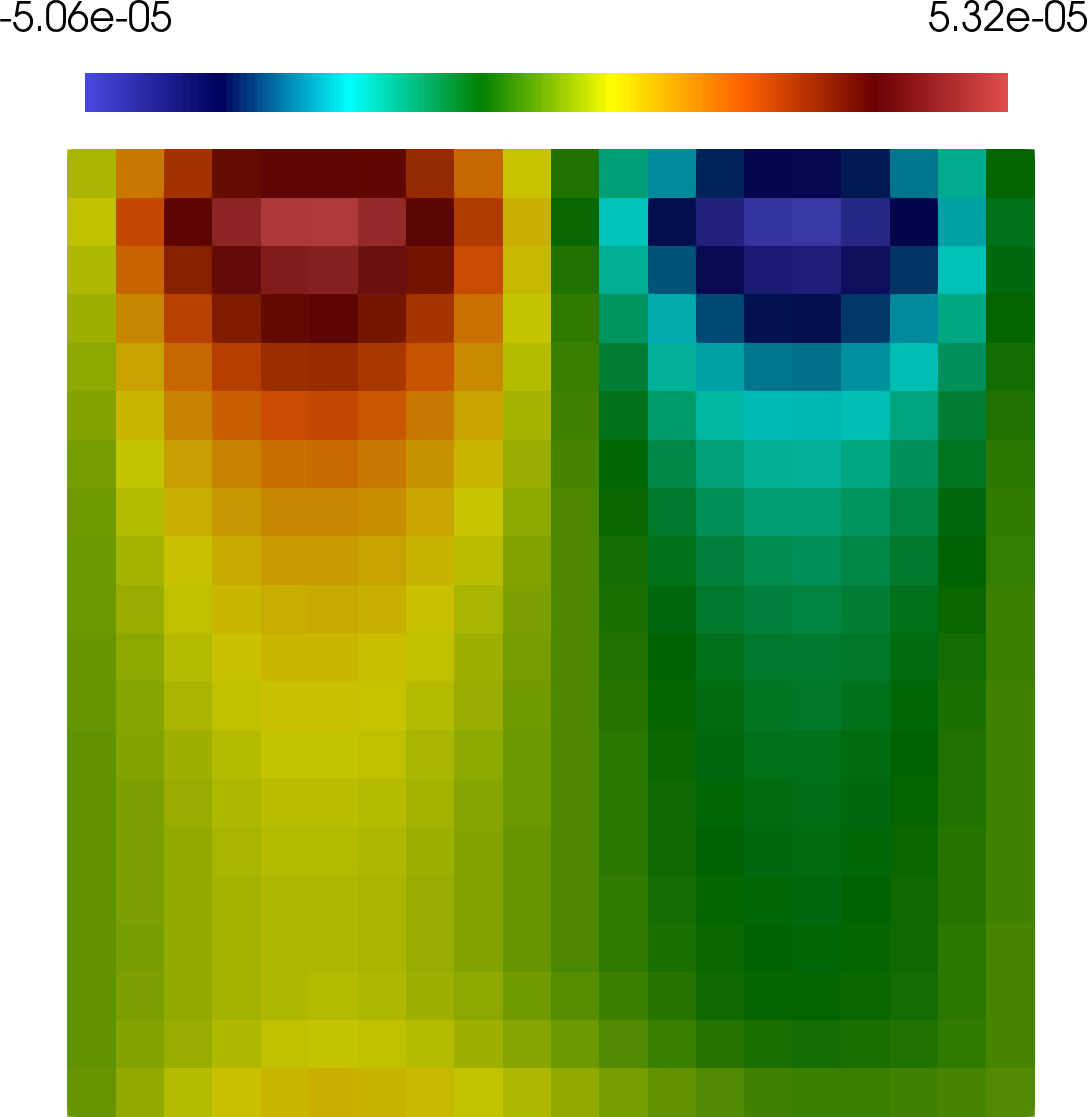}
\includegraphics[height=0.24\textheight]{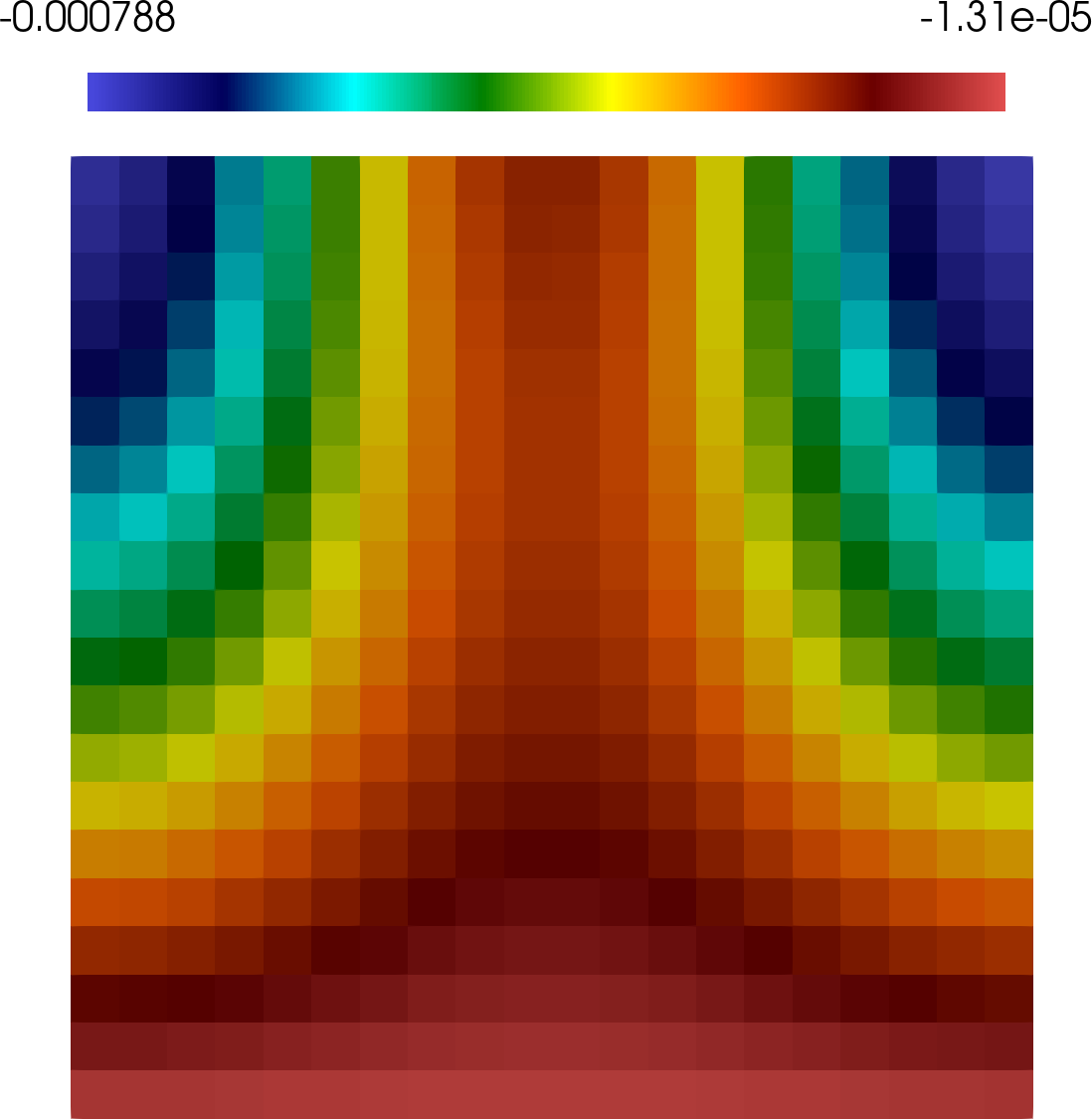}
\caption{Reference average solution}
\end{subfigure}
\begin{subfigure}{\textwidth}
\centering
\includegraphics[height=0.24\textheight]{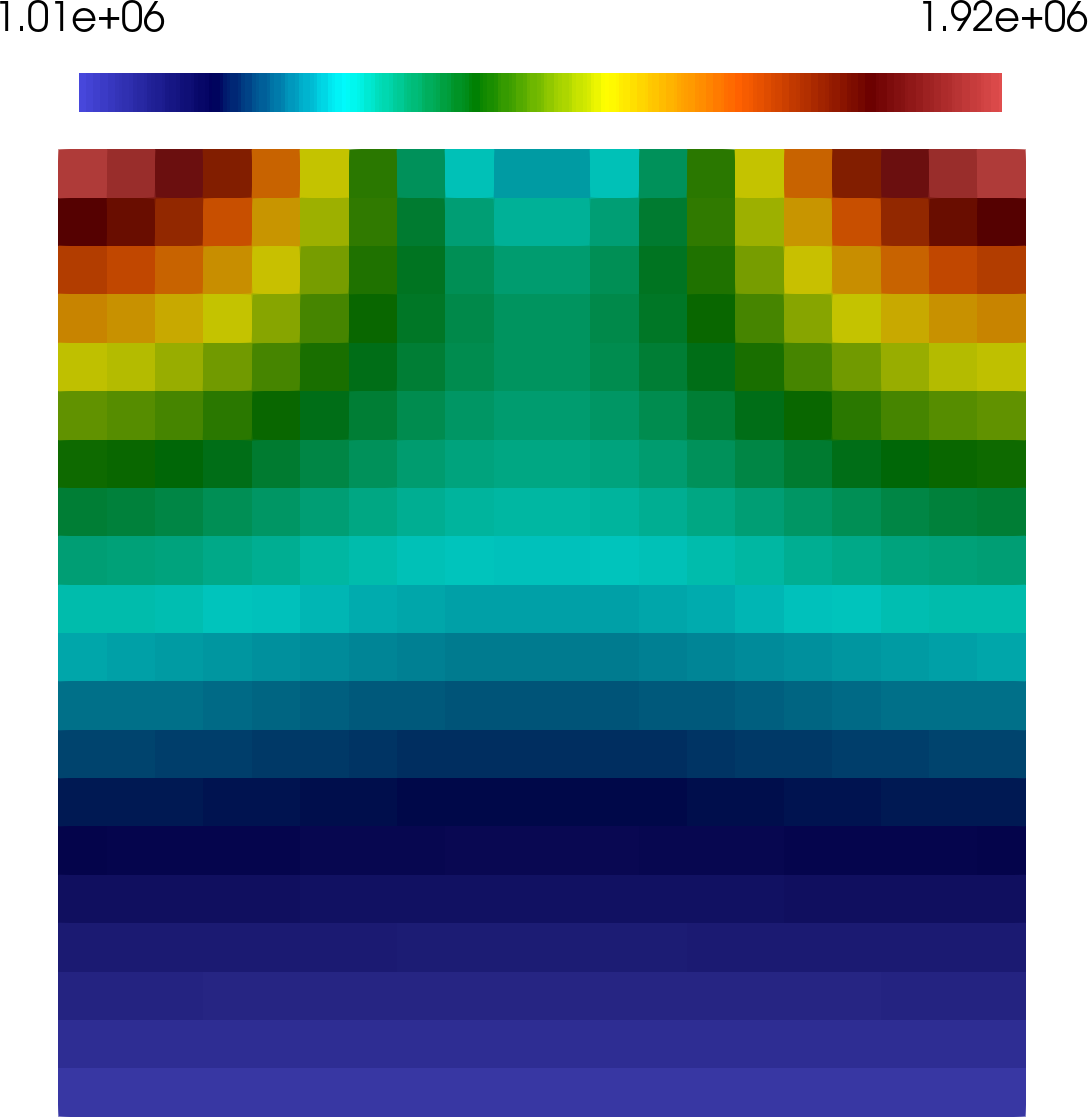}
\includegraphics[height=0.24\textheight]{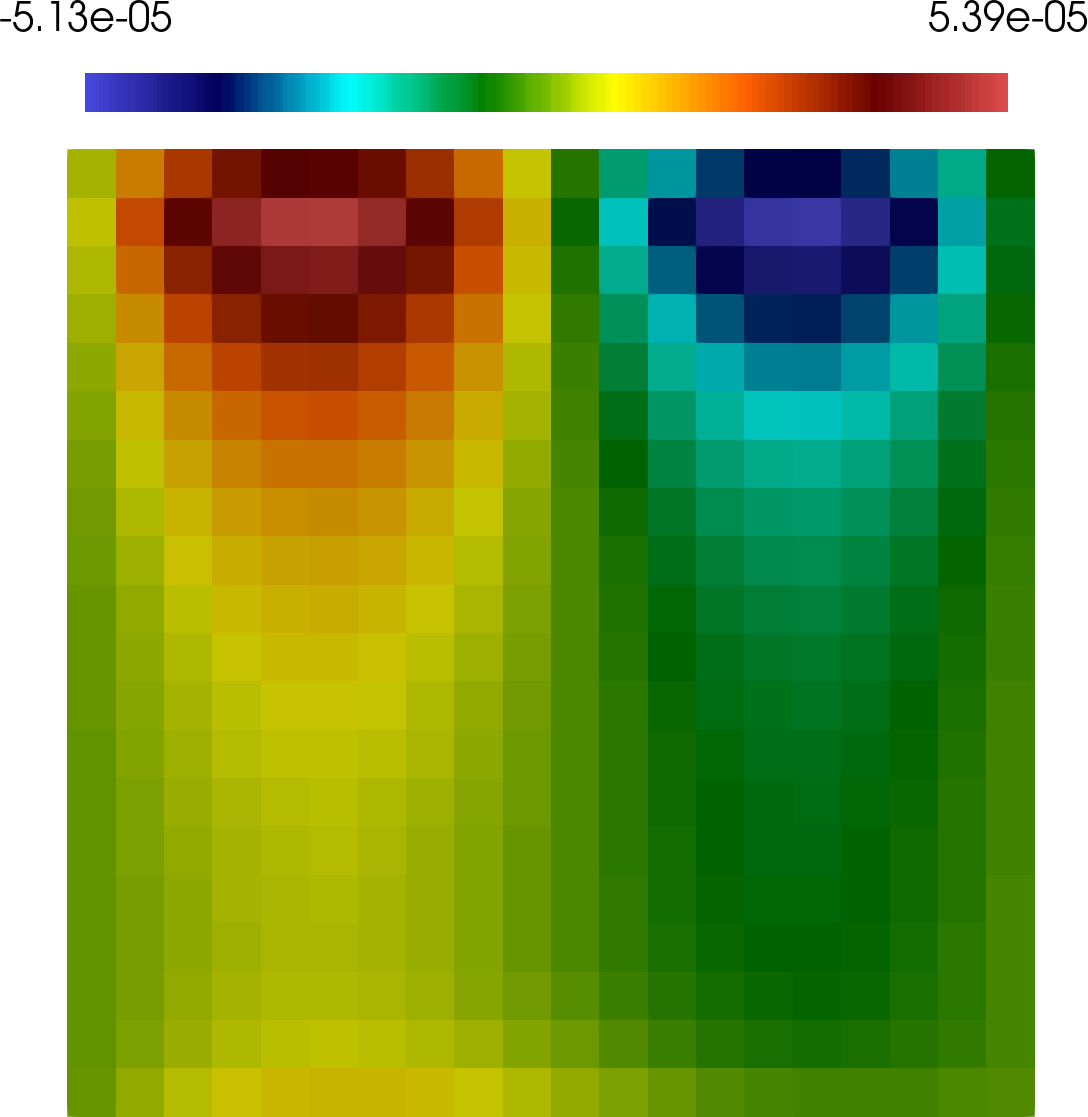}
\includegraphics[height=0.24\textheight]{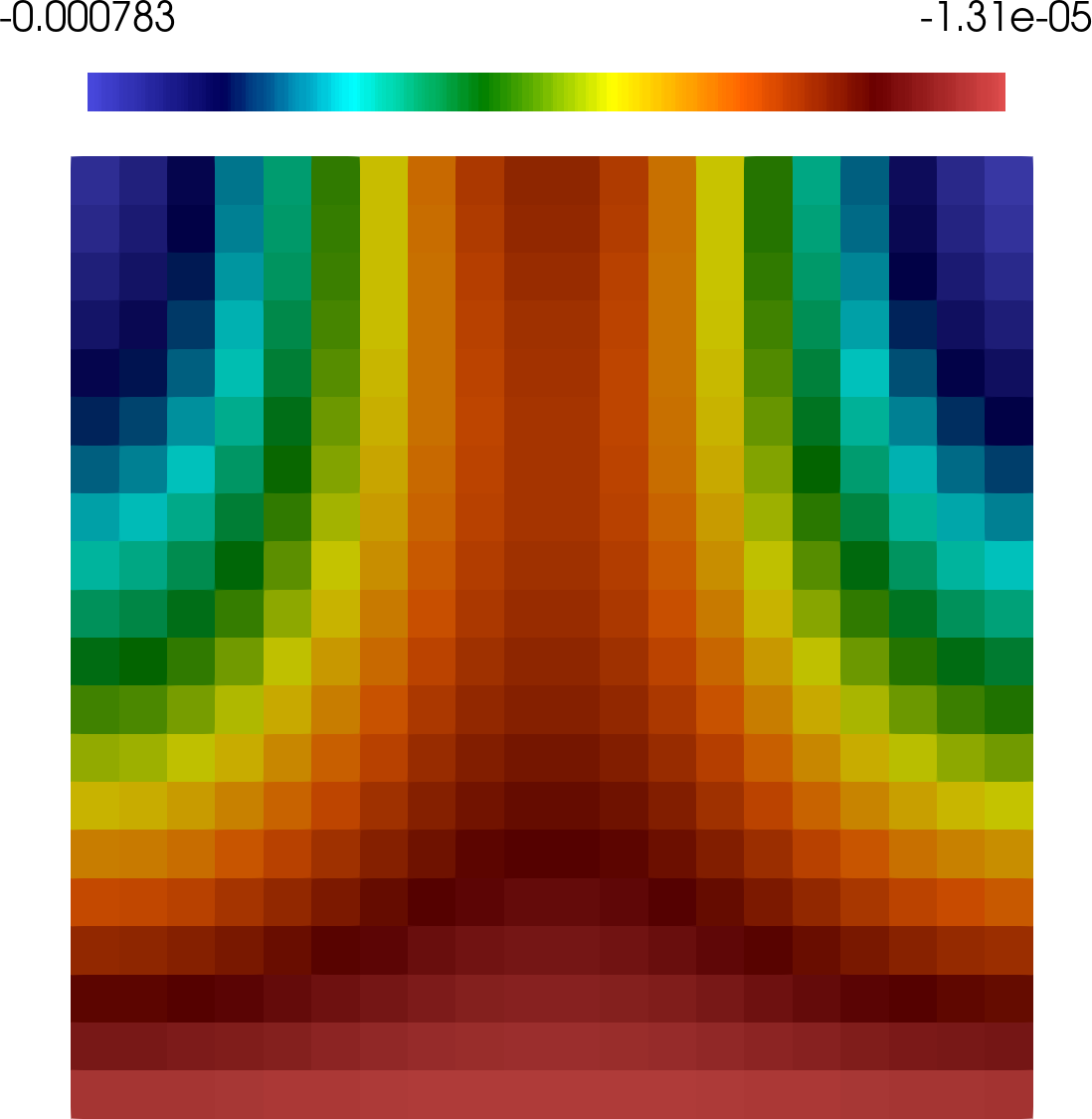}
\caption{Homogenized average solution using the full multicontinuum model}
\end{subfigure}
\begin{subfigure}{\textwidth}
\centering
\includegraphics[height=0.24\textheight]{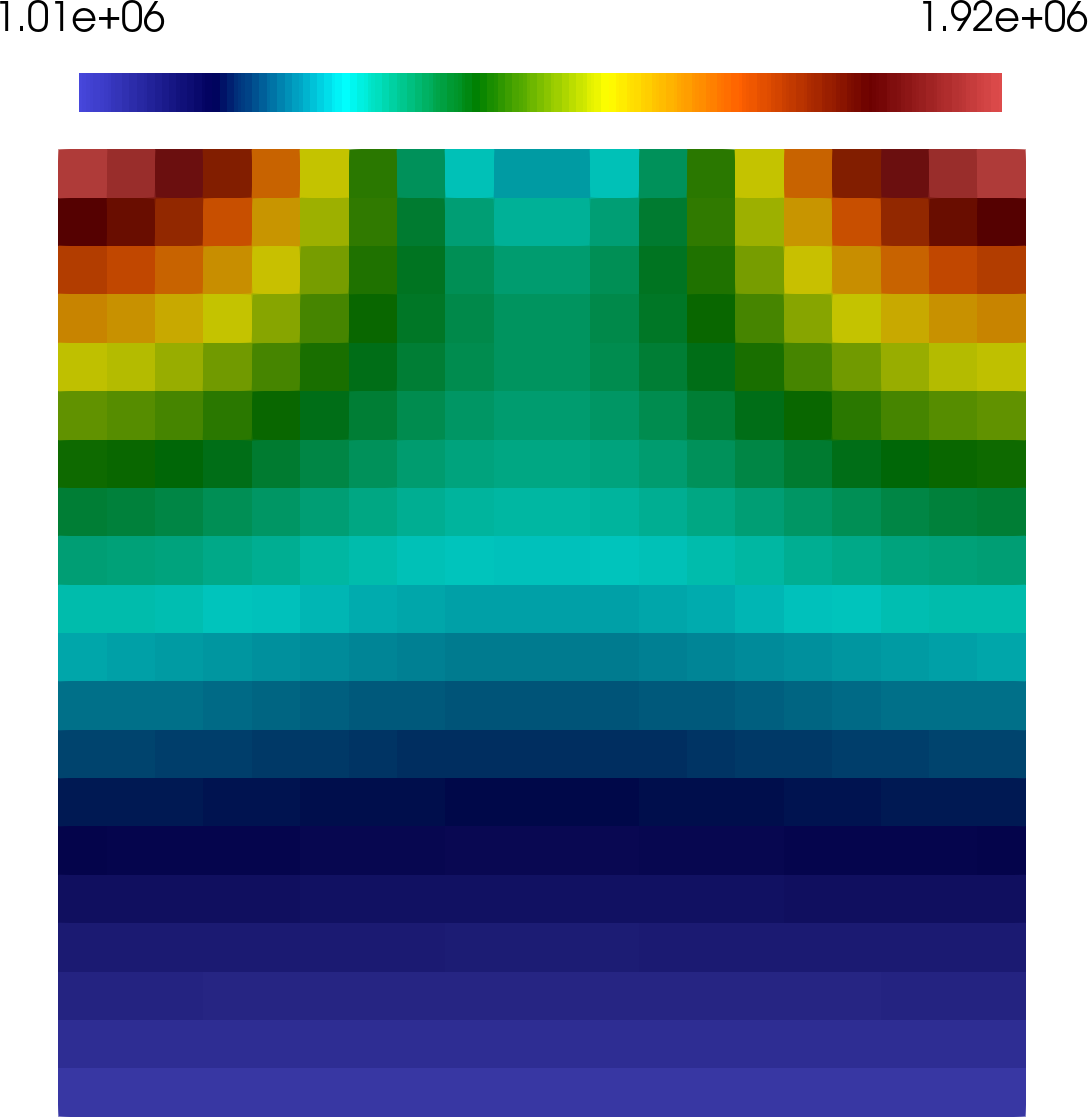}
\includegraphics[height=0.24\textheight]{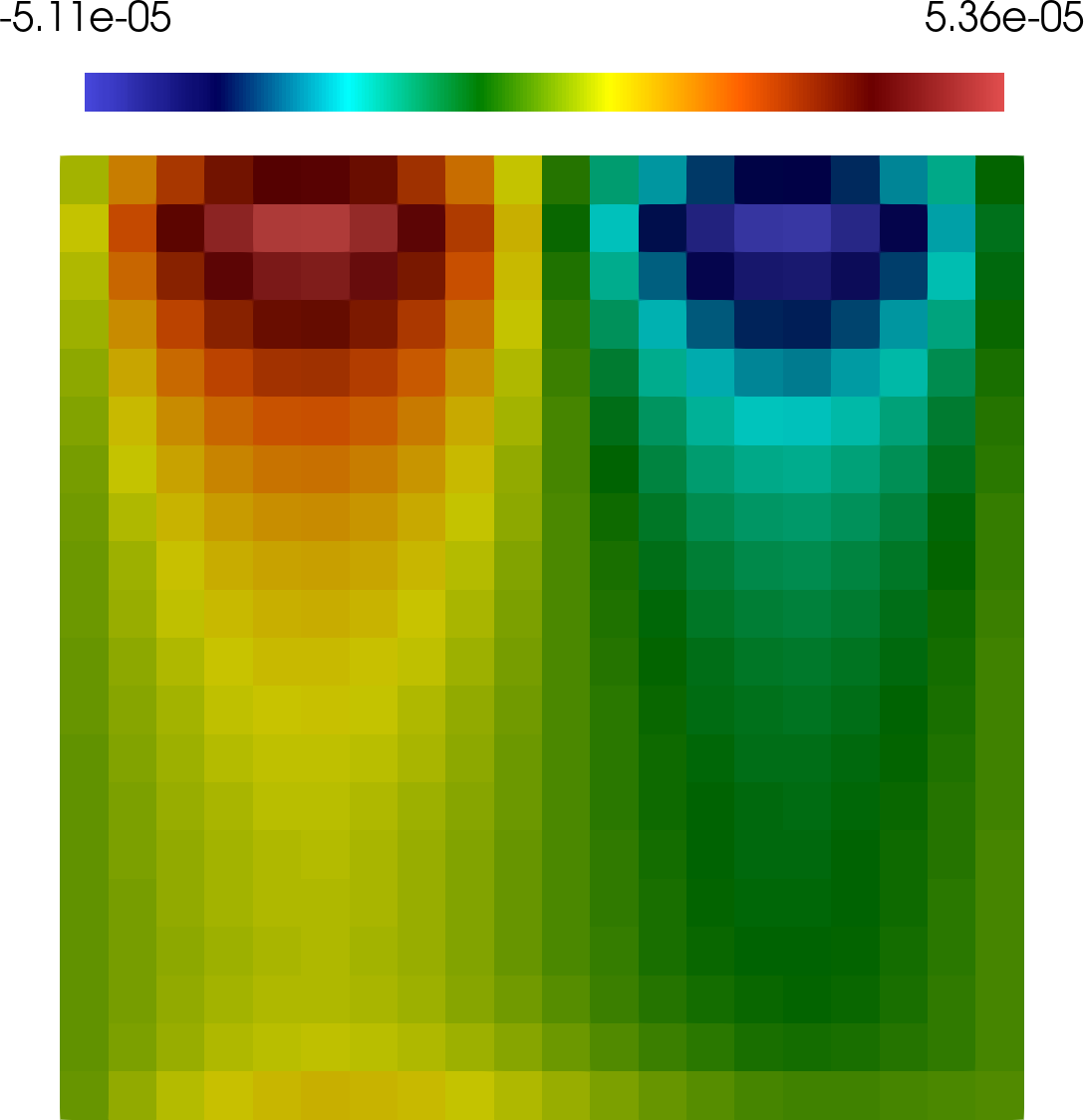}
\includegraphics[height=0.24\textheight]{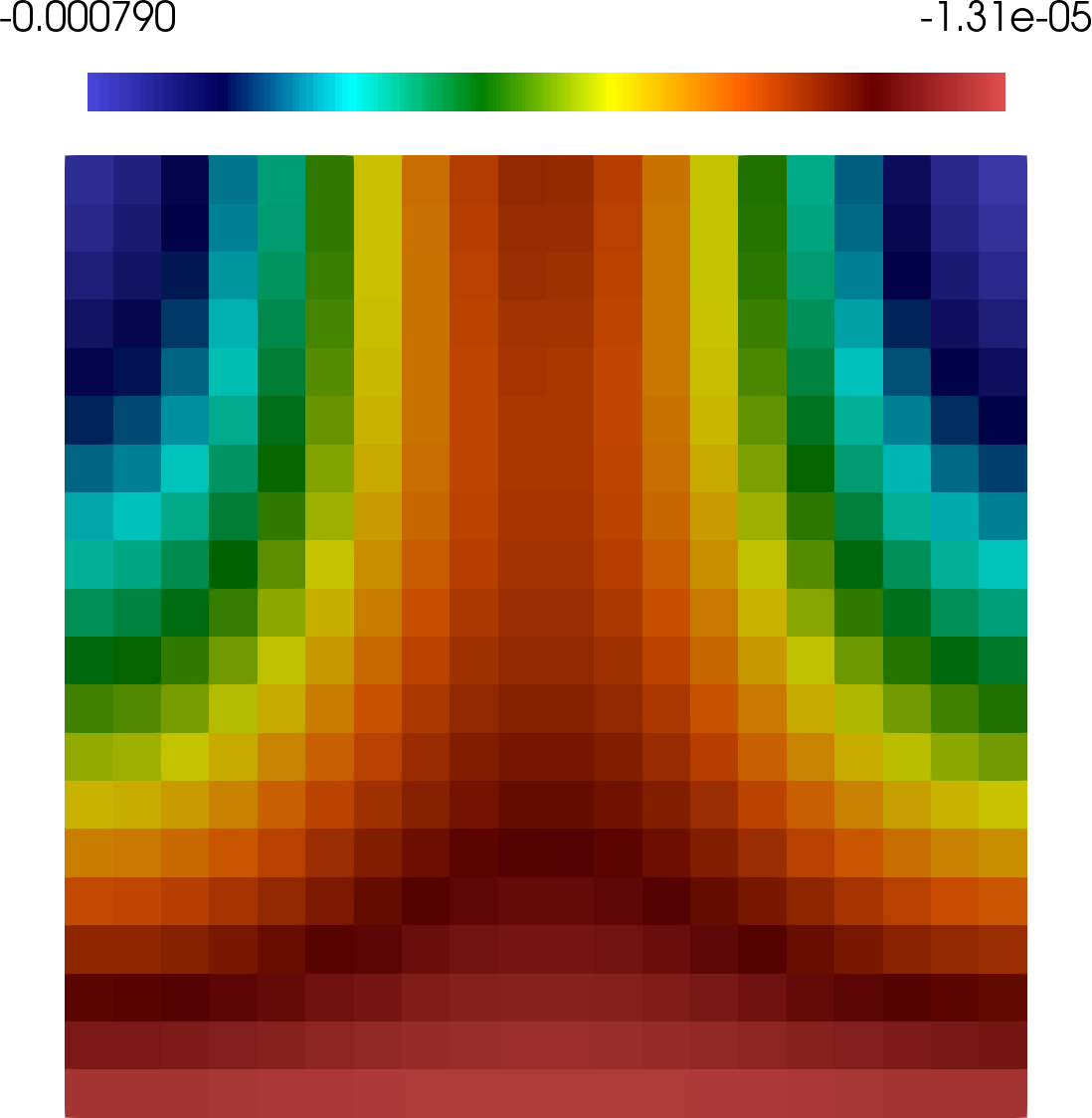}
\caption{Homogenized average solution using the simplified multicontinuum model 2}
\end{subfigure}
\caption{Distributions of average pressure and displacements in $x_1$ and $x_2$ directions (from left to right) in $\Omega_2$ at the final time on the coarse grid $20 \times 20$ for Example 2.}
\label{fig:coarse_results_2_ex_2}
\end{figure}

Figures \ref{fig:coarse_results_1_ex_2}-\ref{fig:coarse_results_2_ex_2} present the reference and homogenized average solutions at the final time in $\Omega_1$ and $\Omega_2$, respectively, using the coarse grid $20 \times 20$. We depict distributions of average pressure and displacements in the $x_1$ and $x_2$ directions from left to right. From top to bottom, we present the reference average solution and the homogenized average solutions using the coarse grid $20 \times 20$. One can see that the obtained solutions are very similar, which indicates the high accuracy of our proposed multicontinuum approaches. 

In terms of the simulated process, we observe fluid diffusion from the top boundary, especially from the top left and right corners. The average displacements in the $x_1$ direction are directed inward due to the influence of the fluid diffusion. At the same time, the average displacements in the $x_2$ direction are directed downward along the left and right boundaries. One can see that the average displacements in $\Omega_1$ are more significant than in $\Omega_2$ due to softer elastic properties.

Let us consider the errors of our proposed multicontinuum approaches. Table \ref{tab:errors_ex_2} presents the relative $L_2$ errors at the final time for different models and coarse grids. One can see that all the errors are minor, but the full and the first simplified models are more accurate.

\begin{table}[hbt!]
\centering
\begin{tabular}{c|c|c|c|c}
Multicontinuum model & $e^{(1)}_p$ & $e^{(2)}_p$ & $e^{(1)}_u$ & $e^{(2)}_u$\\ \hline
\multicolumn{5}{c}{Coarse grid $10 \times 10$} \\ \hline
Full         & 1.64e-03 & 1.36e-03 & 2.26e-02 & 1.91e-02 \\
Simplified 1 & 1.64e-03 & 1.36e-03 & 2.27e-02 & 1.91e-02 \\
Simplified 2 & 1.71e-03 & 1.36e-03 & 2.58e-02 & 2.09e-02 \\ \hline
\multicolumn{5}{c}{Coarse grid $20 \times 20$} \\ \hline
Full         & 7.32e-04 & 3.59e-04 & 8.33e-03 & 5.11e-03 \\
Simplified 1 & 7.32e-04 & 3.59e-04 & 8.33e-03 & 5.11e-03 \\
Simplified 2 & 8.45e-04 & 3.50e-04 & 1.62e-02 & 1.11e-02 \\
\end{tabular}
\caption{Relative $L_2$ errors for different multicontinuum models and coarse grids. Example 2.}
\label{tab:errors_ex_2}
\end{table}

\FloatBarrier
\subsection{Example 3}

\begin{figure}[hbt!]
\centering
\includegraphics[width=0.3\textwidth]{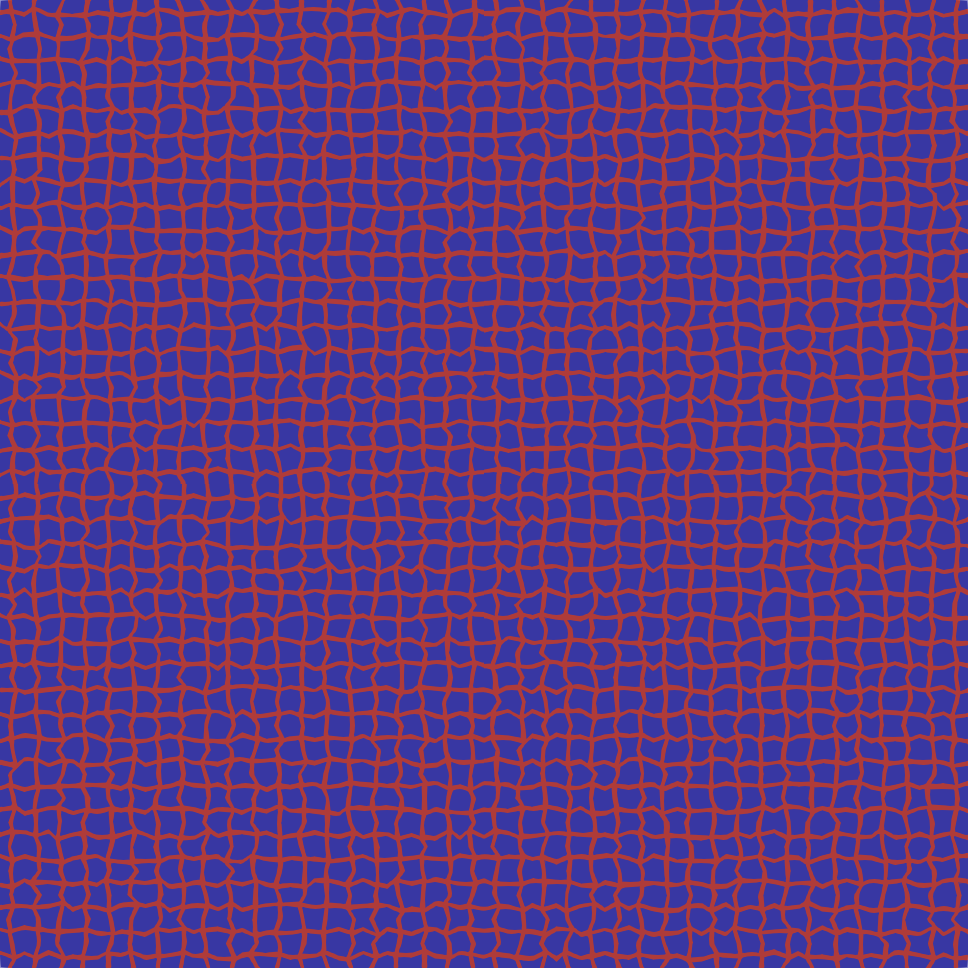}
\caption{Microstructure ($\Omega_1$ is blue, $\Omega_2$ is red) for Example 3.}
\label{fig:microstructure_ex_3}
\end{figure}

In the last representative example, we consider a non-periodic high-contrast heterogeneous medium. Figure \ref{fig:microstructure_ex_3} depicts the microstructure, where $\Omega_1$ is blue, and $\Omega_2$ is red. For heterogeneous coefficients, we set \eqref{eq:material_parameters_ex_2}. We apply the boundary conditions \eqref{eq:bcs_ex_1} and set the following right-hand sides
\begin{equation}
\begin{gathered}
f_{1} = -\sin\left(2\pi\, x_1\right)\sin\left(\pi\, x_2\right) \, 10^7, \quad
f_{2} = \sin\left(\pi\, x_1\right)\sin\left(\pi\, x_2\right) \, 10^7,\\
g = \exp(-40[(x_1-0.5)^2+(x_2-0.5)^2]).
\end{gathered}
\end{equation}
We simulate for $t_{max} = 5$ with 50 time steps. For initial conditions, we set $u_0 = (0, 0)$ and $p_0 = 10^6$.

\begin{figure}[hbt!]
\centering
\includegraphics[height=0.24\textheight]{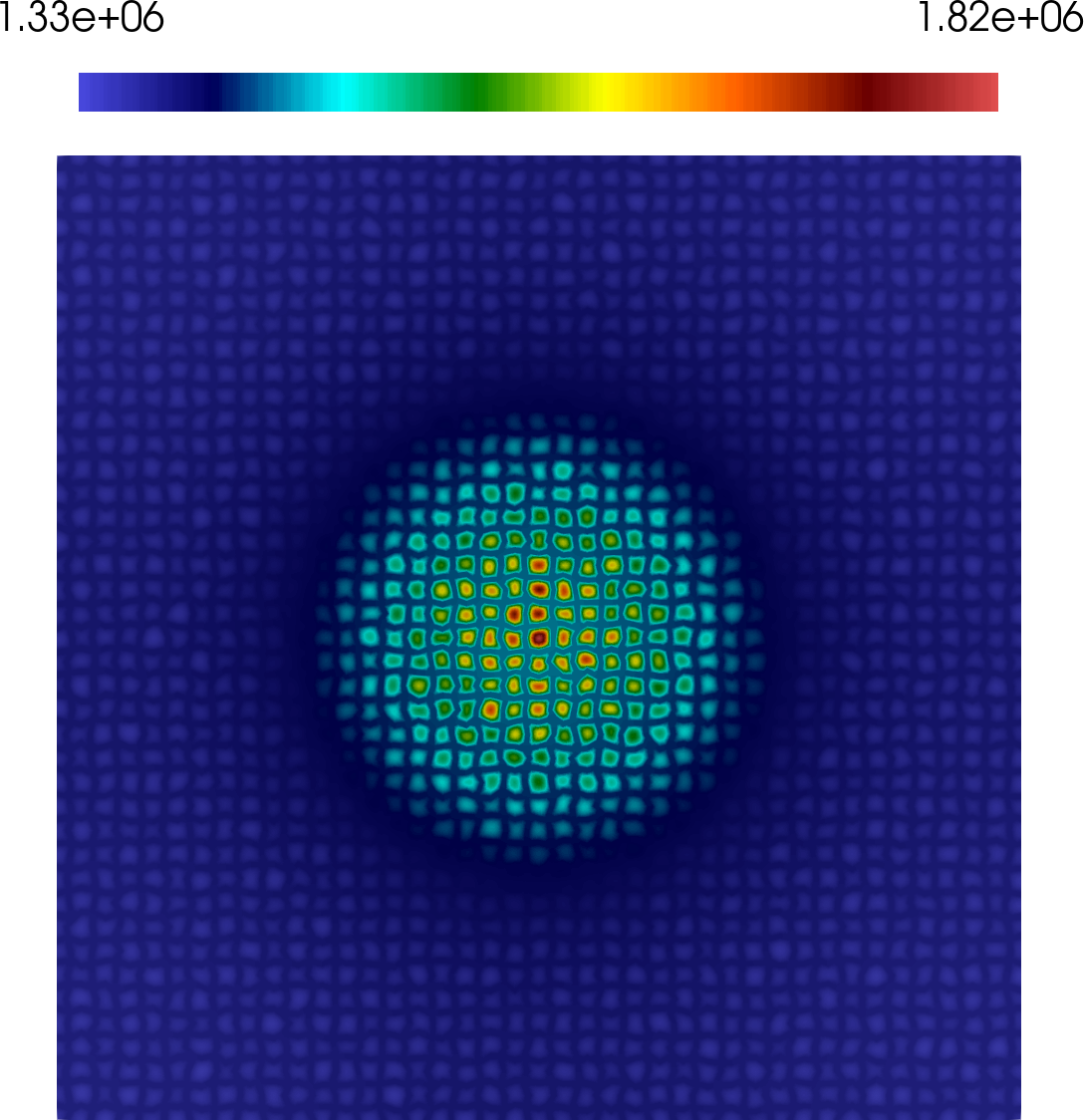}
\includegraphics[height=0.24\textheight]{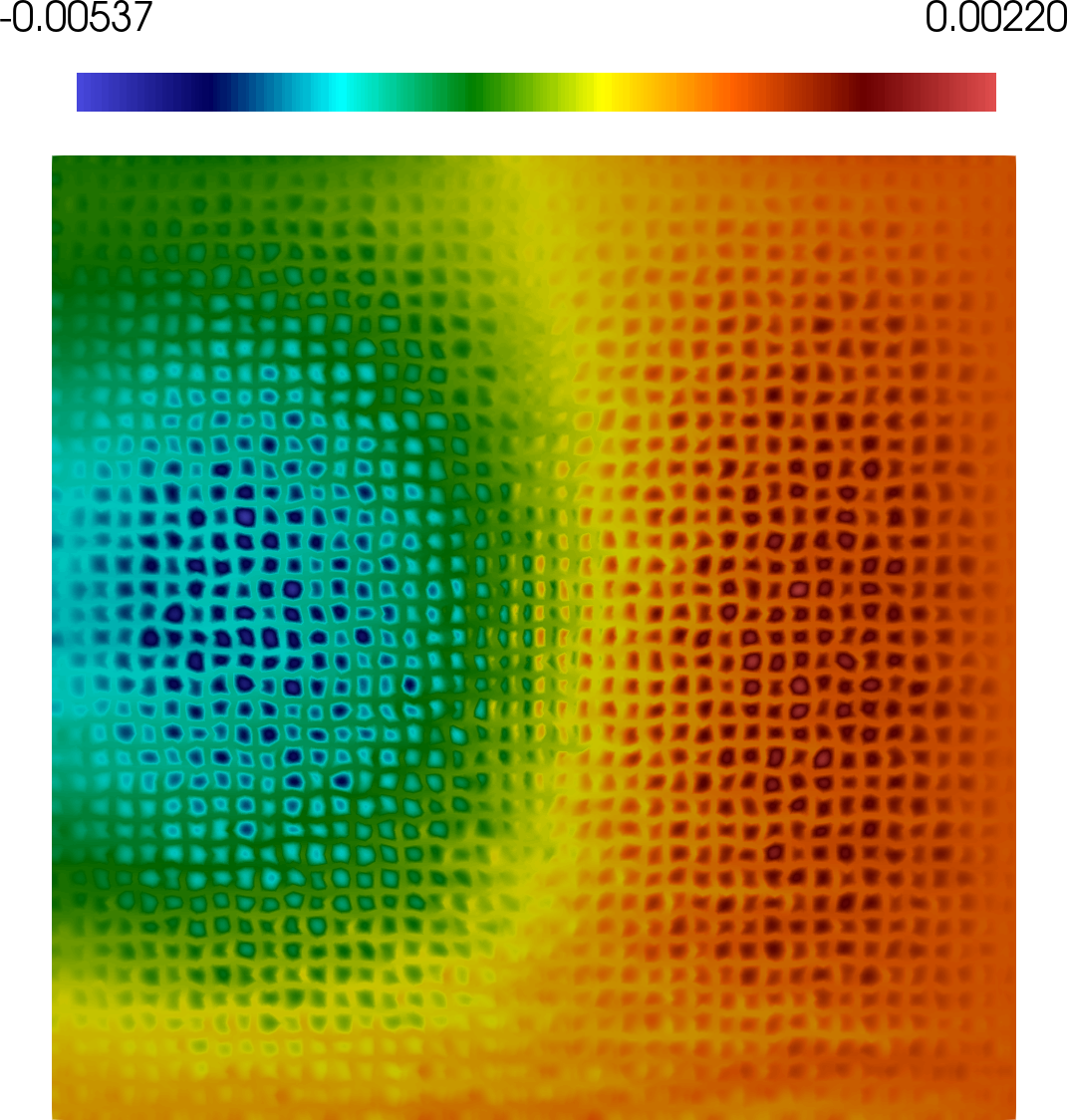}
\includegraphics[height=0.24\textheight]{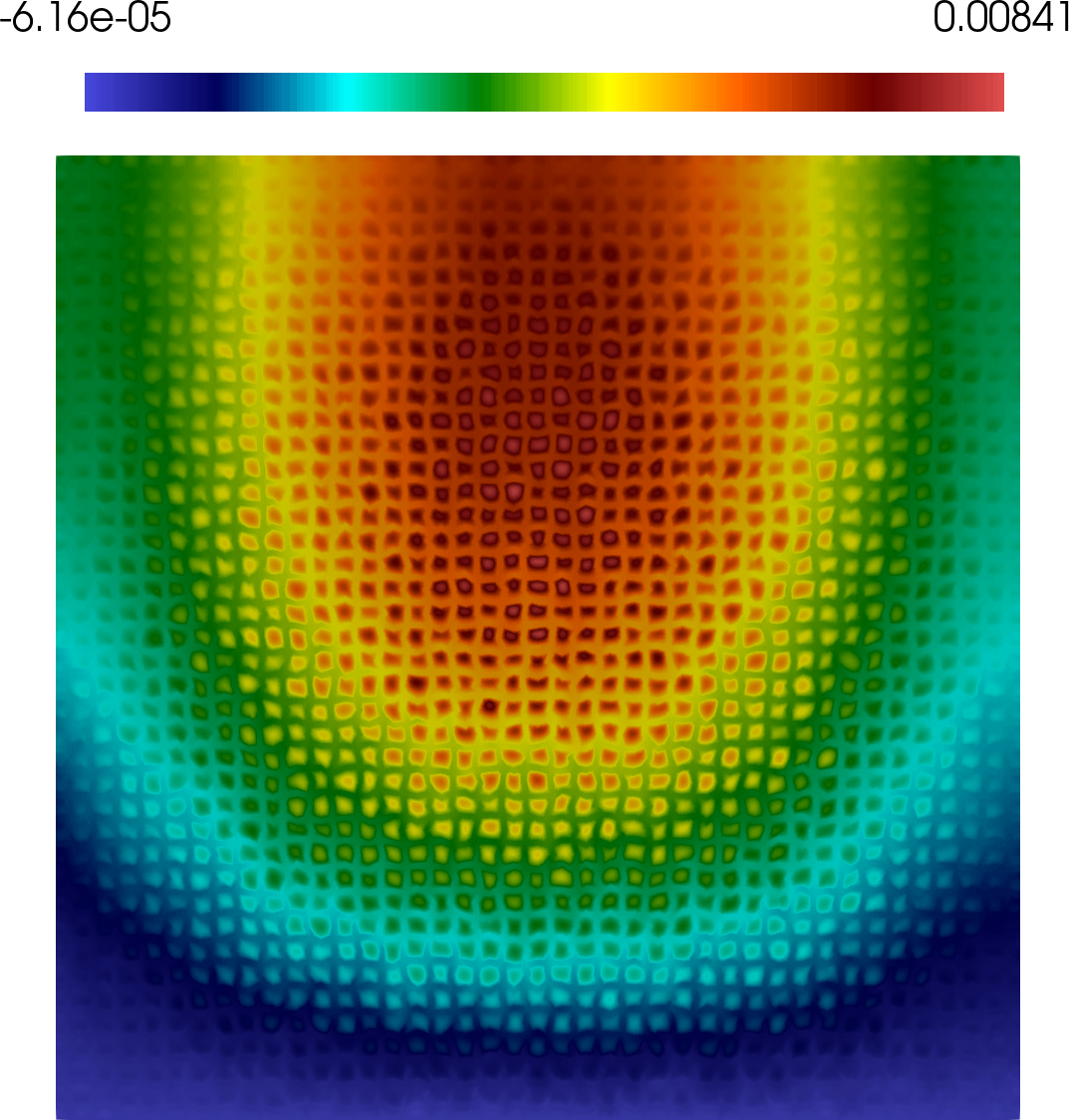}
\caption{Distributions of pressure and displacements in $x_1$ and $x_2$ directions (from left to right) at the final time on the fine grid for Example 3.}
\label{fig:fine_results_ex_3}
\end{figure}

In Figure \ref{fig:fine_results_ex_3}, we depict distributions of pressure and displacements in $x_1$ and $x_2$ directions (from left to right) at the final on the fine grid. One can see the influence of the non-periodic heterogeneous microstructure on the obtained numerical solutions. We observe that the fluid diffuses in an isotropic manner from the source in the middle of the domain. We see that displacement in the $x_1$ direction corresponds to the horizontal stretching. One can also observe upward stretching in the distribution of displacement in the $x_2$ direction.

\begin{figure}[hbt!]
\centering
\begin{subfigure}{\textwidth}
\centering
\includegraphics[height=0.24\textheight]{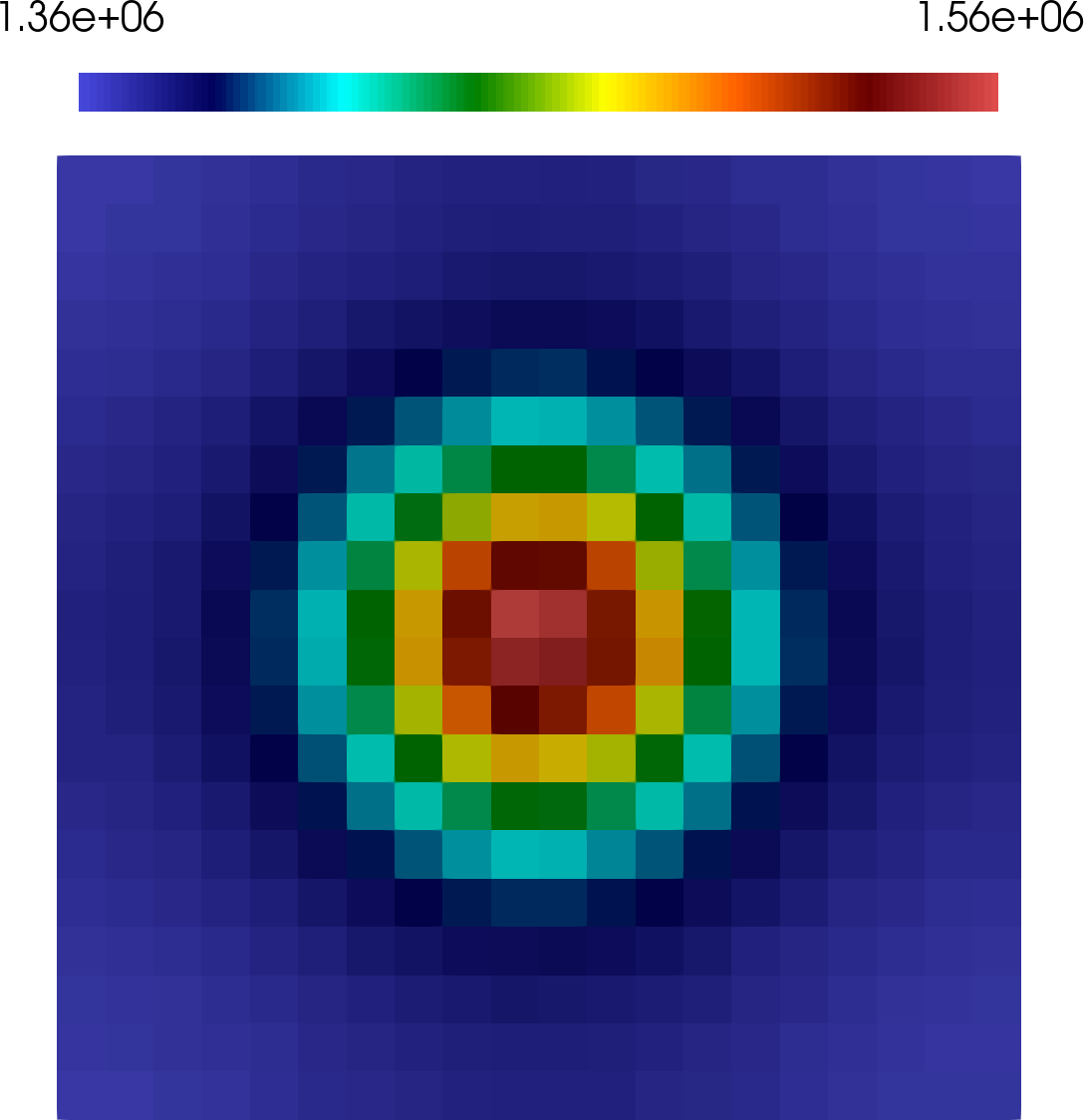}
\includegraphics[height=0.24\textheight]{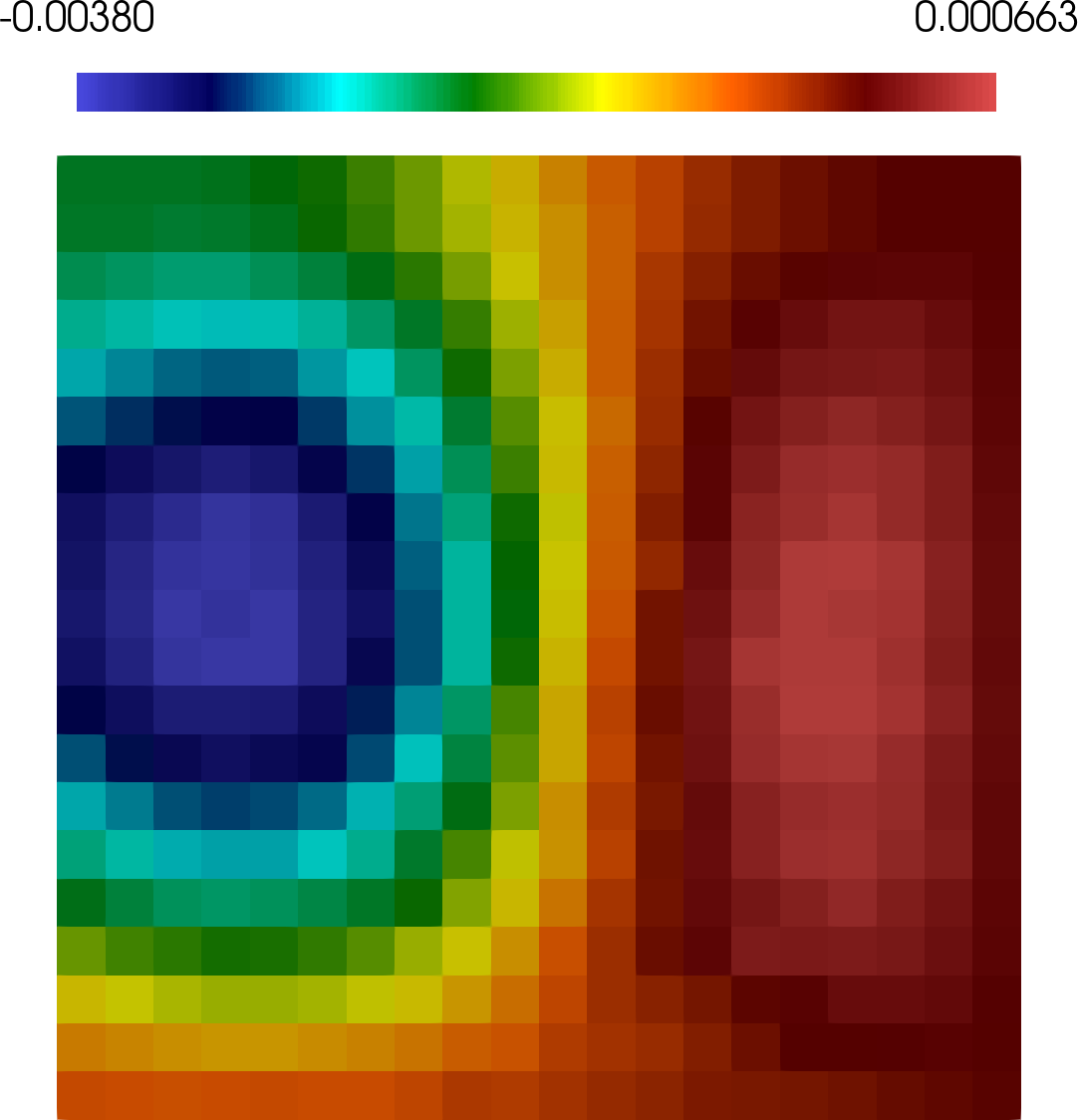}
\includegraphics[height=0.24\textheight]{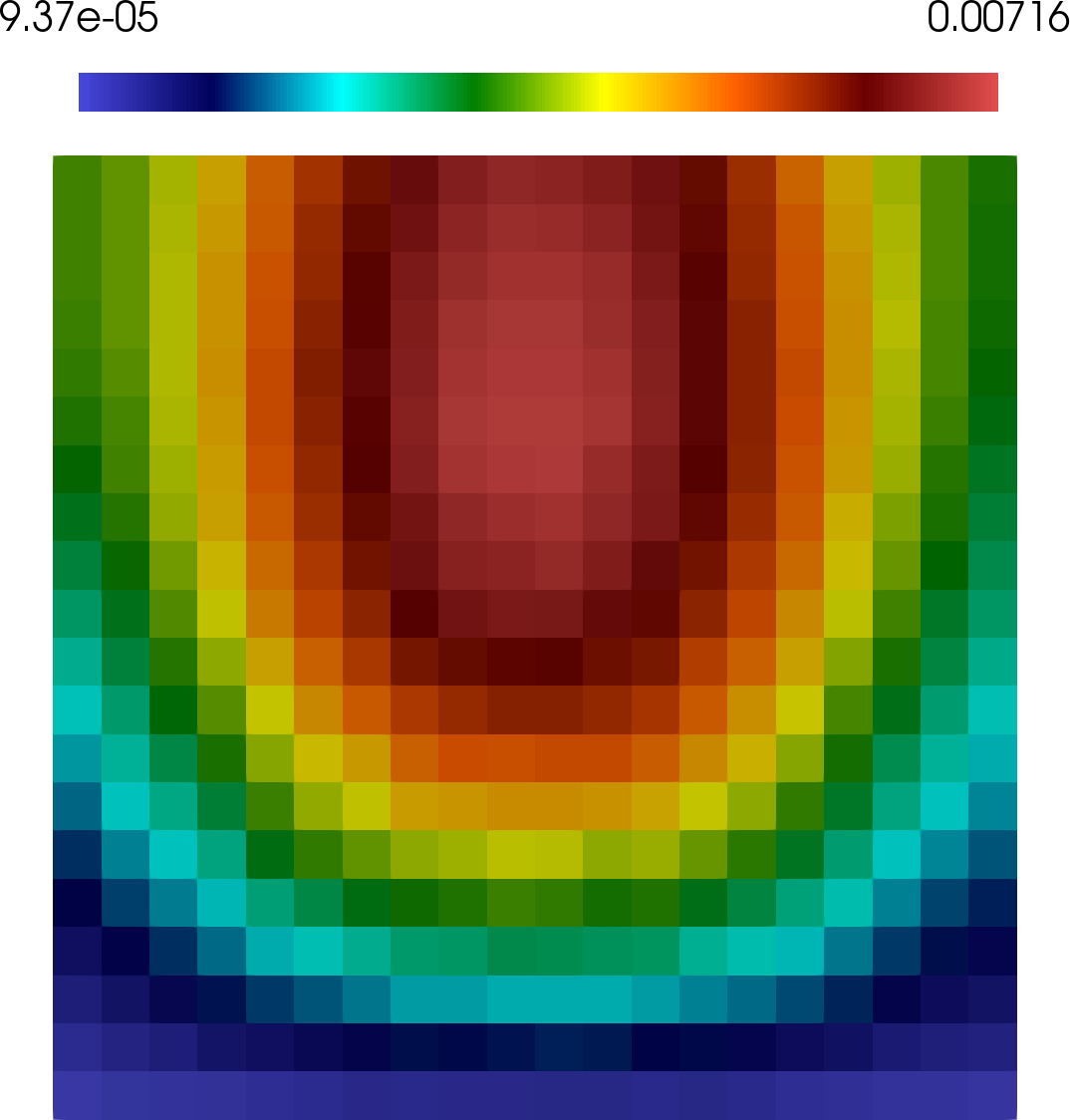}
\caption{Reference average solution}
\end{subfigure}
\begin{subfigure}{\textwidth}
\centering
\includegraphics[height=0.24\textheight]{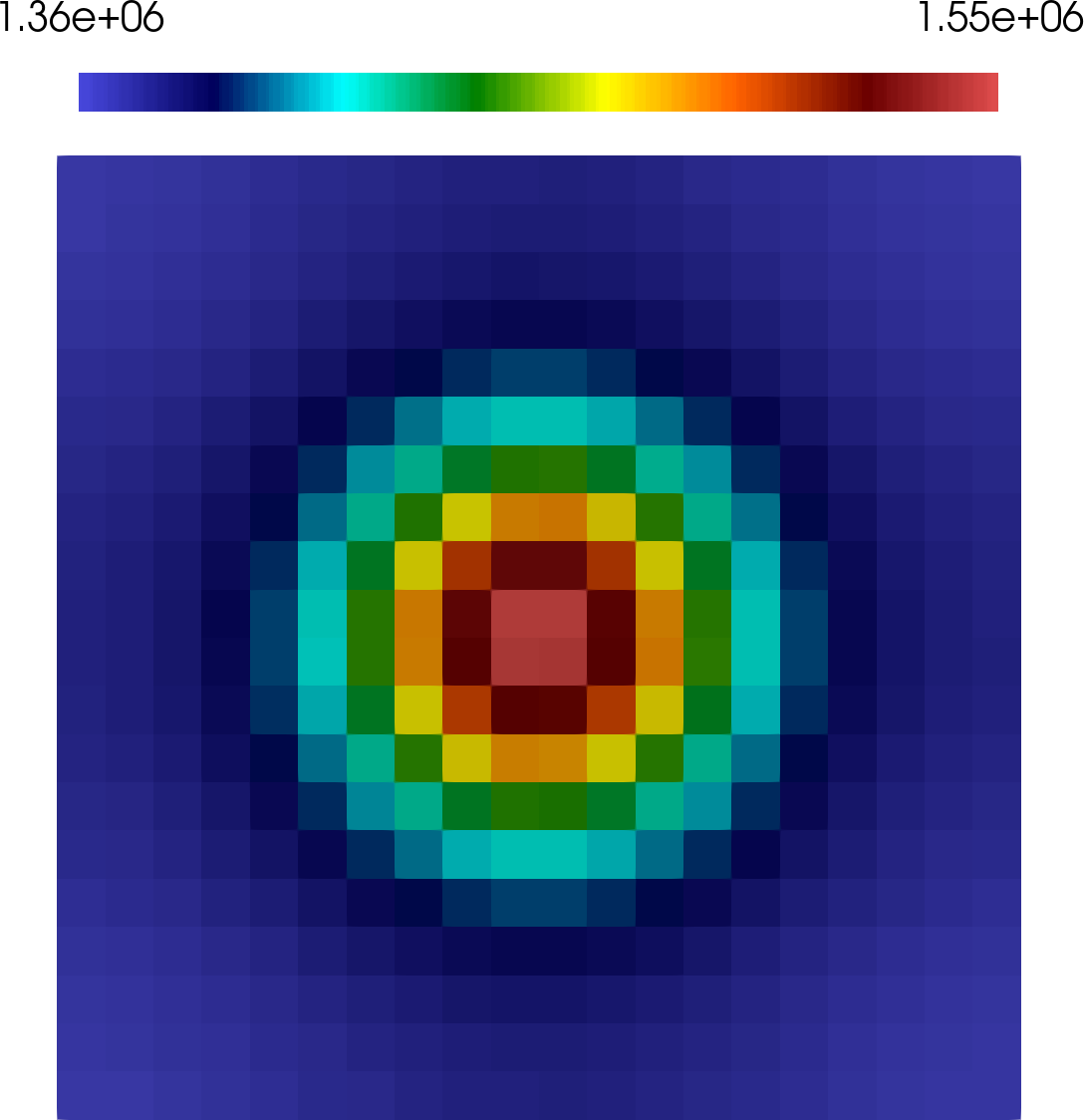}
\includegraphics[height=0.24\textheight]{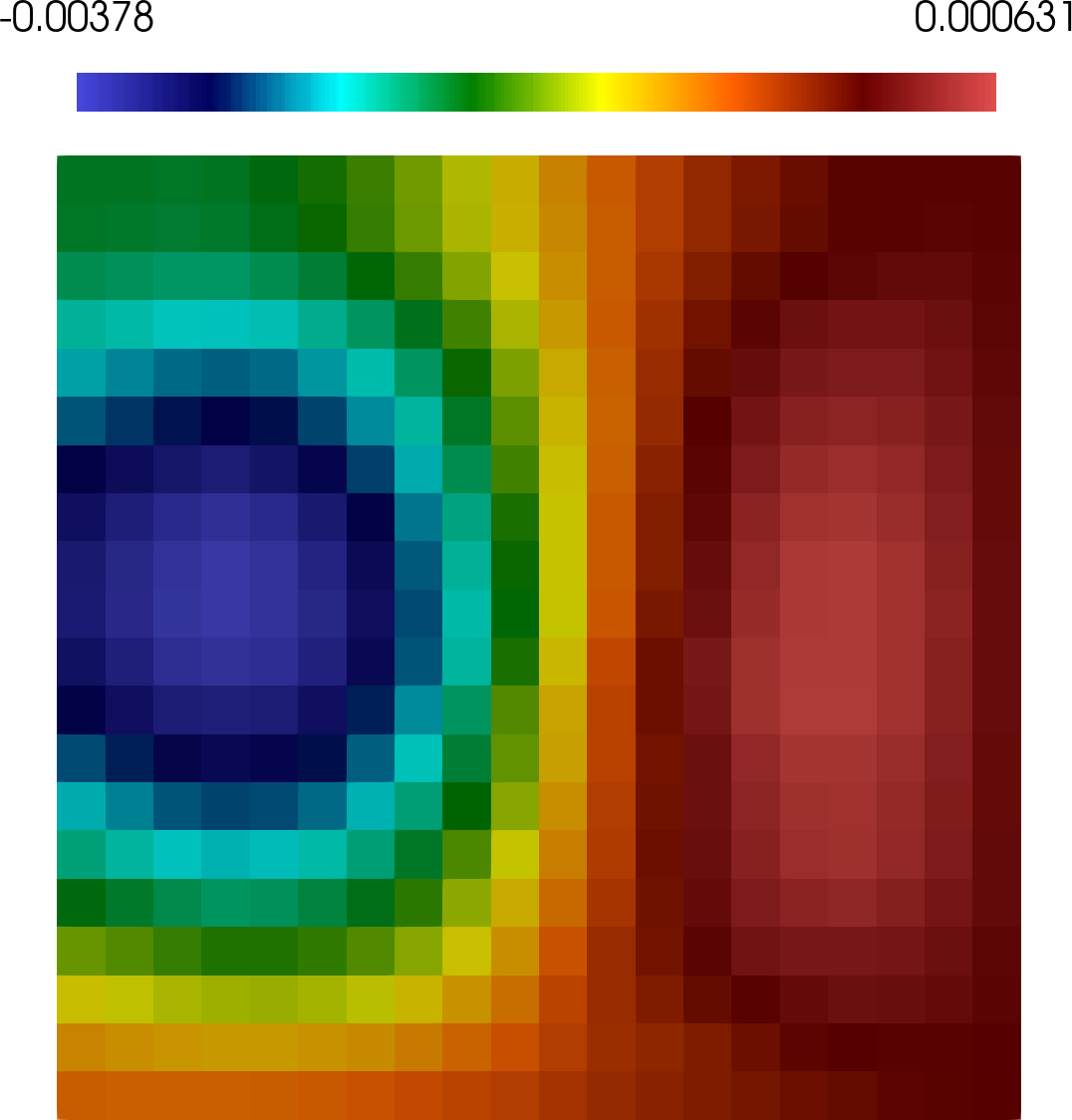}
\includegraphics[height=0.24\textheight]{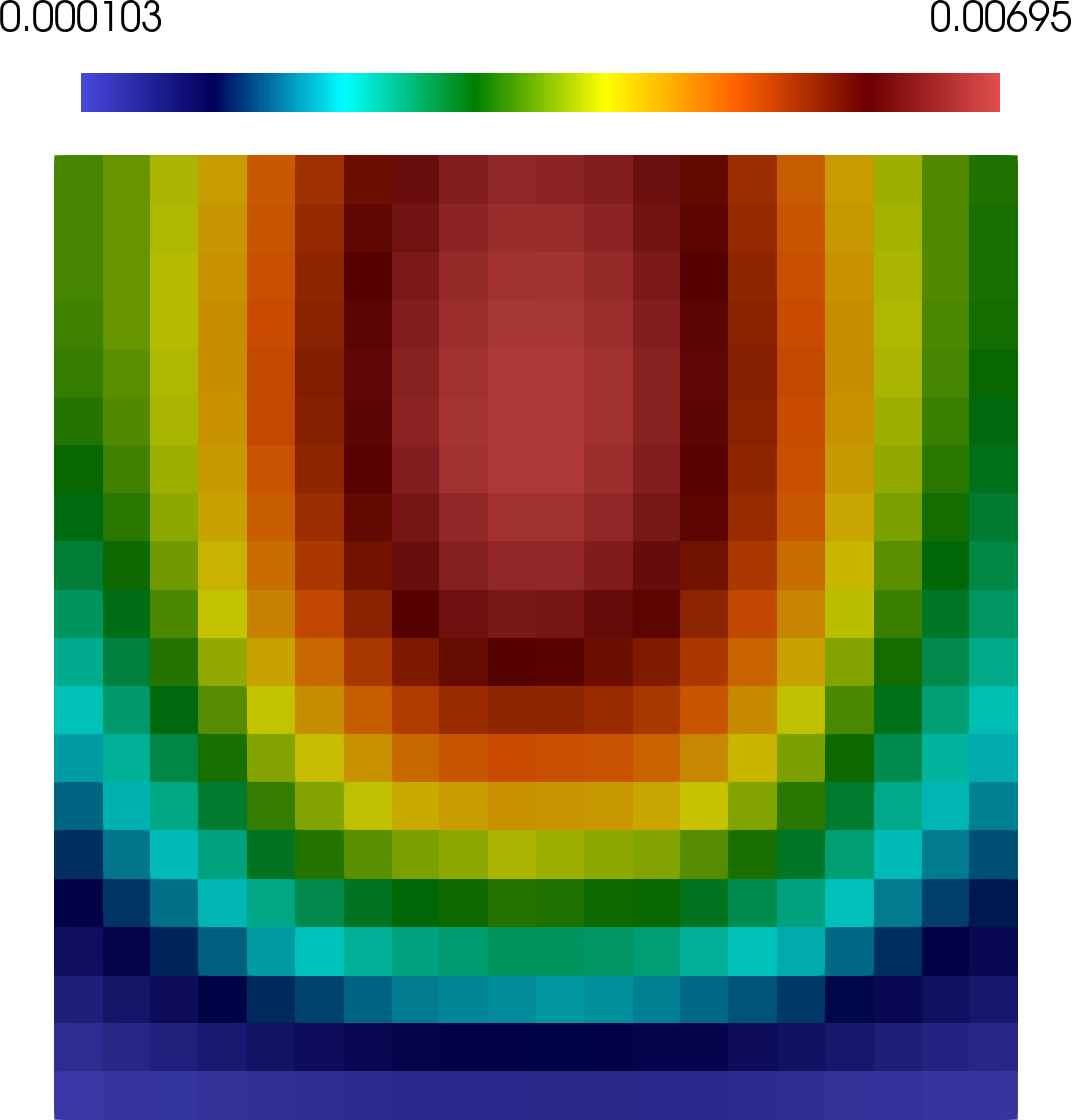}
\caption{Homogenized average solution using the full multicontinuum model}
\end{subfigure}
\begin{subfigure}{\textwidth}
\centering
\includegraphics[height=0.24\textheight]{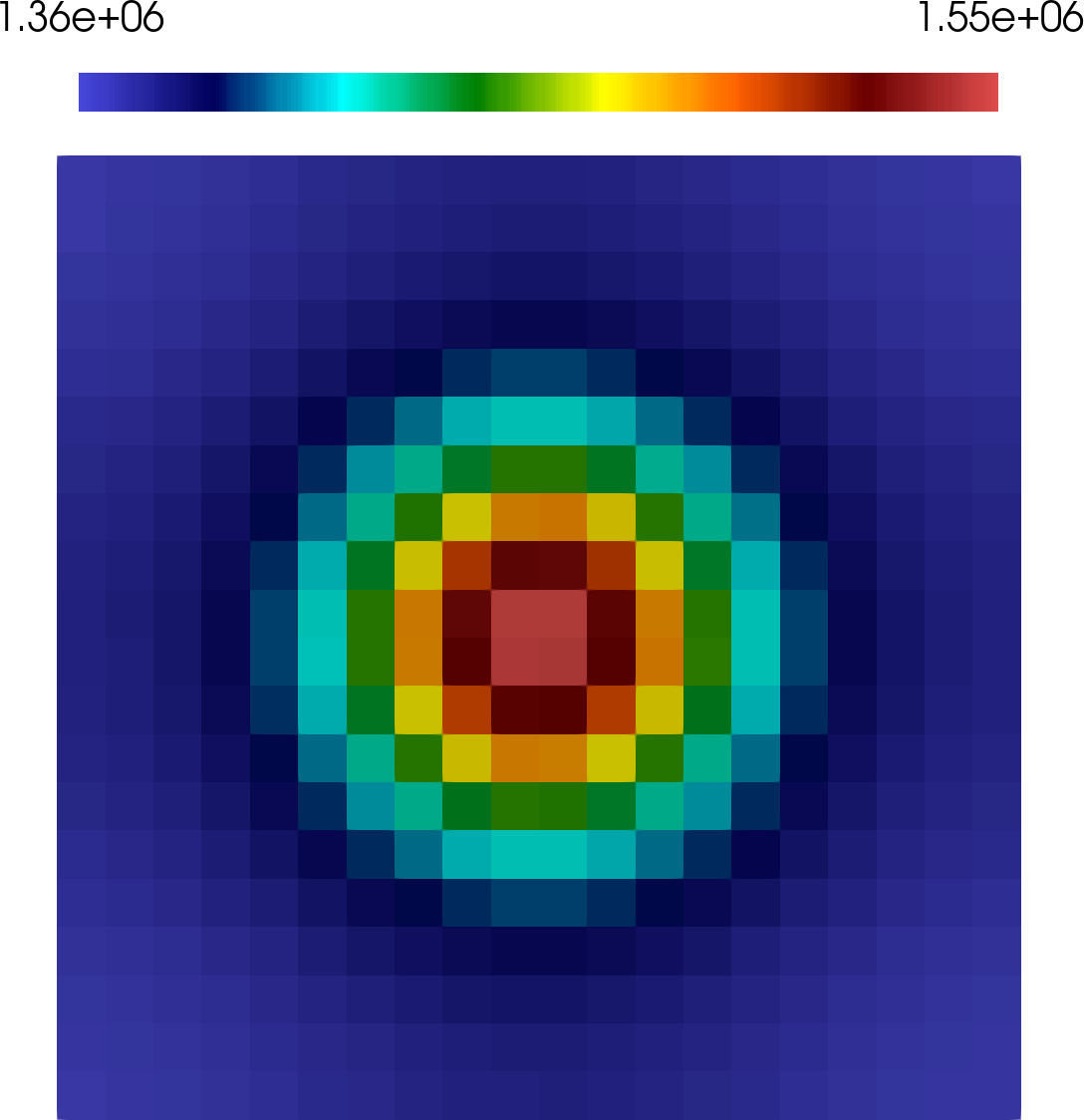}
\includegraphics[height=0.24\textheight]{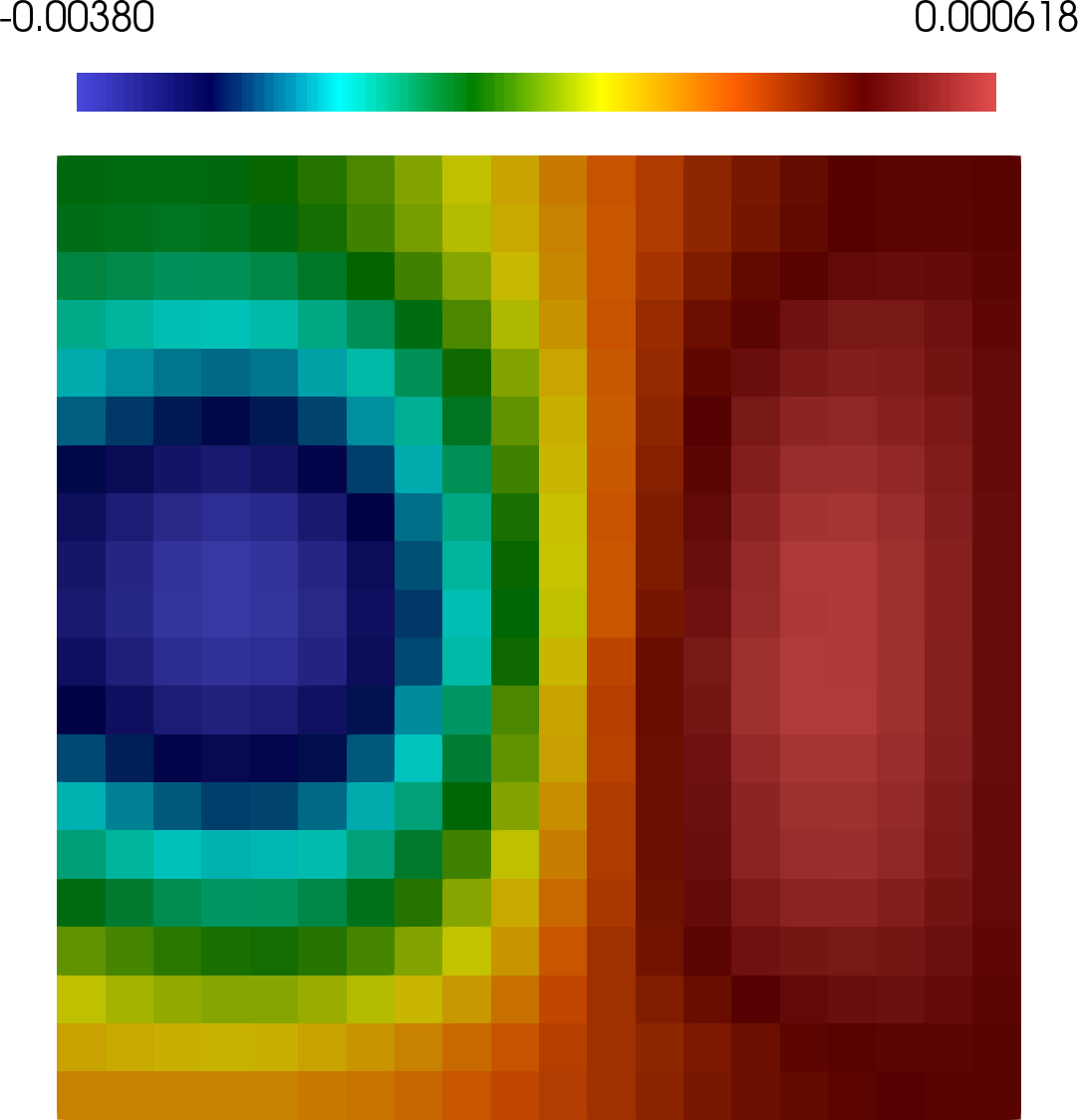}
\includegraphics[height=0.24\textheight]{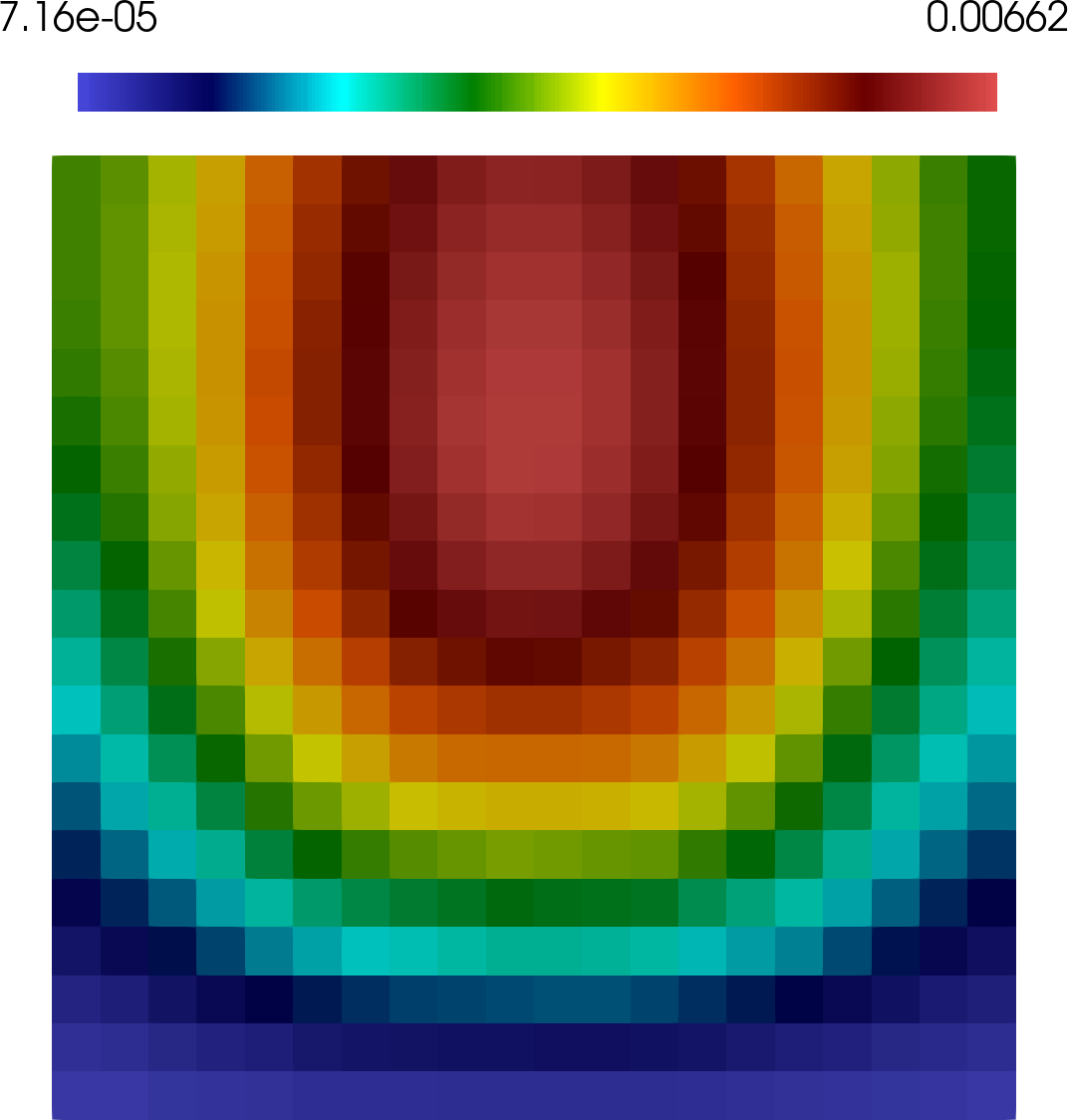}
\caption{Homogenized average solution using the simplified multicontinuum model 2}
\end{subfigure}
\caption{Distributions of average pressure and displacements in $x_1$ and $x_2$ directions (from left to right) in $\Omega_1$ at the final time on the coarse grid $20 \times 20$ for Example 3.}
\label{fig:coarse_results_1_ex_3}
\end{figure}

\begin{figure}[hbt!]
\centering
\begin{subfigure}{\textwidth}
\centering
\includegraphics[height=0.24\textheight]{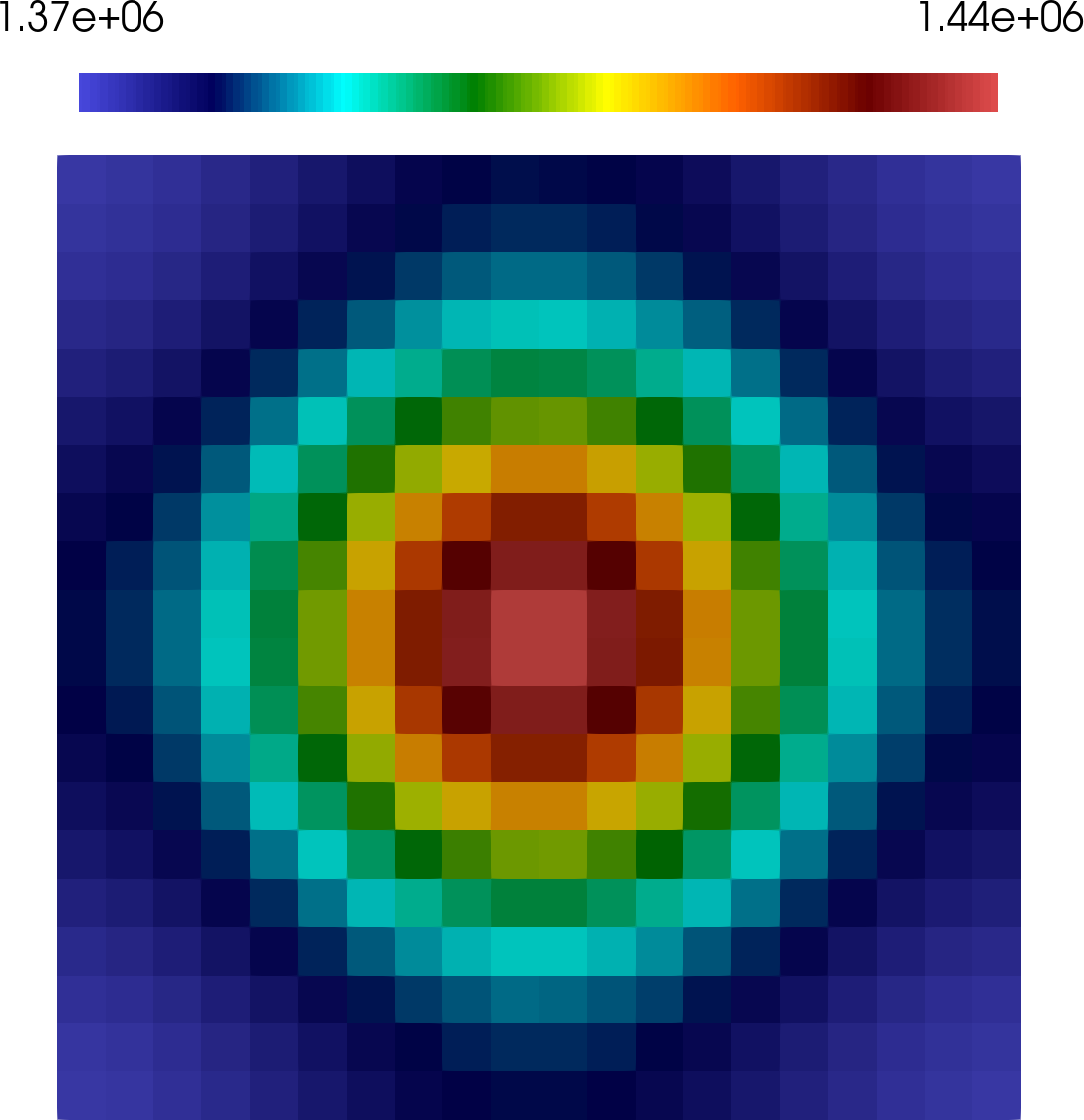}
\includegraphics[height=0.24\textheight]{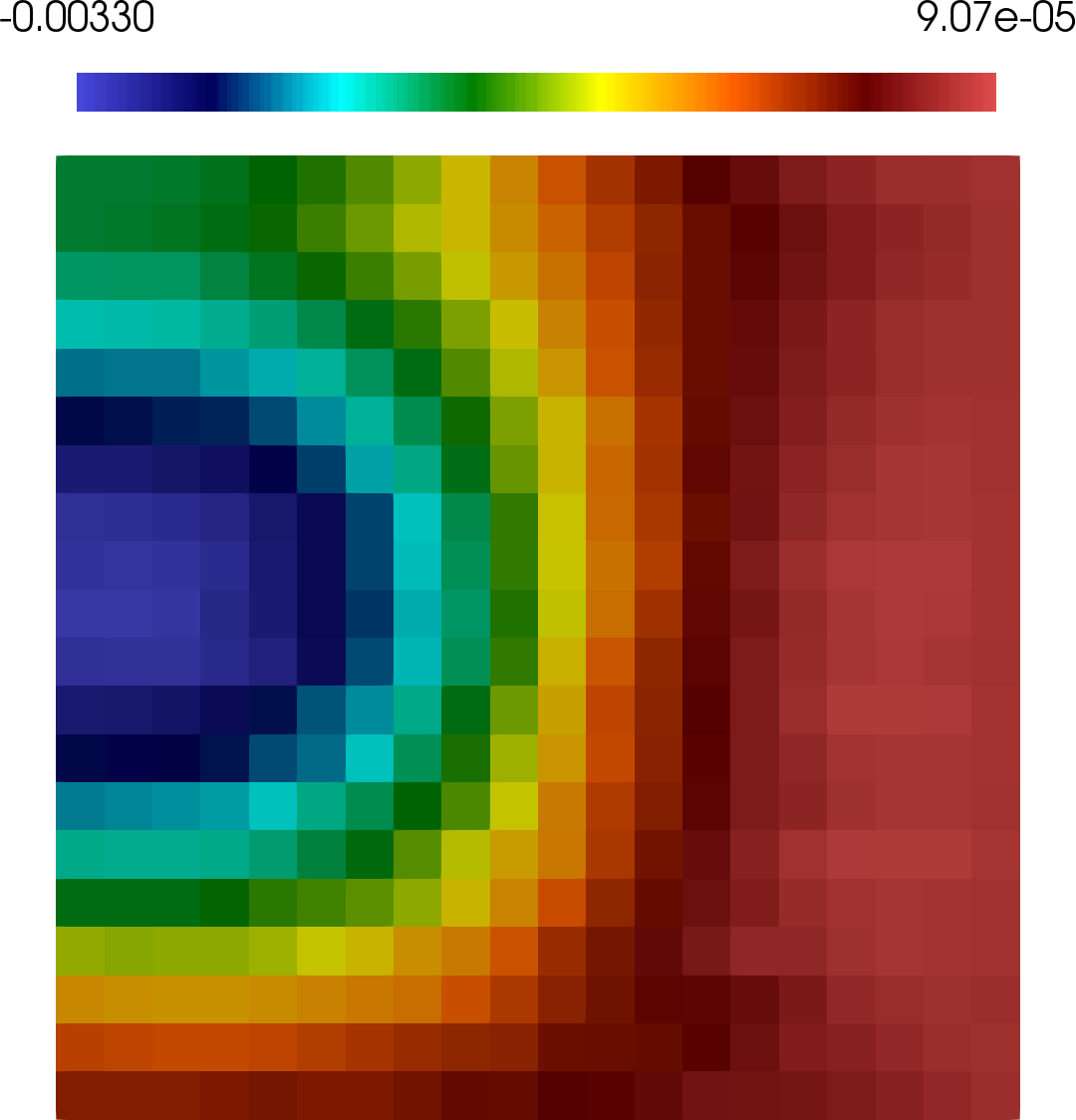}
\includegraphics[height=0.24\textheight]{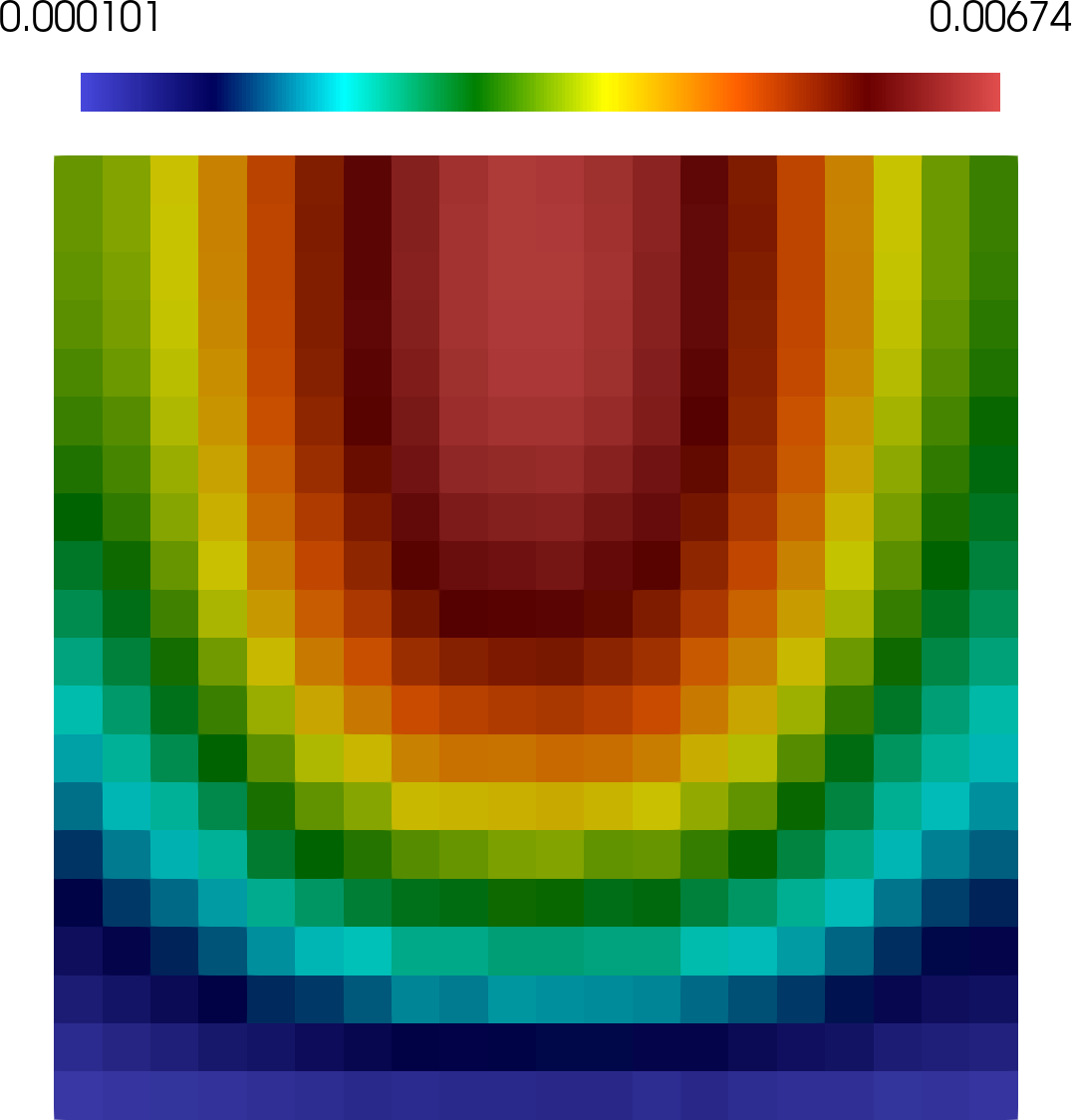}
\caption{Reference average solution}
\end{subfigure}
\begin{subfigure}{\textwidth}
\centering
\includegraphics[height=0.24\textheight]{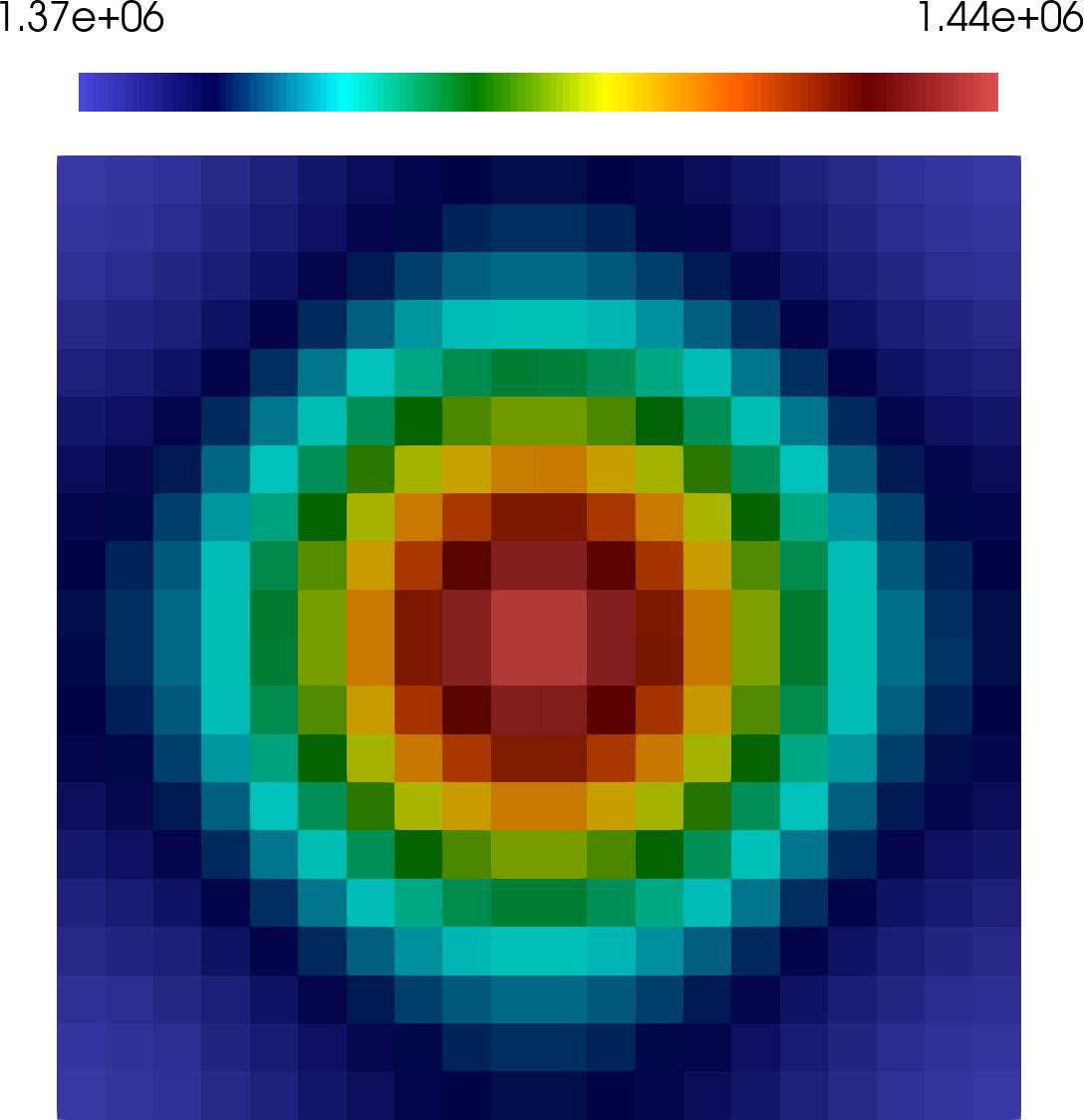}
\includegraphics[height=0.24\textheight]{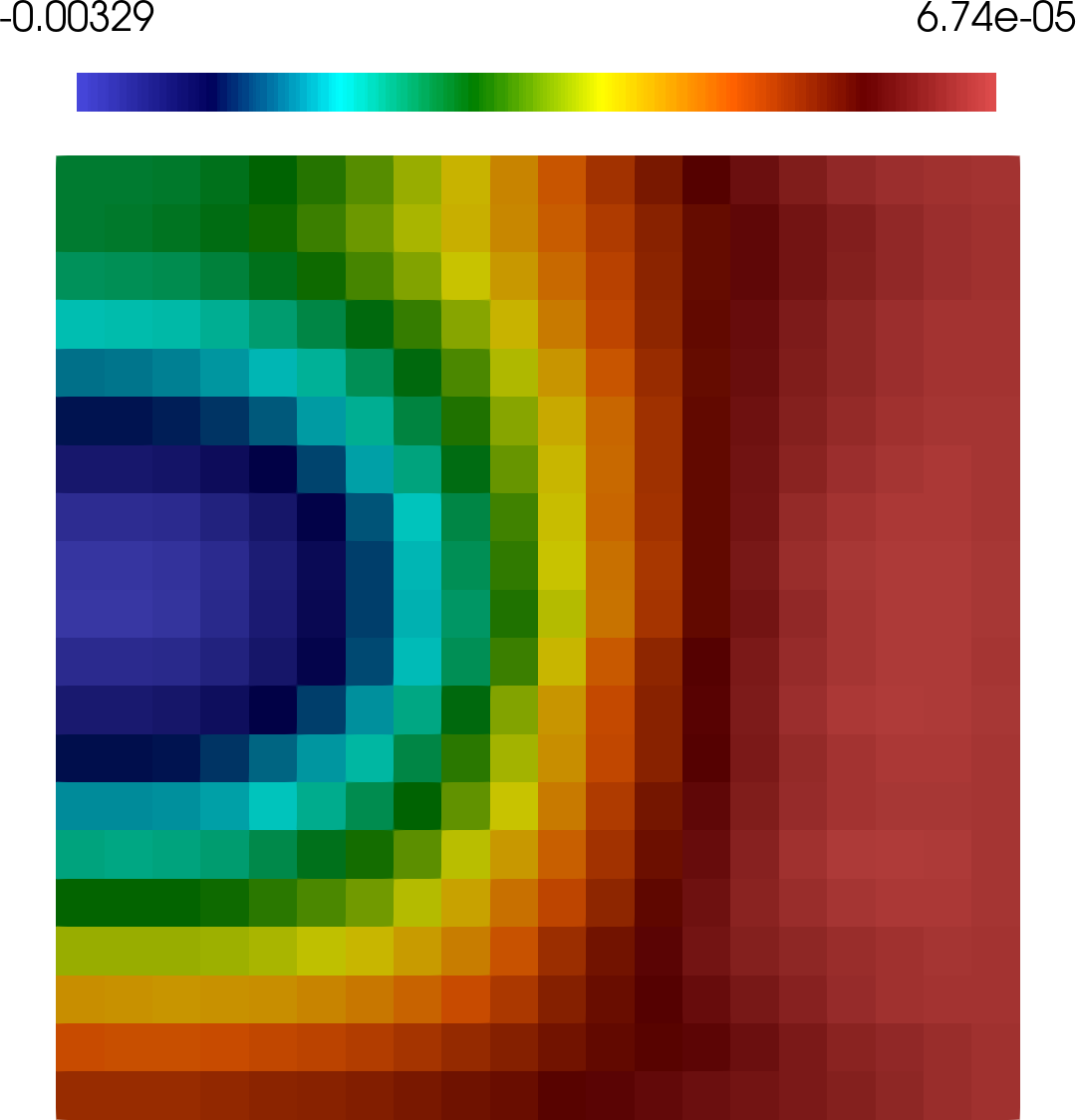}
\includegraphics[height=0.24\textheight]{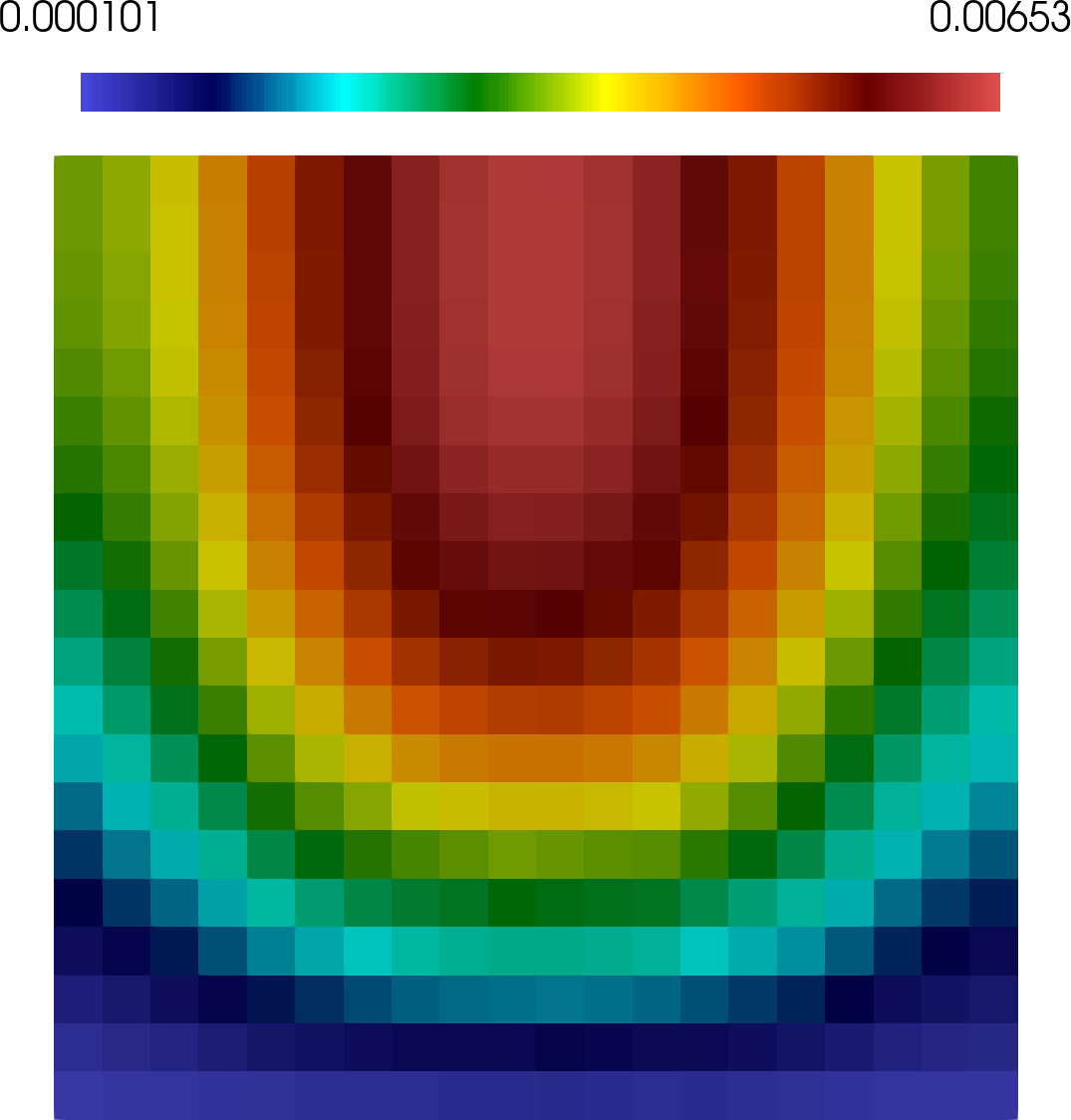}
\caption{Homogenized average solution using the full multicontinuum model}
\end{subfigure}
\begin{subfigure}{\textwidth}
\centering
\includegraphics[height=0.24\textheight]{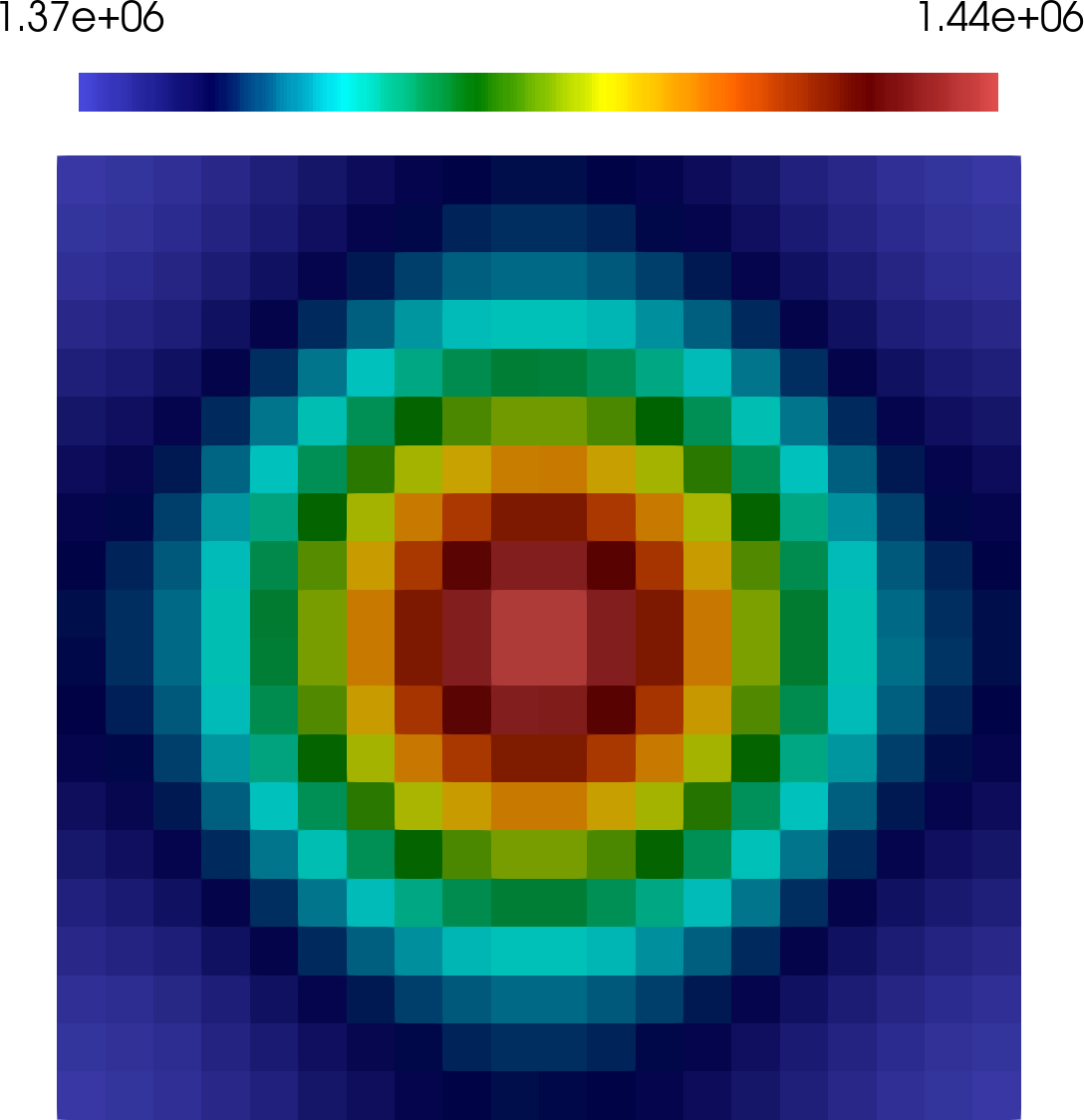}
\includegraphics[height=0.24\textheight]{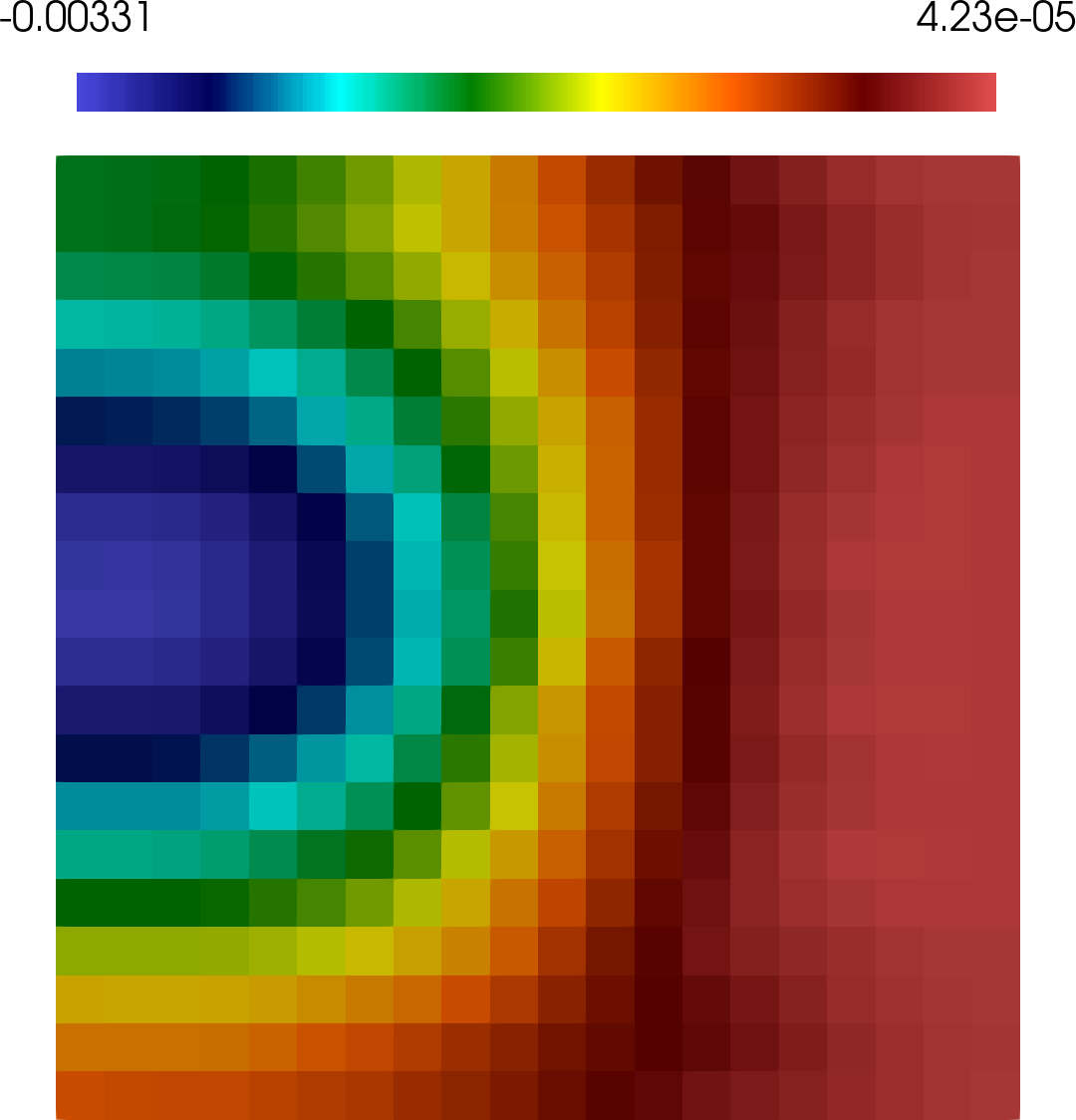}
\includegraphics[height=0.24\textheight]{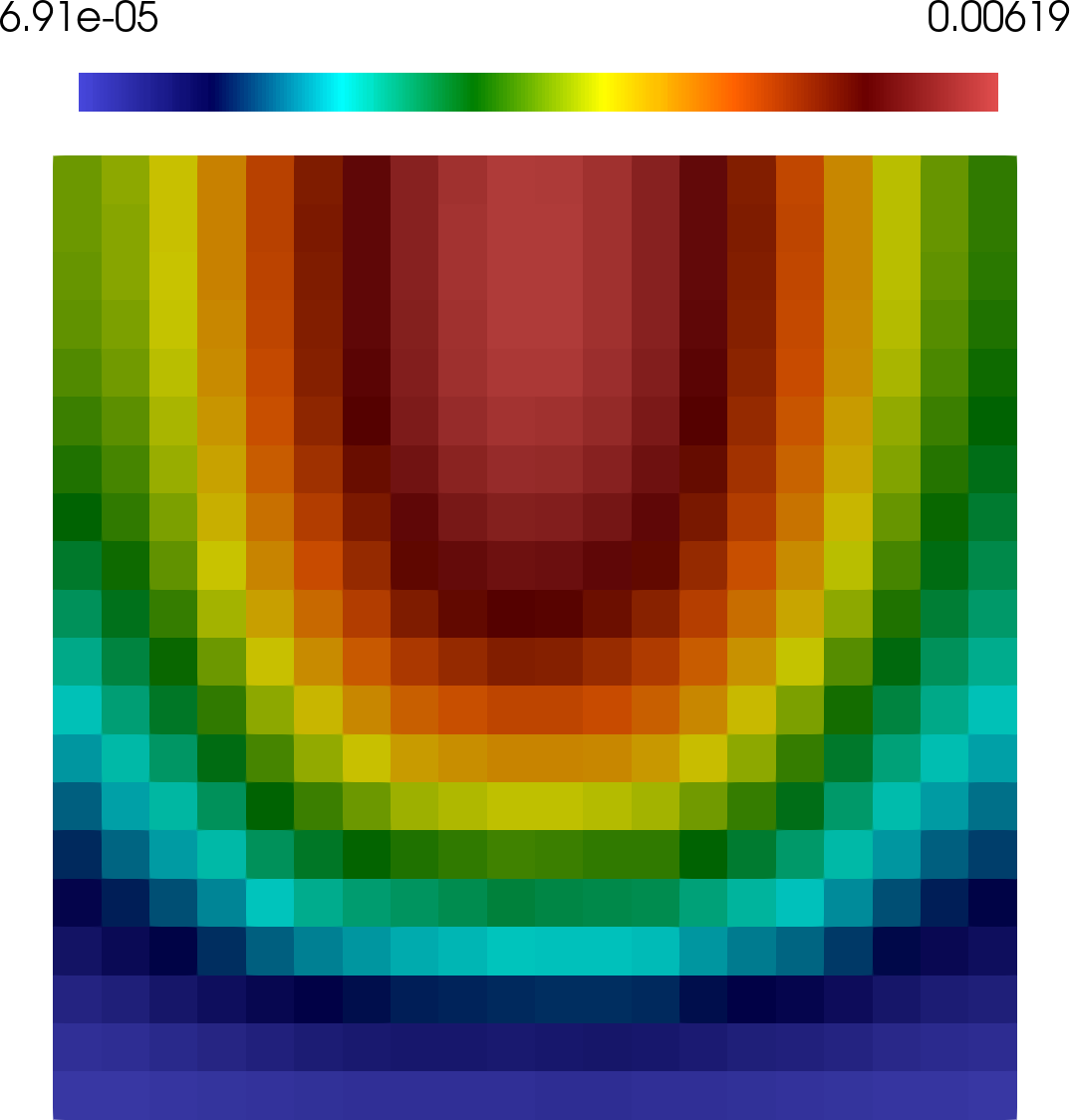}
\caption{Homogenized average solution using the simplified multicontinuum model 2}
\end{subfigure}
\caption{Distributions of average pressure and displacements in $x_1$ and $x_2$ directions (from left to right) in $\Omega_2$ at the final time on the coarse grid $20 \times 20$ for Example 3.}
\label{fig:coarse_results_2_ex_3}
\end{figure}

Figures \ref{fig:coarse_results_1_ex_3}-\ref{fig:coarse_results_2_ex_3} depict average pressure and displacements in the $x_1$ and $x_2$ directions (from left to right) at the final time in $\Omega_1$ and $\Omega_2$, respectively, using the coarse grid $20 \times 20$. As in the previous examples, we depict the reference and homogenized solutions using the full and second simplified models from top to bottom. One can see that all the average solutions are very similar. However, one can notice that the value ranges of the displacements in the $x_2$ direction of the second simplified model are noticeably different from the reference ones. The second simplified model underestimates the displacements in $x_2$ direction.

In terms of the simulated process, we observe fluid diffusion from the source in the middle of the domain. As expected, the diffusion is faster in $\Omega_2$ due to high permeability. Displacements in the $x_1$ direction correspond to the horizontal stretching, where we have larger displacement values in $\Omega_1$ due to softer elastic properties. Displacements in the $x_2$ direction represent the vertical stretching with larger values in $\Omega_1$. The obtained results correspond to the fine-scale solutions in Figure \ref{fig:fine_results_ex_3}.

Let us now consider the errors of our multicontinuum approaches. Table \ref{tab:errors_ex_3} presents the relative $L_2$ errors at the final time using different multicontinuum models and coarse grids. The full and first simplified multicontinuum models provide high accuracy for all the solution fields and coarse grids. The second simplified model has more significant errors for displacements. For example, on the coarse grid $10 \times 10$, the second simplified model has $16.8 \%$ errors for the first continuum displacements and $18.2\%$ errors for the second continuum displacements. In contrast, the first simplified model has $4.19\%$ errors for the first continuum displacements and $4.4\%$ errors for the second continuum displacements.

\begin{table}[hbt!]
\centering
\begin{tabular}{c|c|c|c|c}
Multicontinuum model & $e^{(1)}_p$ & $e^{(2)}_p$ & $e^{(1)}_u$ & $e^{(2)}_u$\\ \hline
\multicolumn{5}{c}{Coarse grid $10 \times 10$} \\ \hline
Full         & 2.40e-03 & 7.65e-04 & 4.17e-02 & 4.39e-02 \\
Simplified 1 & 2.40e-03 & 7.64e-04 & 4.19e-02 & 4.40e-02 \\
Simplified 2 & 2.43e-03 & 8.65e-04 & 1.68e-01 & 1.82e-01 \\ \hline
\multicolumn{5}{c}{Coarse grid $20 \times 20$} \\ \hline
Full         & 9.50e-04 & 2.25e-04 & 3.07e-02 & 3.28e-02 \\
Simplified 1 & 9.89e-04 & 2.51e-04 & 3.11e-02 & 3.30e-02 \\
Simplified 2 & 1.02e-03 & 3.53e-04 & 9.14e-02 & 9.84e-02 \\
\end{tabular}
\caption{Relative $L_2$ errors for different multicontinuum models and coarse grids. Example 3.}
\label{tab:errors_ex_3}
\end{table}

In this way, our proposed multicontinuum approaches can also provide accurate solutions for non-periodic microstructures.

\FloatBarrier
\section{Conclusion}\label{sec:conclusion}

In this work, we have considered poroelasticity problems in high-contrast heterogeneous media. To handle multiscale coefficients, we have proposed a multicontinuum homogenization approach. We formulated coupled cell problems in the oversampled RVEs with constraints to consider different homogenized effects. By solving these problems, we obtained multicontinuum expansions of the fine-scale solution fields over macroscopic variables. Then, we rigorously derived the general multicontinuum poroelastic model for an arbitrary number of continua. In addition, we presented two simplified multicontinuum models. Note that the standard multiple-network poroelasticity model can be considered as a particular case of our second simplified multicontinuum model. We presented representative numerical examples for different microstructures. The numerical results show that our proposed multicontinuum approach can approximate the reference average solutions with high accuracy. However, the second simplified model has more significant errors and does not always accurately capture the reference solution's features, especially in complex heterogeneous media.

\section*{Acknowledgement}\label{sec:acknowledgement}

The authors gratefully acknowledge the support of ADNOC Grant No. 8434000476. The authors are grateful to Dr. Wing Tat Leung for valuable comments and discussions.

\bibliographystyle{unsrt}
\bibliography{lit}

\end{document}